\documentclass{article}
\usepackage[T1]{fontenc}
\usepackage{amsthm, fullpage}
\usepackage{amsmath}
\usepackage{amssymb}
\usepackage{amsfonts}
\usepackage{eucal}
\usepackage[all]{xy}
\usepackage[frenchb, english]{babel}
\usepackage{pslatex}
\usepackage[pagebackref]{hyperref}

\newtheorem{prop}{Proposition}[subsection]
\newtheorem{theo}[prop]{Théor\`eme}
\newtheorem*{theo**}{Théorème}
\newtheorem{coro}[prop]{Corollaire}
\newtheorem*{conj*}{Conjecture}
\newtheorem{lemm}[prop]{Lemme}
\newtheorem{lemm*}{Lemme}[prop]

\theoremstyle{definition}

\newtheorem{vide}[prop]{}
\newtheorem{defi}[prop]{Définition}
\newtheorem*{defi*}{Définition}

\theoremstyle{remark}
\newtheorem{rema}[prop]{Remarques}

\newtheorem{nota}[prop]{Notations}

\numberwithin{equation}{prop}

\newcommand{\riso}{ \overset{\sim}{\longrightarrow}\, }
\newcommand{\liso}{ \overset{\sim}{\longleftarrow}\, }

\newcommand{\Spec}{\mathrm{Spec}\,}
\newcommand{\Spf}{\mathrm{Spf}\,}

\renewcommand{\sp}{\mathrm{sp}}

\renewcommand{\AA}{{\mathcal{A}}}

\newcommand{\FF}{{\mathcal{F}}}
\newcommand{\B}{{\mathcal{B}}}

\newcommand{\E}{{\mathcal{E}}}

\newcommand{\M}{{\mathcal{M}}}
\newcommand{\NN}{{\mathcal{N}}}
\newcommand{\D}{{\mathcal{D}}}

\newcommand{\PP}{{\mathcal{P}}}

\renewcommand{\O}{{\mathcal{O}}}

\newcommand{\V}{\mathcal{V}}

\newcommand{\Y}{\mathcal{Y}}
\newcommand{\ZZ}{\mathcal{Z}}
\newcommand{\X}{\mathfrak{X}}

\newcommand{\DD}{\mathbb{D}}
\renewcommand{\L}{\mathbb{L}}
\newcommand{\R}{\mathbb{R}}
\newcommand{\Q}{\mathbb{Q}}
\newcommand{\Z}{\mathbb{Z}}
\newcommand{\N}{\mathbb{N}}

\newcommand{\hdag}{  \phantom{}{^{\dag} }    }

\begin{document}
\selectlanguage{frenchb}

\title{Sur la compatibilité à Frobenius de l'isomorphisme de dualité relative}
\author{Daniel Caro \footnote{L'auteur a bénéficié du soutien du réseau européen TMR \textit{Arithmetic Algebraic Geometry}
(contrat numéro UE MRTN-CT-2003-504917).}}


\maketitle

\begin{abstract}
\selectlanguage{english}
Let $\V$ be a mixed characteristic complete discrete valuation ring,
let $\X$ and $\Y$ be two smooth formal $\V$-schemes,
let $f _0$ : $X \rightarrow Y$ be a projective morphism between their special fibers,
let $T$ be a divisor of $Y$ such that $T _X := f_0 ^{-1} ( T) $ is a divisor of $X$
and let $\M \in D ^\mathrm{b} _\mathrm{coh} (\D ^\dag _{\X} (\hdag T _X) _{\Q})$.
We construct the relative duality isomorphism
$ f _{0T +} \circ \DD _{\X, T _X} (\M ) \riso \DD _{\Y, T }  \circ f _{0T +}  (\M )$.
This generalizes the known case when there exists a lifting $f\,:\,\X \to \Y$ of $f _{0}$.
Moreover, when $f _0$ is a closed immersion,
we prove that this isomorphism commutes
with Frobenius.
\end{abstract}

\selectlanguage{frenchb}
\date
\tableofcontents

\section*{Introduction}

Soient $\V$ un anneau de valuation discrète complet d'inégales caractéristiques $(0,\,p)$,
$\X$ et $\Y$ deux $\V$-schémas formels lisses,
$X$ et $Y$ les fibres spéciales correspondantes,
$u _0$ : $X \rightarrow Y$ un morphisme propre,
$T$ un diviseur de $Y$ tel que $T _X := u _0 ^{-1} (T)$ soit un diviseur de $X$.
Les modules sont par défaut à droite.
Dans cette introduction, on désigne par $F$ la puissance $s$-ème du Frobenius absolu de $X$ ou de $Y$,
avec $s$ un entier fixé.

On désigne par $\D ^\dag _{\X} (\hdag T _X) _{\Q}$, l'{\it anneau des opérateurs différentiels de niveau fini
à singularités surconvergentes le long de $T _X$} (voir \cite[4.2.5]{Be1}).
On note $D ^\mathrm{b} _\mathrm{coh} (\D ^\dag _{\X} (\hdag T _X) _{\Q} \overset{^d}{})$, la catégorie dérivée des complexes de $\D ^\dag _{\X} (\hdag T _X) _{\Q}$-modules à droite à cohomologie cohérente
et bornée, l'indice $d$ précisant qu'il s'agit de modules à droite. De même en remplaçant
$\X$ par $\Y$ et $T_X$ par $T$.
On note $\DD _{\Y ,T }$ (ou $\DD _{T}$) le foncteur dual
$\D ^\dag _{\Y} (\hdag T ) _{\Q}$-linéaire,
$u _{0T+}$ (resp. $u ^! _{0 T}$) l'image directe (resp. l'image inverse extraordinaire)
par $u_0$ à singularités surconvergentes le long de $T$.
Lorsque $u_0$ se relève en un morphisme $u$ : $\X \rightarrow \Y$ de $\V$-schémas formels,
on les désigne respectivement par $u _{T+}$ et $u ^! _T$.

Supposons dans ce paragraphe que $u_0$ se relève en un morphisme $u$ : $\X \rightarrow \Y$ de $\V$-schémas formels lisses.
Pour tout $\M \in D ^\mathrm{b} _\mathrm{coh} (\D ^\dag _{\X} (\hdag T _X) _{\Q}\overset{^d}{})$,  
on dispose de l'isomorphisme de dualité relative
$\chi$ : $u _{T+} \circ \DD _{\X, T _X} (\M) \riso \DD _{\Y ,T } \circ u _{T+} (\M) $
fonctoriel en $\M$ (voir \cite[1.2.7]{caro_courbe-nouveau}).
Lorsque le diviseur $T$ est vide,
on retrouve par construction l'isomorphisme de Virrion de \cite{Vir04} (en effet, on vérifie dans ce cas que les catégories notées $\widetilde{D} _{\mathrm{parf}}$ et utilisées dans \cite{Vir04} correspondent aux complexes à cohomologie bornée et cohérente).
Pour tout complexe $\NN \in D ^\mathrm{b} _\mathrm{coh} (\D ^\dag _{\Y} (\hdag T ) _{\Q})$,
l'isomorphisme de dualité relative induit l'isomorphisme d'adjonction
$\mathrm{Hom} _{\D ^\dag _{\Y } (\hdag T ) _{\Q}} ( u _{T+}  (\M ), \NN )
\riso
\mathrm{Hom} _{\D ^\dag _{\X } (\hdag T _{X}) _{\Q}} ( \M , u _{T} ^{!} \NN )$.
Il en découle des morphismes d'adjonction entre image directe et image inverse extraordinaire.
Ceux-ci constituent
un ingrédient dans la preuve de la cohérence différentielle des $F$-isocristaux surconvergents unités
(\cite{caro_unite}).
Plus précisément, ceux-ci permettent de faire de la {\it descente propre, génériquement finie et étale}
au niveau des schémas.
De telles descentes propres génériquement finies et étales 
apparaissent si on utilise par exemple le théorème de désingularisation de de Jong (\cite{dejong}),
le théorème de monodromie génériquement finie et étale de Tsuzuki
pour les $F$-isocristaux surconvergents unités (\cite[1.3.1]{tsumono}) ou encore plus généralement le théorème de la réduction semistable de Kedlaya (voir \cite{kedlaya-semistableI,kedlaya-semistableII,kedlaya-semistableIII,kedlaya-semistableIV}).

Donnons à présent deux manières d'étendre cet isomorphisme de dualité relative :

{\it $1 ^{\circ}$) Commutation à Frobenius.}
Berthelot définit la catégorie
$F\text{-}D ^\mathrm{b} _\mathrm{coh} (\D ^\dag _{\X} (\hdag T _X) _{\Q})$ de la façon suivante :
les objets sont les couples $(\M,\phi)$, où $\M$ est un objet de
$D ^\mathrm{b} _\mathrm{coh} (\D ^\dag _{\X} (\hdag T _X) _{\Q})$
et $\phi$ est un isomorphisme
$\D ^\dag _{\X} (\hdag T _X) _{\Q}$-linéaire $\M \riso F ^*  \M$.
Les flèches $(\M,\phi) \rightarrow (\M',\phi')$ sont les morphismes
$\D ^\dag _{\X} (\hdag T _X) _{\Q}$-linéaires $\M \rightarrow \M '$ commutant à Frobenius.
Les objets de cette catégorie sont par définition les {\og $F$-complexes\fg} ou 
{\og complexes munis d'une structure de Frobenius \fg}. 
Pour l'instant, quelques théorèmes parmi les plus importants nécessitent une structure de Frobenius. 
Par exemple, le théorème de pleine fidélité de Kedlaya du foncteur restriction de la catégorie des $F$-isocristaux surconvergents
dans celle des $F$-isocristaux convergents sur les $k$-variétés lisses (voir \cite{kedlaya_full_faithfull} ou \cite{kedlaya-semistableII} pour une version relative),
théorème fréquemment utilisé dans \cite{caro_devissge_surcoh} et \cite{caro-2006-surcoh-surcv},
n'est à l'heure actuelle prouvé qu'avec des structures de Frobenius (on conjecture néanmoins cette pleine fidélité sans structure de Frobenius).
D'où notre intérêt concernant
l'extension de $D ^\mathrm{b} _\mathrm{coh} (\D ^\dag _{\X} (\hdag T _X) _{\Q})$ à
$F\text{-}D ^\mathrm{b} _\mathrm{coh} (\D ^\dag _{\X} (\hdag T _X) _{\Q})$
de l'isomorphisme $\chi $, i.e., sa {\it compatibilité à Frobenius}.

{\it $2 ^{\circ}$) Cas non relevable.}
Il est naturel (voir aussi la fin de l'introduction) de s'intéresser au cas non relevable, i.e.,
à l'établissement d'un isomorphisme de la forme
$u _{0T+} \circ \DD _{\X, T _X} (\M) \riso \DD _{\Y,T } \circ u _{0T+} (\M) $
redonnant $\chi$ lorsque $u _0$ se relève.
\\

Cet article s'attaque à ces deux extensions.
En dégageant ses résultats principaux,
détaillons les deux parties qui le composent.

Dans une première partie,
nous travaillons au niveau des schémas.
Dans le cas d'un morphisme {\it non relevable}, quasi-séparé et quasi-compact de schémas lisses sur une base noethérienne de dimension de Krull finie,
nous construisons, sauf pour le morphisme trace,
des analogues de tous les morphismes intervenant dans la construction de l'isomorphisme de dualité relative
de Virrion (\cite{Vir04}).
Surtout, on prouve leur compatibilité à Frobenius et, dans le cas d'une immersion fermée, au changement de niveau.
Rappelons que, par construction de la théorie des $\D$-modules arithmétiques de Berthelot,
les propriétés formelles, i.e., au niveau des schémas formels, se déduisent (notamment par complétion
puis passage à la limite sur le niveau)
de celles au niveau des schémas.

Passons maintenant à la seconde partie.
Soient $d _X$ la dimension de $X$,
$\omega _{\X}= \Omega _\X ^{d _X}$ le faisceau des formes différentielles de degré maximum relatives à
$\X \rightarrow \Spf \V$,
$\O _{\X} (\hdag T _X) _{\Q}$ le faisceau des fonctions sur $\X$ à singularités surconvergentes le long
de $T _X$ et $\omega _{\X} (\hdag T _X) _{\Q}:= \omega _{\X} \otimes _{\O _{\X}} \O _{\X} (\hdag T _X) _{\Q}$.
De même en changeant $\X$ par $\Y$ et $T _X$ par $T$.

$(a)$ Supposons dans ce paragraphe que $u _0$ est une immersion fermée. 
Nous construisons dans ce papier un morphisme trace
$\mathrm{Tr} _{T+} \ : \
u _{0T+} (\omega _{\X} (\hdag T _X) _{\Q} [d _X])
\rightarrow
\omega _{\Y} (\hdag T) _{\Q} [d _Y]$
compatible à Frobenius (voir \ref{can<->tr-gen}).
Lorsque le diviseur est vide et que $u _0$ se relève, on retrouve le morphisme trace de Virrion (\cite{Vir04}).
En utilisant la première partie de ce travail, ce morphisme trace permet de construire l'isomorphisme de dualité relative
$\chi$ : $ u _{0T +} \circ \DD _{\X, T _X} (\M ) \riso \DD _{\Y, T }  \circ u _{0T +}  (\M )$
fonctoriel en $\M \in D ^\mathrm{b} _\mathrm{coh} (\D ^\dag _{\X} (\hdag T _X) _{\Q})$
et vérifiant la condition de transitivité pour le composé de
deux immersions fermées.
On établit de plus la compatibilité à Frobenius de cet isomorphisme $\chi$ (\ref{isodualrelfrob}).

$(b)$ Supposons dans ce paragraphe que $u _{0}$ est projectif ou que $X$ est quasi-projectif. Dans ce cas, 
on construit l'isomorphisme canonique de dualité relative
$\chi\ : \ u _{0T +} \circ \DD _{\X, T _X} ( \M) \riso \DD _{\Y, T }  \circ u _{0 T +}  ( \M)$
fonctoriel en $ \M\in D ^{\mathrm{b}} _\mathrm{coh}(\D ^\dag _{\X } (\hdag T _X) _{\Q} )$
(voir \ref{isodualrelpropreadj}).

Dans les situations $(a)$ ou $(b)$,
il en résulte des morphismes d'adjonction entre image directe et image inverse extraordinaire. 
Dans le cas $(a)$, ceux-ci sont compatibles à Frobenius.

Terminons cette introduction par une application, qui a d'ailleurs été la source
de motivation de cet article.
Revenons à la situation où $u _0$ est une immersion fermée. Sans supposer que
$X$ se relève en un $\V$-schéma formel lisse $\X$,
nous avons construit dans \cite{caro-construction} (et \cite{caro_unite}) un foncteur, noté $\sp _{X\hookrightarrow \Y, T,+}$,
associant, aux isocristaux sur $X\setminus T _X$ surconvergents le long de $T _X$,
des $\D ^\dag _{\Y } (\hdag T ) _{\Q}$-modules à gauche cohérents et à support dans $X$.
Or, la construction de l'isomorphisme de dualité relative dans le cas non relevable (et donc son indépendance par rapport au relèvement choisi)
permet de vérifier que ce foncteur $\sp _{X\hookrightarrow \Y, T,+}$ (construit par recollement via des choix de relèvements d'immersions fermées)
commute aux foncteurs duaux respectifs
(voir la preuve de \cite[2.4.3]{caro-construction}, ce qui permet d'établir \cite[4.3.1]{caro-construction}).
De plus,
afin d'établir la compatibilité à Frobenius de l'isomorphisme de commutation de $\sp _{X\hookrightarrow \Y, T,+}$
aux foncteurs duaux respectifs,
par construction de $\sp _{X\hookrightarrow \Y, T,+}$ qui utilise les images directes,
la compatibilité à Frobenius des isomorphismes de dualité relative
par une immersion fermée de \ref{isodualrelfrob}
est requise.
Rappelons que les foncteurs de la forme $\sp _+$ fournissent le lien entre les $F$-isocristaux surconvergents et
les $F$-$\D$-modules arithmétiques. Ces foncteurs nous ont ainsi permis d'établir que la catégorie des $F$-isocristaux surconvergents
sur les schémas séparés, lisses (sur le corps résiduel de $\V$) est incluse dans celle des $\D$-modules arithmétiques
(pour plus de détails, voir \cite{caro_devissge_surcoh} et \cite{caro-2006-surcoh-surcv})).

\section*{Conventions}

\indent 

$\bullet$ Les fibres spéciales des $\V$-schémas formels lisses seront désignées par les lettres romanes correspondantes.
De plus, si $u $ est un morphisme de $\V$-schémas formels lisses, le morphisme induit
au niveau des fibres spéciales sera noté $u _0$. Si $X \rightarrow S$ est un morphisme lisse de schémas, on note
$d _X$ la dimension de Krull de $X$, $d _{X/S}:=d _X - d_S$ la dimension relative,
$\omega _{X/S}= \Omega _{X/S} ^{d _X}$ ou plus simplement $\omega _X$
le faisceau des formes différentielles de degré maximum relative à $X \rightarrow S$.

$\bullet$ Les indices $\mathrm{qc}$, $\mathrm{tdf}$, $\mathrm{coh}$ et $\mathrm{parf}$
signifient respectivement \textit{quasi-cohérent}, \textit{de Tor-dimension finie}, {\it cohérent}
et {\it parfait}
tandis que $D ^\mathrm{b}$, $D ^+$ et $D ^-$ désignent respectivement les catégories dérivées
des complexes à cohomologie bornée,
bornée inférieurement et bornée supérieurement. Enfin, si $\mathcal{A}$ est un faisceau d'anneaux,
les symboles $\overset{ ^g}{} \mathcal{A}$, $ \mathcal{A} \overset{^d}{}$ ou $\overset{^*}{}\mathcal{A}$
se traduisent respectivement par
$\mathcal{A}$-module {\og à gauche\fg}, {\og à droite\fg} ou {\og à droite ou à gauche\fg}
(par exemple,
$D (\overset{ ^g}{} \mathcal{A})$ indique la catégorie dérivée des complexes de $\mathcal{A}$-modules à gauche).
Si $\AA$ et $\B$ sont deux faisceaux d'anneaux,
par abus de notations, on désignera par $D  (\overset{^*}{}\AA,\overset{^*}{}\B )$ la catégorie
des complexes résolubles (voir \cite[2.1.9]{caro_comparaison}) de $(\AA,\B)$-bimodules (ou $(\AA,\B)$-bimodules à droite ou $(\AA,\B)$-bimodules à gauche selon
la valeur des étoiles).
De plus,
on notera $D _{(.,\,\mathrm{qc})} ( \overset{^*}{}\AA ,\overset{^*}{} \B)$,
la sous catégorie pleine de
$D   ( \overset{^*}{} \AA  ,\overset{^*}{} \B)$ formée des complexes quasi-cohérents pour la structure de $\B$-module induite. 
De même en mettant le point à droite et en remplaçant {\og structure de $\B$-module\fg} par {\og structure de $\AA$-module\fg} ; en remplaçant respectivement
{\og $\mathrm{qc}$\fg} et {\og quasi-cohérent\fg} par {\og $\mathrm{tdf}$\fg} et {\og de Tor-dimension finie\fg} etc.
Comme les catégories de complexes sont par défaut les catégories dérivées,
on écrit sans ambiguïté $\mathrm{Hom} _{\AA}(-,-)$ pour $\mathrm{Hom} _{D (\AA)}(-,-)$.

$\bullet$ Nous travaillerons de préférence avec des modules à droite. Tous les théorèmes démontrés dans cet article
sont valables par passage de droite à gauche, e.g. on dispose d'un isomorphisme de dualité relative pour les modules à gauche.

$\bullet$ L'isomorphisme canonique de commutation à Frobenius de l'image directe
que nous utiliserons dans ce papier sera,
par passage de gauche à droite, celui construit par Berthelot (voir \cite[3.4.4]{Be2} pour la version \og à gauche\fg).

Cependant, comme cela est expliqué dans la remarque \ref{remacompisofrobu+}, dans le cas d'une {\it immersion fermée},
la compatibilité à Frobenius de l'isomorphisme de dualité relative de \ref{isodualrelfrob}
implique que l'isomorphisme de commutation à Frobenius de l'image directe de \cite[3.4.4]{Be2}
 est identique à celui construit dans \cite[1.2.12]{caro_courbe-nouveau}.
Cela fournit donc une unification de ces deux constructions.

\section{\label{ca}Préliminaire sur la construction de l'isomorphisme de dualité relative d'un morphisme non relevable de schémas et sa compatibilité à Frobenius }

Soient deux entiers $m,\,s\geq 0$, $S$ un schéma noethérien de dimension de Krull finie et
muni d'un $m$-PD-idéal quasi-cohérent
$(\mathfrak{a},\,\mathfrak{b},\,\alpha)$ et $m$-PD-nilpotent,
$S _0 = V (\mathfrak{a})$ et
$u _0$ : $X _0\rightarrow Y _0$ un morphisme
quasi-séparé et quasi-compact
de $S _0$-schémas lisses. On suppose de plus que $p\in \mathfrak{a}$
et $p$ est nilpotent sur $S$ (l'idéal $\mathfrak{a}$ est ainsi un nilidéal).
On notera $F ^s _{S _0}$ : $S _0 \rightarrow S _0$ la puissance $s$-ème du morphisme de Frobenius absolu,
$u _0 '$ : $X _0 ^{(s)} \rightarrow Y _0 ^{(s)}$ le morphisme déduite de $u _0$
par le changement de base de $F ^s _{S _0}$,
$F ^s _{X _0/S _0}$ : $X _0 \rightarrow X _0 ^{(s)}$ le morphisme de Frobenius relatif (de même en remplaçant
$X $ par $Y$). 

Considérons les deux hypothèses suivantes.
\begin{itemize}
  \item {\it Cas général} : On suppose qu'il existe respectivement des relèvements $X$, $X'$, $Y$ et $Y'$
de $X _0$, $X _0 ^{(s)}$, $Y _0$ et $Y _0 ^{(s)}$.
  \item  {\it Cas relevable} : On suppose qu'il existe respectivement des relèvements
$u$ : $X \rightarrow Y$, $u '$ : $X ' \rightarrow Y '$,
$F _X$ : $X \rightarrow X'$ et $F _Y$ : $Y \rightarrow Y'$ de
$u _0$, $ u _0 '$, $F ^s _{X _0/S _0}$ et $F ^s _{Y _0/S _0}$ tels que $F _Y \circ u = u' \circ F _X$.
On remarque que comme $\mathfrak{a}$ est un nilidéal, les foncteurs $u _*$, $u ' _*$ et $u _{0*}$
sont égaux.
\end{itemize}

Nous nous plaçons ici dans le {\og cas général \fg}. 
Si aucune confusion n'est à craindre, on notera $F$ pour $F ^s _{X _0/S _0}$ ou $F ^s _{Y _0/S _0}$.
Nous disposons du foncteur $u _0 ^*$
(construit par recollement et localement canoniquement isomorphe à $u^*$)
de la catégorie des $\D _Y ^{(m)}$-modules à gauche dans celle des $\D _X ^{(m)}$-modules à gauche ainsi
que du foncteur $F ^*$ (resp. $F ^\flat$) de la catégorie des
$\D _{X'} ^{(m)}$-modules à gauche (resp. à droite) dans celle des $\D _X ^{(m+s)}$-modules à gauche (resp. à droite).

Pour tout entier $m$, soient $\B _{Y'}^{(m)}$ une $\O _{Y'}$-algèbre munie d'une structure compatible de 
$\D ^{(m)} _{Y'}$-module à gauche tels que, pour tout $m' \geq m$, le morphisme canonique 
 $\B _{Y'}^{(m)}\to \B _{Y'}^{(m')}$ soit un morphisme $\D ^{(m)} _{Y'}$-linéaire. 
On pose $\B _{Y}^{(m+s)}:= F ^{*} (\B _{Y'} ^{(m)})$, $\B _{X'} ^{(m)}:= u _{0} ^{*} (\B _{Y'})$,
$\B _{X}^{(m+s)}:= F ^{*} (\B _{X'}^{(m)})$.
On note $\widetilde{\omega} _{X'} ^{(m)} := \B _{X'} ^{(m)} \otimes _{\O _{X'}} \omega _{X'}$, 
$\widetilde{\omega} _{X} ^{(m+s)} := \B _{X} ^{(m+s)} \otimes _{\O _{X}} \omega _{X}$, de même en remplaçant $X$ par $Y$.
Pour alléger les notations (les faisceaux de la forme $\B$ étant fixés dans tout ce chapitre), 
pour tout $ \M ' \in D  ( \B _{X'} ^{(m)} \otimes _{\O _{X'}} \D _{X'}^{(m)} \overset{^{\mathrm{d}}}{})$, le foncteur dual de niveau $m$ de $\M'$
se note simplement
$$ \DD ^{(m)} (\M') :=
 \widetilde{\omega} _{X'} ^{(m)}
 \otimes _{\B _{X'} ^{(m)}} 
 \R \mathcal{H} om _{\B _{X'} ^{(m)} \otimes _{\O _{X'}} \D _{X'}^{(m)}} 
 ( \M ', \B _{X'} ^{(m)} \otimes _{\O _{X'}} \D _{X'}^{(m)})
 [d _X].$$

De plus, le foncteur image directe (relativement à $\B ^{(m)} _{X'}$) de niveau $m$ 
de $ \M ' \in D ( \B _{X'} ^{(m)} \otimes _{\O _{X'}} \D _{X'}^{(m)}\overset{^{\mathrm{d}}}{})$ par $u ' _{0}$
se note simplement : 
$$u ^{\prime (m)} _{0+} ( \M ') := \R u _{0*} ( \M ' \otimes ^\L _{\B _{X'} ^{(m)} \otimes _{\O _{X'}}\D ^{(m)} _{X'}} (\B _{X'} ^{(m)} \otimes _{\O _{X'}} \D ^{(m)} _{X '\rightarrow Y'})),$$
où $\D ^{(m)} _{X' \rightarrow Y'}:=u _0 ^* \D ^{(m)} _{Y'}$ est, par passage de gauche à droite de la remarque \cite[3.4.3]{Be2},
de Tor-dimension finie comme $\D ^{(m)} _{X'}$-module.

{\bf Résumé du chapitre:} 
Pour tout $ \M ' \in D ^- ( \B _{X'} ^{(m)} \otimes _{\O _{X'}} \D _{X'}^{(m)}\overset{^{\mathrm{d}}}{})$,
nous construisons le morphisme canonique 
\begin{equation}
\label{1}
u _{0+} ^{\prime (m)}\circ  \DD ^{(m)} (  \M ')
\to \R \mathcal{H} om _{ \B _{Y'} ^{(m)} \otimes _{\O _{Y'}} \D _{Y'}^{(m)}}
(u _{0+} ^{\prime (m)} ( \M ' ), u _{0+} ^{\prime (m)} (\widetilde{\omega} _{X'} ) \otimes _{\B ^{(m)} _{Y'}} (\B _{Y'} ^{(m)} \otimes _{\O _{Y'}} \D _{Y'}^{(m)}))[d _{X}].
\end{equation}
Nous vérifions de plus sa commutation à Frobenius (voir \ref{theo1-7}) et au changement de niveaux (voir
\ref{theo17niveau}).

Pour déduire de \ref{1} la construction du morphisme (qui est un isomorphisme si $\M$ est un complexe parfait) de dualité relative
$u _{0+} ^{\prime (m)} \circ \DD ^{(m)} (  \M ') \to \DD ^{(m)} \circ  u _{0+} ^{\prime (m)} (  \M ')$, il ne reste plus qu'à construire le morphisme 
trace $u _{0+} ^{\prime (m)} ( \widetilde{\omega} _{X'} ) [d _{X}] \to ( \widetilde{\omega} _{Y'} ) [d _{Y}]$, ce que nous traiterons dans le prochain chapitre dans le cas d'une immersion fermée.

\subsection{Rappels et notations}

Jusqu'à la section \ref{ccn} incluse, le niveau est égal à $m$ sur $X'$ ou $Y'$ et à $m+s$ sur $X$ ou $Y$. 
Nous noterons alors dans toutes ces sections, $\D _{?}:= 
\B _{X} ^{(m+s)} \otimes _{\O _{X}} \D ^{(m+s)} _{?}$, 
$\D _{?}:= 
\B _{X'} ^{(m)} \otimes _{\O _{X'}} \D ^{(m)} _{?}$, où $\D ^{(m+s)} _{?}$ et $\D ^{(m)} _{?}$ sont respectivement des faisceaux sur $X$ et sur $X'$. De même, en remplaçant $X$ par $Y$. 
De plus, comme le niveau est fixé, on omettra de l'indiquer et on notera simplement $u _{0+} ^{\prime}$, $u _{0+} $, $\DD$, $\B _{X'}$,  $\B _{X}$ etc.

La construction du morphisme \ref{1}
et sa compatibilité à Frobenius 
se décomposent en trois parties 
(voir \ref{prop1-2}, \ref{prop2-4} et \ref{prop4-7}).

\begin{vide}
  [Convention] \label{convention}
Soient $\M '\in D ^- (\D _{X'} \overset{^{\mathrm{d}}}{})$, $\NN ' \in D ^+ (\D _{X'} \overset{^{\mathrm{d}}}{},\D _{X'} \overset{^{\mathrm{d}}}{}) $.
Pour calculer $\R \mathcal{H} om _{\D _{X'}} (  \M ', \NN ')$ on utilisera la structure {\it gauche} (de complexe de $\D _{X'}$-modules à droite) de $\NN '$.

Par exemple,
l'isomorphisme de transposition de \cite[1.3.3]{Be2} induit fonctoriellement le suivant
$\DD ( \M ' ) \riso
\R \mathcal{H} om _{\D _{X'}} (  \M ', \widetilde{\omega} _{X'} \otimes _{\B _{X'}}\D _{X'}) [d _X]$
(la structure gauche de $\widetilde{\omega} _{X'} \otimes _{\B _{X'}}\D _{X'}$ est par convention la structure induite par le produit tensoriel).
\end{vide}

\begin{vide}
Nous utiliserons fréquemment,
pour tous $\M '\in D ^- (\D _{X'} \overset{^{\mathrm{d}}}{})$
et $\E '\in D ^- (\overset{^{\mathrm{g}}}{} \D _{X'} )$,
les isomorphismes
$\D _X \riso F ^* F ^\flat \D _{X'}$ (cas relevable : \cite[2.5.2.1]{Be2}, le cas non relevable découle de
\cite[3.4.2.(ii)]{Be2} appliqué à l'identité),
$\mu _{\E' }$ :
$\widetilde{\omega} _{X} \otimes _{\B _X} F ^* \E' \riso F ^\flat (\widetilde{\omega} _{X'} \otimes _{\B _{X'}} \E')$ (\cite[2.4.5.1]{Be2})
et
$F ^\flat  \M ' \otimes _{\D _X} F ^* \E ' \riso  \M ' \otimes _{\D _{X'}} \E '$ (\cite[2.5.7]{Be2} reste
valable dans le cas non relevable).

Pour tout $\NN ' \in D (u ^{\text{-}1} _0 \D _{Y'}  \overset{^{\mathrm{d}}}{})$, on pose
$F ^\flat (\NN ') = \NN ' \otimes _{u ^{\text{-}1} _0 \D _{Y'}} u ^{\text{-}1} _0 F ^\flat \D _{Y'}$ (de même avec des tildes).
On obtient l'isomorphisme
$\D _{X \rightarrow Y} \riso u ^* _0 F ^* F ^\flat (\D _{Y'})
\riso F ^* F ^\flat \D _{X' \rightarrow Y'}$ (voir \cite[3.4.2.3]{Be2}).

Pour tous $\FF '\in D  _{\mathrm{tdf},\mathrm{qc}} (\overset{^g}{} \D _{Y'})$
(resp. $\FF '\in D ^- _{\mathrm{qc}} (\overset{^g}{} \D _{Y'})$) et
 $\NN ' \in D ( u ^{-1}\D _{Y'} \overset{^d}{})$ (resp. $\NN ' \in D ^-( u ^{-1}\D _{Y'} \overset{^d}{})$),
on notera {\og$\mathrm{proj}$\fg}
 l'isomorphisme de projection induit par $u $ (ou ceux qui s'en déduisent par fonctorialité) :
\begin{equation}
\label{iso-proj-D}
\mathrm{proj}  \ : \ \R u _{0*} \NN '\otimes _{\D _{Y'}} ^\L \FF '\riso
 \R u _{0*} (\NN '\otimes ^\L _{u _0^{-1}\D _{Y'}} u _0 ^{-1} \FF '),
\end{equation}
 dont la construction
 est analogue à \cite[II.5.6]{HaRD}.
Et de même au niveau supérieur, i.e.,
en substituant $\D _{Y}$ à $\D _{Y'}$.
Par exemple, 
en remplaçant respectivement $\NN'$ par $\M ' \otimes ^\L _{\D_{X'}} \D _{X '\rightarrow Y'}$
et $\FF'$ par $F ^\flat \D _{Y'}$, on obtient l'isomorphisme en bas à droite de \ref{u+comFrob}
(de manière analogue à \cite[2.4]{Be2}, le foncteur $F ^\flat $ est exact et
$\D _{X '\rightarrow Y'}\otimes ^\L _{u ^{\text{-}1} _0 \D _{Y'}} u ^{\text{-}1} _0 F ^\flat \D _{Y'}
\riso F ^\flat \D _{X '\rightarrow Y'} $).

On désignera encore par $\mathrm{proj}$, le morphisme de projection
$$\mathrm{proj}  \ : \ \R u _{0*} \NN '\otimes _{\B _{Y'}} ^\L \FF '\riso
 \R u _{0*} (\NN '\otimes ^\L _{u _0^{-1}\B _{Y'}} u _0 ^{-1} \FF '),$$
 qui sera notamment utilisé dans \ref{u+!proj}.

\end{vide}

\begin{vide}
  \label{videu+comFrob}
Berthelot a construit l'isomorphisme de commutation à Frobenius de l'image directe par $u_0$ de
$\M '\in D (\D _{X'} \overset{^{\mathrm{d}}}{})$ (voir \cite[3.4.4]{Be2} pour la version \og à gauche\fg)
via le composé :
\begin{equation}
  \label{u+comFrob}
\xymatrix @R=0,3cm {
{u _{0+} (F ^\flat  \M ')}
\ar@{.>}[dd] ^-\sim
\ar[r] _-\sim
&
{\R u _{0*} ( F ^\flat  \M ' \otimes ^\L _{\D _X} F ^* F ^\flat \D _{X '\rightarrow Y'})}
\ar[d] _-\sim
\\
&
{\R u _{0*} (  \M ' \otimes ^\L _{\D _{X'}} F ^\flat \D _{X '\rightarrow Y'})}
\\
{F ^\flat u _{0+} '  ( \M ')}
\ar@{=}[r]
&
{\R u _{0*} (  \M ' \otimes ^\L _{\D_{X'}} \D _{X '\rightarrow Y'}) \otimes _{\D _{Y'}} F ^\flat \D _{Y'}.}
\ar[u] _-\sim ^-{\mathrm{proj}}
}
\end{equation}
\end{vide}

\begin{vide}
  \label{defDualcomFrob}
Rappelons la construction de Virrion (voir \cite[II.3.2.ii)]{virrion})
de l'isomorphisme de commutation à Frobenius du foncteur dual.
  Pour tout $ \M ' \in D ( \D _{X'} \overset{^{\mathrm{d}}}{})$,
on définit d'abord,
par commutativité du diagramme
\begin{equation}
  \label{DualcomFrob}
\xymatrix @R=0,3cm {
{ \R \mathcal{H} om _{\D _X} ( F ^\flat  \M ', \widetilde{\omega} _X \otimes _{\B _X} \D _X)}
\ar[r] _-\sim \ar@{.>}[d] ^-\sim &
{\R \mathcal{H} om _{\D _X} ( F ^\flat  \M ', F ^\flat _{\mathrm{g}} F ^\flat _{\mathrm{d}} (\widetilde{\omega} _{X'} \otimes _{\B _{X'}}\D _{X'}))}
\\
{F ^\flat \R \mathcal{H} om _{\D _{X'}} (  \M ', \widetilde{\omega} _{X'} \otimes _{\B _{X'}}\D _{X'} )}
&
{\R \mathcal{H} om _{\D _{X'}} (  \M ', F ^\flat _{\mathrm{d}} (\widetilde{\omega} _{X'} \otimes _{\B _{X'}}\D _{X'})),}
\ar[l] ^-\sim \ar[u] ^-{F ^\flat} _-\sim
}
\end{equation}
où les symboles {\og $\mathrm{d}$\fg} et {\og $\mathrm{g}$\fg} signifient que l'on utilise respectivement la structure droite et la structure gauche du $\D _{X'}$-bimodule à droite $\widetilde{\omega} _{X'} \otimes _{\B _{X'}}\D _{X'}$, 
l'isomorphisme
$ \R \mathcal{H} om _{\D _X} ( F ^\flat  \M ', \widetilde{\omega} _X \otimes _{\B _X} \D _X)
\riso F ^\flat \R \mathcal{H} om _{\D _{X'}} (  \M ', \widetilde{\omega} _{X'} \otimes _{\B _{X'}}\D _{X'} )$.
Il en dérive, via les isomorphismes de transposition de \cite[1.3.3]{Be2},
$\DD F ^\flat ( \M ' ) \riso F ^\flat \DD ( \M ')$.

\end{vide}

\subsection{Première partie de la construction de l'isomorphisme de dualité et sa compatibilité à Frobenius}
\begin{vide}
  \label{preprop1-2}
Par \ref{videu+comFrob} et \ref{DualcomFrob},
pour tout $ \M '  \in D ( \D _{X'} \overset{^{\mathrm{d}}}{})$,
on dispose par construction du diagramme commutatif :
\begin{equation}
  \label{1-2bis}
\xymatrix @R=0,3cm {
{u _{0+} \DD ( F ^\flat  \M ')}
\ar[r] ^-\sim
\ar[d] ^-\sim
&
{\R u _{0*} \left( \R \mathcal{H} om _{\D _X} ( F ^\flat  \M ', \widetilde{\omega} _X \otimes _{\B _X} \D _X)
\otimes ^\L _{\D _X} \D _{X\rightarrow Y} \right) [d _{X}]}
\ar[d] ^-\sim
\\
{u _{0+} F ^\flat  \DD (  \M ')}
\ar[r] ^-\sim
\ar[d] ^-\sim
&
{\R u _{0*} \left( F ^\flat   \R \mathcal{H} om _{\D _{X'}} ( \M ', \widetilde{\omega} _{X' }\otimes _{\B _{X'}} \D _{X'})
\otimes ^\L _{\D _X} \D _{X\rightarrow Y} \right) [d _{X}]}
\ar[d] ^-\sim
\\
{F ^\flat  u _{0+} ' \DD (  \M ')}
\ar[r] ^-\sim
&
{F ^\flat  \R u _{0*} \left ( \R \mathcal{H} om _{\D _{X'}} (  \M ', \widetilde{\omega} _{X'} \otimes _{\B _{X'}}\D _{X'})
\otimes ^\L _{\D _{X'}} \D _{X '\rightarrow Y'}\right )[d _{X}],}
}
\end{equation}
dont les isomorphismes horizontaux sont
induits par les isomorphismes de transposition de $\widetilde{\omega} _X \otimes _{\B _X} \D _X$
ou de $\widetilde{\omega} _{X'} \otimes _{\B _{X'}}\D _{X'}$.
\end{vide}

\begin{prop}
 \label{prop1-2}
Pour tout $ \M ' \in D ( \D _{X'} \overset{^{\mathrm{d}}}{})$,
on bénéficie du diagramme commutatif suivant :
\begin{equation}
  \label{1-2}
\xymatrix @R=0,3cm {
{\R u _{0*} ( \R \mathcal{H} om _{\D _X} ( F ^\flat  \M ', \widetilde{\omega} _X \otimes _{\B _X} \D _X)
\otimes ^\L _{\D _X} \D _{X\rightarrow Y} ) }
\ar[r]
\ar[d] ^-\sim
&
{\R u _{0*} \R \mathcal{H} om _{\D _X} ( F ^\flat  \M ', \widetilde{\omega} _X \otimes _{\B _X} \D _{X\rightarrow Y})}
\ar[d] ^-\sim
\\
{\R u _{0*} ( \R \mathcal{H} \underset{\D _X}{om}
( F ^\flat  \M ', F ^\flat _{\mathrm{d}} F ^\flat _{\mathrm{g}} (\overset{^\mathrm{t}}{} \D _{X'}))
\underset{\D _X }{\overset{\L}{\otimes}}
F ^* F ^\flat \D _{X '\rightarrow Y'} )}
\ar[r]
&
{\R u _{0*} \R \mathcal{H} \underset{\D _X}{om}
( F ^\flat  \M ', F ^\flat _{\mathrm{d}} F ^\flat _{\mathrm{g}} (\overset{^\mathrm{t}}{} \D _{X'})
\underset{\D _X }{\otimes} F ^* F ^\flat \D _{X '\rightarrow Y'})}
\\
{\R u _{0*} ( \R \mathcal{H} om _{\D _{X'}} (  \M ', \overset{^\mathrm{t}}{} \D _{X'})
\otimes ^\L _{\D _{X'}} F ^\flat \D _{X '\rightarrow Y'})}
\ar[u] ^-\sim \ar[r]
&
{\R u _{0*} \R \mathcal{H} om _{\D _{X'}}
(  \M ', \overset{^\mathrm{t}}{} \D _{X'} \otimes _{\D _{X'} } F ^\flat \D _{X '\rightarrow Y'})}
\ar[u] ^-\sim
\\
{F ^\flat  \R u _{0*} ( \R \mathcal{H} om _{\D _{X'}} (  \M ', \overset{^\mathrm{t}}{} \D _{X'})
\otimes ^\L _{\D _{X'}} \D _{X '\rightarrow Y'})}
\ar[r] \ar[u] ^-\sim _{\mathrm{proj}}
&
{F ^\flat  \R u _{0*} \R \mathcal{H} om _{\D _{X'}}
(  \M ', \overset{^\mathrm{t}}{} \D _{X'} \otimes _{\D _{X'} } \D _{X '\rightarrow Y'}),}
\ar[u] ^-\sim _{\mathrm{proj}}
}
\end{equation}
\normalsize
où $\overset{^\mathrm{t}}{} \D _{X'}$ est le $\D _{X'}$-bimodule à droite
$\widetilde{\omega} _{X'} \otimes _{\B _{X'}} \D _{X'}$.
\end{prop}
\begin{proof}
  Pour tous $ \M ' \in D ( \D _{X'} \overset{^{\mathrm{d}}}{})$,
$\NN ' \in D ^+ ( \D _{X'} \overset{^{\mathrm{d}}}{}, \D _{X'} \overset{^{\mathrm{d}}}{})$
et $\E ' \in D _{\mathrm{tdf}} ( \overset{^{\mathrm{g}}}{} \D _{X'} )$,
il résulte de \cite[2.1.12(iii)]{caro_comparaison} la commutativité du diagramme ci-dessous :
\begin{equation}
  \label{pre1-2}
\xymatrix @R=0,3cm {
{\R \mathcal{H} om _{\D _X} ( F ^\flat  \M ' , F ^\flat _{\mathrm{g}} F ^\flat _{\mathrm{d}} \NN ' )
\otimes ^\L _{\D _X} F ^* \E '}
\ar[r]
&
{\R \mathcal{H} om _{\D _X}
( F ^\flat  \M ' , F ^\flat _{\mathrm{g}} F ^\flat _{\mathrm{d}} \NN ' \otimes ^\L _{\D _X} F ^* \E ')}
\ar[d] ^-\sim
\\
{\R \mathcal{H} om _{\D _{X'}} (  \M ' , F ^\flat _{\mathrm{d}} \NN ' )
\otimes ^\L _{\D _X} F ^* \E '}
\ar[u] ^-{F ^\flat} _-\sim
\ar[d] ^-\sim
&
{\R \mathcal{H} om _{\D _X}
( F ^\flat  \M ' , F ^\flat _{\mathrm{g}} \NN ' \otimes ^\L _{\D _{X'}} \E ')}
\\
{\R \mathcal{H} om _{\D _{X'}} (  \M ' ,  \NN ' )
\otimes ^\L _{\D _{X'}} \E '}
\ar[r]
&
{\R \mathcal{H} om _{\D _{X'}}
(  \M ' , \NN ' \otimes ^\L _{\D _{X'}} \E ').}
\ar[u] ^-{F ^\flat} _-\sim
}
\end{equation}
La commutativité de \ref{pre1-2} entraîne celle du carré du milieu de \ref{1-2}.
Celle du carré du haut de \ref{1-2} s'acquiert par fonctorialité
tandis que celle du carré du bas se vérifie par fonctorialité et transitivité des morphismes de la forme
de ceux horizontaux de \ref{pre1-2}.
 \end{proof}

\begin{rema}
\label{rema1-2}
Le morphisme de droite de \ref{1-2bis} est (au décalage $[d _{X}]$ près)
le composé de gauche de \ref{1-2}. On pourra donc composer les contours de ces deux diagrammes.
\end{rema}

\subsection{Deuxième partie de la construction de l'isomorphisme de dualité et sa compatibilité à Frobenius}

\begin{vide}\label{pre24}
  Soient $ \M ' \in D ^- ( \D _{X'} \overset{^{\mathrm{d}}}{})$ et
  $\NN ' \in D ^+ ( \D _{X'} \overset{^{\mathrm{d}}}{})$.
  Comme $\D _{X '\rightarrow Y'}$ est de Tor-dimension finie sur $\D _{X '}$, on définit en résolvant platement
$\D _{X '\rightarrow Y'}$, injectivement $\NN'$, via
$\mathcal{H} om _{u ^{-1} _0 \D _{Y'} } (-,-)\rightarrow \R \mathcal{H} om _{u ^{-1} _0 \D _{Y'} } (-,-)$,
 la flèche :
$\R \mathcal{H} om _{\D _{X'}} ( \M ' ,\NN ')  \rightarrow
\R \mathcal{H} om _{u ^{-1} _0 \D _{Y'} }
   ( \M ' \otimes ^\L _{\D _{X'}} \D _{X '\rightarrow Y'},\NN '\otimes ^\L _{\D _{X'}} \D _{X '\rightarrow Y'})$.
On note $\rho _{2}$ l'application définie via le diagramme commutatif :
  $$\xymatrix @R=0,3cm {
{\R u _{0*} \R \mathcal{H} om _{\D _{X'}} ( \M ' ,\NN ') }
\ar[r]
\ar@{.>}[rd] _-{\rho _{2}}
&
{   \R u _{0*} \R \mathcal{H} om _{u ^{-1} _0 \D _{Y'} }
   ( \M ' \otimes ^\L _{\D _{X'}} \D _{X '\rightarrow Y'},\NN '\otimes ^\L _{\D _{X'}} \D _{X '\rightarrow Y'})}
\ar[d] ^-{\R u _{0*}}
\\
&
{ \R \mathcal{H} om _{\D _{Y'} } (u _{0+} '  ( \M ' ),u _{0+} '   (\NN '))}
},$$
où $\R u _{0*} $désigne le foncteur
$\R u _{0*} \R \mathcal{H} om _{u ^{-1} _0 \D _{Y'} }(-,-)
\rightarrow \R \mathcal{H} om _{\D _{Y'} }(\R u _{0*} -,\R u _{0*} -)$ construit de manière analogue à 
\cite[II.5.5]{HaRD}.
De même en remplaçant $\D _{X'}$ par $\D _X$, $\D _{X '\rightarrow Y'}$
par $\D _{X \rightarrow Y}$ etc.

Afin de prouver la proposition \ref{prop2-4}, nous aurons besoin des lemmes \ref{lemmpre24}
et \ref{lemmu_+fonccomflatappl} ci-après.
\end{vide}

\begin{lemm}\label{lemmpre24}
  Avec les notations de \ref{pre24}, le diagramme canonique
\begin{equation}\label{lemmpre24diaglem}
  \xymatrix @R=0,3cm {
  {\R u _{0*} \R \mathcal{H} om _{\D _X} (F ^\flat  \M ' ,F ^\flat \NN ')}
  \ar[r] ^-{\rho _{2}}
  &
  {\R \mathcal{H} om _{\D _Y } (u _{0+} (F ^\flat  \M ' ),u _{0+}  (F ^\flat  \NN '))}
  \ar[d] ^-\sim
  \\
  &
  {\R \mathcal{H} om _{\D _Y } (F ^\flat  u _{0+} '  ( \M ' ),F ^\flat  u _{0+} '   ( \NN '))}
  \\
  {\R u _{0*} \R \mathcal{H} om _{\D _{X'}} ( \M ' ,\NN ')}
  \ar[uu] ^-{F ^\flat} _-\sim
  \ar[r] ^-{\rho _2}
  &
  {\R \mathcal{H} om _{\D _{Y'} } (u _{0+} '  ( \M ' ),u _{0+} '   (\NN '))}
  \ar[u] ^-{F ^\flat} _-\sim
  }
\end{equation}
est commutatif.
\end{lemm}
\begin{proof}
Via l'isomorphisme canonique
\begin{gather}\notag
  F ^\flat  \M '  \otimes ^\L _{\D _X} \D _{X\rightarrow Y}
\riso F ^\flat  \M '  \otimes ^\L _{\D _X} F ^* F ^\flat \D _{X '\rightarrow Y'}
\riso  \M '  \otimes ^\L _{\D _{X'}} F ^\flat \D _{X '\rightarrow Y'}
= F ^\flat ( \M '  \otimes ^\L _{\D _{X'}} \D _{X '\rightarrow Y'}),
\end{gather}
\normalsize
de même pour $\NN'$ à la place de $\M'$, 
on obtient la flèche de droite du haut du diagramme :
\begin{equation}
\label{lemmpre24diag1}
  \xymatrix @R=0,3cm @C=0cm {
  {\R u _{0*} \R \mathcal{H} om _{\D _X} (F ^\flat  \M ' ,F ^\flat \NN ')}
  \ar[r]
  &
  {\R u _{0*} \R \mathcal{H} om _{u ^{\text{-}1} _0\D _Y }
  (F ^\flat  \M '  \otimes ^\L _{\D _X} \D _{X\rightarrow Y},
  F ^\flat \NN '  \otimes ^\L _{\D _X} \D _{X\rightarrow Y})}
  \ar[d] ^-\sim
  \\
  &
  {\R u _{0*} \R \mathcal{H} om _{u ^{\text{-}1} _0\D _Y }
  (F ^\flat ( \M '  \otimes ^\L _{\D _{X'}}\! \! \! \D _{X '\rightarrow Y'}),
  F ^\flat (\NN '  \otimes ^\L _{\D _{X'}} \! \! \!\D _{X '\rightarrow Y'}))}
  \\
  {\R u _{0*} \R \mathcal{H} om _{\D _{X'}} ( \M ' ,\NN ')}
  \ar[uu] ^-{F ^\flat} _-\sim
  \ar[r]
  &
  {\R u _{0*} \R \mathcal{H} om _{u ^{\text{-}1} _0 \D _{Y'} }
   ( \M ' \otimes ^\L _{\D _{X'}} \D _{X '\rightarrow Y'},\NN '\otimes ^\L _{\D _{X'}} \D _{X '\rightarrow Y'}).}
  \ar[u] ^-{F ^\flat} _-\sim
  }
\end{equation}
\normalsize
Celui-ci est commutatif. En effet, il suffit de l'établir sans $\R u _{0*}$. Résolvons injectivement
$\NN'$ et platement $\D _{X '\rightarrow Y'}$. Via l'isomorphisme canonique $\D _{X\rightarrow Y}
\riso  F ^* F ^\flat \D _{X '\rightarrow Y'} $, il en résulte une résolution plate de $\D _{X\rightarrow Y}$ (les foncteurs $F  ^{*}$ et $F ^{\flat}$ préservent la platitude). Par construction et comme $F ^\flat$ préserve l'injectivité
(on s'en sert notamment pour obtenir la commutation à $F ^\flat$ de
$\mathcal{H} om _{u ^{-1} _0 \D _{Y'} } (-,-)\rightarrow \R \mathcal{H} om _{u ^{-1} _0 \D _{Y'} } (-,-)$),
on se ramène à le vérifier sans les $\R$ ni $\L$, ce qui est immédiat.
En notant $\phi $ le foncteur
$- \otimes ^\L _{\D _{X'}} \D _{X '\rightarrow Y'}$,
considérons le diagramme suivant : 
\begin{equation}
\label{lemmpre24diag2}
  \xymatrix @R=0,3cm @C=0,5cm {
  {\R u _{0*} \R \mathcal{H} \underset{u ^{\text{-}1} _0\D _Y }{om}
  (F ^\flat  \M '  \underset{\D _X}{\overset{\L}{\otimes}} \D _{X\rightarrow Y},
  F ^\flat \NN '  \underset{\D _X}{\overset{\L}{\otimes}} \D _{X\rightarrow Y})}
  \ar[r] ^-{\R u _{0*}}
  \ar[d] ^-\sim
  &
 {\R \mathcal{H} om _{\D _Y }
  (\R u _{0*}  (F ^\flat  \M '  \underset{\D _X}{\overset{\L}{\otimes}} \D _{X\rightarrow Y}),
  \R u _{0*}  (F ^\flat \NN '  \underset{\D _X}{\overset{\L}{\otimes}} \D _{X\rightarrow Y}))}
  \ar[d] ^-\sim
  \\
  {\R u _{0*} \R \mathcal{H} \underset{u ^{\text{-}1} _0\D _Y }{om}
  (F ^\flat   \phi ( \M ' ) ,
  F ^\flat \phi (\NN '  ))}
  \ar[r] ^-{\R u _{0*}}
  &
  {\R \mathcal{H} om _{\D _Y }
  (\R u _{0*} F ^\flat \phi  ( \M '  ),
  \R u _{0*} F ^\flat \phi  (\NN '  ))}
  \\
  {\R u _{0*} \R \mathcal{H} \underset{u ^{\text{-}1} _0\D _{Y'} }{om}
   ( \phi  (\M ') , \phi  (\NN '))}
  \ar[u] ^-{F ^\flat} _-\sim \ar[rd] ^-{\R u _{0*}}
  &
  {\R \mathcal{H} om _{\D _Y }
  (F ^\flat \R u _{0*} \phi  ( \M '  ),
  F ^\flat \R u _{0*} \phi  (\NN '  ))}
  \ar[u] ^-{\mathrm{proj}} _-\sim
  \\
  &
  {\R \mathcal{H} om _{\D _{Y'} }
  (\R u _{0*} \phi  ( \M '  ),
  \R u _{0*} \phi  (\NN '  )).}
  \ar[u] ^-{F  ^\flat} _-\sim
  }
\end{equation}
\normalsize
On obtient par fonctorialité la commutativité du carré (en haut) de \ref{lemmpre24diag2}.
En résolvant $\phi  (\NN ')$ par des $u ^{\text{-}1} _0\D _{Y'} $-modules injectifs et $\phi  (\M ' )$ par des $u ^{\text{-}1} _0\D _{Y'} $-modules 
$u _{0*}$-acycliques, on résout le complexe $\R u _{0*} \R \mathcal{H} om _{u ^{\text{-}1} _0\D _{Y'} }
   ( \phi  (\M ') , \phi  (\NN '))$. Comme $F ^{\flat}$ préserve l'injectivité et la $u _{0*}$-acyclicité, 
   il en est de même du deuxième terme de gauche du trapèze de \ref{lemmpre24diag2}.
   Comme ce trapèze est commutatif sans les symboles $\R$ et $\L$, on vérifie alors sa commutativité. 
Le diagramme \ref{lemmpre24diag2} est donc commutatif.
Comme le composé de \ref{lemmpre24diag1} et \ref{lemmpre24diag2}
donne \ref{lemmpre24diaglem}, on conclut la preuve.
 \end{proof}

\begin{vide}  \label{u_+fonccomflatnot}
  Soit
$\NN ' \in D ^+ ( \D _{X'} \overset{^{\mathrm{d}}}{}, u _0 ^{-1} \D _{Y'} \overset{^{\mathrm{d}}}{})$.
Par fonctorialité
$\NN '\otimes ^\L _{\D _{X'}} \D _{X '\rightarrow Y'}\in
D ^+ ( u _0 ^{-1} \D _{Y'} \overset{^{\mathrm{d}}}{}, u _0 ^{-1} \D _{Y'} \overset{^{\mathrm{d}}}{})$
et $ u _{0+} '  (\NN ') = \R u _{0*} (\NN '\otimes ^\L _{\D _{X'}} \D _{X '\rightarrow Y'})
\in D ^+ ( \D _{Y'} \overset{^{\mathrm{d}}}{}, \D _{Y'} \overset{^{\mathrm{d}}}{})$.
On dispose des isomorphismes canoniques  :
\begin{equation}
  \label{u_+fonccomflat}
\xymatrix @R=0,3cm {
{F ^\flat _{\mathrm{d}} u _{0+} '  (\NN ')}
\ar@{=}[r]
\ar@{.>}[dd]
&
{\R u _{0*} (\NN '\otimes ^\L _{\D _{X'}} \D _{X '\rightarrow Y'})
\underset{\D _{Y'}}{\overset{\mathrm{d}}{\otimes}}
F ^\flat \D _{Y'}}
\ar[d] ^-{\mathrm{proj}} _-\sim
\\
&
{\R u _{0*} ((\NN '\otimes ^\L _{\D _{X'}} \D _{X '\rightarrow Y'})
\underset{u _0 ^{\text{-}1} \D _{Y'}}{\overset{\mathrm{d}}{\otimes}}
u _0 ^{\text{-}1}  F ^\flat \D _{Y'} )}
\ar[d] _-\sim
\\
{ u _{0+} '  (F ^\flat _{\mathrm{d}} \NN ')}
\ar@{=}[r]
&
{\R u _{0*} (F ^\flat _{\mathrm{d}} (\NN ') \otimes ^\L _{\D _{X'}} \D _{X '\rightarrow Y'}),}
}
\end{equation}
l'indice {\og d \fg} au-dessus du symbole $\otimes$ signifiant que, pour calculer le produit tensoriel,
on choisit la structure {\og droite \fg} de respectivement
$\R u _{0*} (\NN '\otimes ^\L _{\D _{X'}} \D _{X '\rightarrow Y'})$
et $\NN '\otimes ^\L _{\D _{X'}} \D _{X '\rightarrow Y'}$.

\end{vide}

\begin{lemm}
  \label{lemmu_+fonccomflatappl}
 Soient $ \M ' \in D ^- ( \D _{X'} \overset{^{\mathrm{d}}}{})$ et
$\NN ' \in D ^+ ( \D _{X'} \overset{^{\mathrm{d}}}{}, u _0 ^{-1} \D _{Y'} \overset{^{\mathrm{d}}}{})$.
On dispose du diagramme commutatif suivant :
\begin{equation}
  \label{u_+fonccomflatappl}
\xymatrix @R=0,3cm @C=0,5cm {
    { \R u _{0*} \R \mathcal{H} om _{\D _{X'}} ( \M ' ,F ^\flat _{\mathrm{d}} \NN ')}
  \ar[r] ^-{\rho _{2}}
  &
 { \R \mathcal{H} om _{\D _{Y'} } (u _{0+} '  ( \M ' ),u _{0+} '   F ^\flat _{\mathrm{d}} (\NN '))}
  \\
  {\R u _{0*} \left ( \R \mathcal{H} om _{\D _{X'}} ( \M ' ,\NN ')
\otimes  _{u _0 ^{\text{-}1} \D _{Y'} } u _0 ^{\text{-}1} F ^\flat \D _{Y'} \right )}
\ar[u] _-\sim
  &
  {\R \mathcal{H} om _{\D _{Y'} } (u _{0+} '  ( \M ' ),F ^\flat _\mathrm{d}  u _{0+} '   (\NN '))}
  \ar[u] _-\sim ^-{\ref{u_+fonccomflat}}
  \\
  {F ^\flat  \R u _{0*} \R \mathcal{H} om _{\D _{X'}} ( \M ' ,\NN ')}
  \ar[u] ^-{\mathrm{proj}} _-\sim
  \ar[r] ^-{\rho _2}
  &
  {F ^\flat  \R \mathcal{H} om _{\D _{Y'} } (u _{0+} '  ( \M ' ),u _{0+} '   (\NN ')),}
  \ar[u] _-\sim
  }
\end{equation}
où la flèche en haut à droite dérive par fonctorialité de \ref{u_+fonccomflat}.
\end{lemm}
\begin{proof}
Considérons le diagramme ci-dessous: 
\begin{equation}
  \label{u_+fonccomflatappldiag1}
\xymatrix @R=0,3cm @C=0,5cm {
    { \R u _{0*} \R \mathcal{H} om _{\D _{X'}} ( \M ' ,F ^\flat _{\mathrm{d}} \NN ')}
  \ar[r]
  &
 { \R u _{0*} \R \mathcal{H} om _{u _0 ^{\text{-} 1} \D _{Y'}}
( \M '\otimes ^\L _{\D _{X'}} \! \! \!  \D _{X '\rightarrow Y'} ,
F ^\flat _{\mathrm{d}} \NN '\otimes ^\L _{\D _{X'}} \! \! \! \D _{X '\rightarrow Y'})}
  \\
  {\R u _{0*} F ^\flat \R \mathcal{H} om _{\D _{X'}} ( \M ' ,\NN ')}
\ar[u] _-\sim \ar[r]
  &
  { \R u _{0*} F ^\flat  \R \mathcal{H} om _{u _0 ^{\text{-} 1} \D _{Y'}}
( \M '\otimes ^\L _{\D _{X'}} \! \! \!  \D _{X '\rightarrow Y'} ,\NN '\otimes ^\L _{\D _{X'}} \! \! \!  \D _{X '\rightarrow Y'} )}
  \ar[u] _-\sim
  \\
  {F ^\flat  \R u _{0*} \R \mathcal{H} om _{\D _{X'}} ( \M ' ,\NN ')}
  \ar[u] ^-{\mathrm{proj}} _-\sim
  \ar[r]
  &
  {F ^\flat  \R u _{0*} \R \mathcal{H} om _{u _0 ^{\text{-} 1} \D _{Y'}}
( \M '\otimes ^\L _{\D _{X'}} \! \! \!  \D _{X '\rightarrow Y'} ,\NN '\otimes ^\L _{\D _{X'}} \! \! \! \D _{X '\rightarrow Y'} ).}
  \ar[u] ^-{\mathrm{proj}} _-\sim
  }
\end{equation}
\normalsize
Le carré du bas est commutatif par fonctorialité. Pour s'assurer de celle du carré du haut, il suffit par fonctorialité
de l'établir sans $\R u _{0*}$. 
On vérifie que le carré du haut est commutatif sans les symboles $\R$ et $\L$.
On obtient alors sa commutativité en résolvant
$\NN'$ par des $(\D _{X'}, u _0 ^{-1} \D _{Y'} )$-bimodules à droite injectifs et via
une résolution finie de $\D _{X '\rightarrow Y'}$ par des $(\D _{X'}, u _0 ^{-1} \D _{Y'})$-bimodules plats.

Considérons à présent le diagramme ci-après: 
\begin{equation}
  \label{u_+fonccomflatappldiag2}
\xymatrix @R=0,3cm {
 { \R u _{0*} \R \mathcal{H} om _{u _0 ^{\text{-} 1} \D _{Y'}}
( \M '\otimes ^\L _{\D _{X'}} \D _{X '\rightarrow Y'} ,
F ^\flat _{\mathrm{d}} \NN '\otimes ^\L _{\D _{X'}} \D _{X '\rightarrow Y'})}
\ar[r] ^-{\R u _{0*}}
&
{ \R \mathcal{H} om _{\D _{Y'} } (u _{0+} '  ( \M ' ),u _{0+} '   F ^\flat _{\mathrm{d}} (\NN '))}
  \\
  { \R u _{0*} F ^\flat  \R \mathcal{H} om _{u _0 ^{\text{-} 1} \D _{Y'}}
( \M '\otimes ^\L _{\D _{X'}} \D _{X '\rightarrow Y'} ,\NN '\otimes ^\L _{\D _{X'}} \D _{X '\rightarrow Y'} )}
  \ar[u] _-\sim
&
  {\R \mathcal{H} om _{\D _{Y'} } (u _{0+} '  ( \M ' ),F ^\flat _{\mathrm{d}}  u _{0+} '   (\NN '))}
  \ar[u] _-\sim ^-{\ref{u_+fonccomflat}}
  \\
  {F ^\flat  \R u _{0*} \R \mathcal{H} om _{u _0 ^{\text{-} 1} \D _{Y'}}
( \M '\otimes ^\L _{\D _{X'}} \D _{X '\rightarrow Y'} ,\NN '\otimes ^\L _{\D _{X'}} \D _{X '\rightarrow Y'} )}
  \ar[u] ^-{\mathrm{proj}} _-\sim
\ar[r]  ^-{\R u _{0*}}
&
  {F ^\flat  \R \mathcal{H} om _{\D _{Y'} } (u _{0+} '  ( \M ' ),u _{0+} '   (\NN ')).}
  \ar[u] _-\sim
  }
\end{equation}
\normalsize
Pour calculer le terme en bas à gauche, il s'agit de résoudre
$\NN '\otimes ^\L _{\D _{X'}} \D _{X '\rightarrow Y'}$
par des $u _0 ^{- 1} \D _{Y'}$-bimodules à droite injectifs 
et $ \M '\otimes ^\L _{\D _{X'}} \D _{X '\rightarrow Y'}$ par des $u _0 ^{- 1} \D _{Y'}$-modules à droite 
$u _{0*}$-acycliques.
Comme le foncteur $F ^{\flat}$ préserve la $u _0 ^{- 1} \D _{Y'}$-injectivité et la $u _{0*}$-acyclicité, comme
$F ^\flat _{\mathrm{d}} (\NN '\otimes ^\L _{\D _{X'}} \D _{X '\rightarrow Y'} )\riso 
F ^\flat _{\mathrm{d}} \NN '\otimes ^\L _{\D _{X'}} \D _{X '\rightarrow Y'}$,
on vérifie que les autres termes de gauche de \ref{u_+fonccomflatappldiag2} se calculent de même via ces résolutions. 
Par construction des foncteurs horizontaux $\R u _{0*}$, pour vérifier que ce diagramme est commutatif, il suffit de l'établir sans les symboles $\R$ et $\L$, ce qui s'obtient à la main.

Comme \ref{u_+fonccomflatappl} est le composé des diagrammes commutatifs \ref{u_+fonccomflatappldiag1} et \ref{u_+fonccomflatappldiag2}, on termine la preuve.
 \end{proof}

\begin{prop}
\label{prop2-4}
Pour tout $ \M ' \in D ^-( \D _{X'} \overset{^{\mathrm{d}}}{})$, on bénéficie du diagramme commutatif suivant :
\begin{equation}
  \label{2-4}
\xymatrix @R=0,3cm {
{\R u _{0*} \R \mathcal{H} om _{\D _X}
( F ^\flat  \M ', \widetilde{\omega} _X \underset{\B _X}{\otimes} \D _{X\rightarrow Y})}
\ar[r] ^-{\rho _2}
\ar[d] ^-\sim
&
{\R \mathcal{H} om _{\D_Y}
(u _{0+}  (F ^\flat  \M ' ),
u _{0+} (\widetilde{\omega} _X \underset{\B _X}{\otimes} \D _{X \rightarrow Y} ))}
\ar[d] ^-\sim
\\
{\R u _{0*} \R \mathcal{H} om _{\D _X}
( F ^\flat  \M ', F ^\flat _{\mathrm{g}} F ^\flat _{\mathrm{d}} (\widetilde{\omega} _{X'} \underset{\B _{X'}}{\otimes} \D _{X '\rightarrow Y'}))}
\ar[r] ^-{\rho _2}
&
{\R \mathcal{H} om _{\D_Y}
(u _{0+}  (F ^\flat \M ' ),
u _{0+} (F ^\flat _{\mathrm{g}} F ^\flat _{\mathrm{d}} ( \widetilde{\omega} _{X'} \underset{\B _{X'}}{\otimes} \D _{X '\rightarrow Y'} )))}
\\
{\R u _{0*} \R \mathcal{H} om _{\D _{X'}}
( \M ', F ^\flat _{\mathrm{d}} (\widetilde{\omega} _{X'} \underset{\B _{X'}}{\otimes} \D _{X '\rightarrow Y'}))}
\ar[u]_-\sim  ^{F ^\flat} \ar[r] ^-{\rho _2}
&
{\R \mathcal{H} om _{\D _{Y'}} (
u _{0+} '   ( \M ' ),u _{0+} '  (F ^\flat _{\mathrm{d}} ( \widetilde{\omega} _{X'} \underset{\B _{X'}}{\otimes} \D _{X '\rightarrow Y'} ))
)}
\ar[u] _-\sim 
\\
{F ^\flat  \R u _{0*} \R \mathcal{H} om _{\D _{X'}}
( \M ', \widetilde{\omega} _{X'} \underset{\B _{X'}}{\otimes} \D _{X '\rightarrow Y'})}
\ar[r] ^-{\rho _2}
\ar[u]  _-\sim
&
{F ^\flat \R \mathcal{H} om _{\D _{Y'}} (
u _{0+} '   ( \M ' ), u _{0+} '  ( \widetilde{\omega} _{X'} \underset{\B _{X'}}{\otimes} \D _{X '\rightarrow Y'} ) ).}
\ar[u] _-\sim
}
\end{equation}
\normalsize
\end{prop}
\begin{proof}
  Par \ref{lemmpre24}, on obtient la commutativité du carré du milieu. Celle du haut se vérifie
par fonctorialité tandis que celle du bas résulte de \ref{u_+fonccomflatappl}.
 \end{proof}

\subsection{Troisième partie de la construction de l'isomorphisme de dualité et sa compatibilité à Frobenius}
\begin{vide}\label{para1.4.1}
Posons $\L u _0 ^{\prime *}: = u _0 ^{\prime !} [-d _{X'/Y'}]$.
Pour tous $\FF ' \in D ^- _{\mathrm{qc}} (\overset{^{\mathrm{g}}}{}  \D _{Y'}     )$,
$\M '\in D ^- (\D _{X'} \overset{^{\mathrm{d}}}{})$,
de manière analogue à \cite[1.2.27]{caro_surcoherent},
on bénéficie de l'isomorphisme canonique
$u ' _{0+} (\M '\otimes ^\L _{\B _{X'} } \L u _0 ^{\prime *}(\FF '))$
\linebreak[1]
$
\riso
u '_{0+} (\M ' ) \otimes ^\L  _{\B _{Y'}} \FF'$ via le composé :
\begin{gather}
\notag
\R u _{0*} (  (\M '\otimes ^\L _{\B _{X'}} \L u _0 ^{\prime *} (\FF ' ))
  \otimes _{\D _{X'}} ^\L \! \! \!  \D _{X '\rightarrow Y'})
\riso
\R u _{0*} (  \M '\otimes _{\D _{X'}} ^\L (\L u _0 ^{\prime *} (\FF ')
  \otimes _{\B _{X'}} ^\L \! \! \!\D _{X '\rightarrow Y'}))
\\
 \label{u+!proj}
\underset{\gamma _{\FF'}}{\liso}
\R u _{0*} ( \M '\otimes ^\L _{\D _{X'}}
 ( \D _{X' \rightarrow Y'} \otimes  _{ u ^{-1} _0 \B _{Y'} } ^\L u ^{-1} _0  \FF '))
\underset{\mathrm{proj}}{\liso}
u '_{0+}   (\M ') \otimes _{\B _{Y'}} ^\L \FF',
\end{gather}
où le premier isomorphisme se vérifie de façon analogue à \cite[1.2.23]{caro_surcoherent}
et où, modulo les isomorphismes canoniques $\L u _0 ^{\prime *} (\FF ')
  \otimes _{\B _{X'}} ^\L \D _{X '\rightarrow Y'}
  \riso \L u _0 ^{\prime *} (\FF '
  \otimes _{\B _{Y'}} ^\L \D _{Y'})$ 
  et
  $ \L u _0 ^{\prime *} (\D _{Y'}\otimes _{\B _{Y'}} ^\L \FF ' )\riso   
    \D _{X' \rightarrow Y'} \otimes  _{ u ^{-1} _0 \B _{Y'} } ^\L u ^{-1} _0  \FF '$,
   $\gamma _{\FF'}$ désigne l'image par $\R u _{0*} (  \M '\otimes _{\D _{X'}} ^\L \L u _0 ^{\prime *} (-))$
de l'isomorphisme de transposition d'une résolution plate de $\FF'$.
\end{vide}

\begin{vide}
Pour tout $ \M ' \in D ^- ( \D _{X'} \overset{^{\mathrm{d}}}{})$,
nous prouvons dans cette sous-section la commutation à Frobenius de
  l'isomorphisme canonique :
$$\R \mathcal{H} om _{\D _{Y'}}
(u _{0+} '   ( \M ' ),
 u _{0+} '  (\widetilde{\omega} _{X'} \otimes _{\B _{X'}} \D _{X '\rightarrow Y'} ) )
\riso
 \R \mathcal{H} om _{\D _{Y'}} (u _{0+} '   ( \M ' ), u _{0+} '  (\widetilde{\omega} _{X'} ) \otimes _{\B _{Y'}} \D _{Y'})$$
 \normalsize
 (i.e., le diagramme \ref{4-7} est commutatif).
Pour cela, nous établirons que l'isomorphisme
$u _{0+} '  (\widetilde{\omega} _{X'} \otimes _{\B _{X'}} \D _{X '\rightarrow Y'} )
\riso
u _{0+} '  (\widetilde{\omega} _{X'} ) \otimes _{\B _{Y'}} \D _{Y'}$
commute à Frobenius (voir \ref{47u*}).
Celui-ci est le composé de trois isomorphismes (voir \ref{u+!proj}).
Via \ref{4-5u*}, \ref{56} et \ref{67}, nous allons successivement vérifier
leur commutation à Frobenius.
\end{vide}

\begin{vide}
Afin de valider la commutativité de \ref{4-5u*}, nous aurons besoin de celle de \ref{4-5diag4}.
À cette fin, établissons d'abord celle des diagrammes \ref{4-5prediag2} et \ref{4-5prediag5} ci-dessous.

Soient $ \M '$ un $\D _{X'}$-module à droite et $\E '$ un $\D _{X'}$-module à gauche.
Par passage de gauche à droite du rectangle du haut de \cite[1.4.15.1]{caro_comparaison}, on obtient
la commutativité du diagramme suivant :
\begin{equation}
\notag  
  \xymatrix @R=0,3cm {
  {F ^\flat  \M ' \otimes _{\D _X} \D _X} \ar[r] _-\sim
&
{F ^\flat  \M ' \otimes _{\D _X}  F ^* F ^\flat  \D _{X'}}
\ar[d] _-\sim
\\
 {F ^\flat  \M '}
  \ar[u] _-\sim
&
    { \M ' \otimes _{\D _{X'}} F ^\flat  \D _{X'}.}
\ar[l] _-\sim
}
\end{equation}

En lui appliquant
$-\otimes _{\B _X} F ^* \E '$,
il arrive le carré du haut de :
\begin{equation}
  \label{4-5prediag2}
  \xymatrix @R=0,3cm {
  {F ^\flat  \M ' \otimes _{\D _X} \D _X \otimes _{\B _X} F ^* \E '}
\ar[rr] _-\sim
&&
 {F ^\flat  \M ' \otimes _{\D _X} F^* F ^\flat \D _{X'} \otimes _{\B _X} F ^* \E '}
\ar[d] _-\sim
\\
{F ^\flat  \M ' \otimes _{\B _X} F ^* \E '}
\ar[u] _-\sim \ar[d] _-\sim
&&
 { \M ' \otimes _{\D _{X'}} F ^\flat \D _{X'} \otimes _{\B _X} F ^* \E ' }
\ar[d] _-\sim \ar[ll] _-\sim
\\
 {F ^\flat ( \M ' \otimes _{\B _{X'}} \E ' )}
\ar[r] _-\sim
&
 {F ^\flat (   \M ' \otimes _{\D _{X'}}  \D _{X'} \otimes _{\B _{X'}} \E ' ) }
&
 { \M ' \otimes _{\D _{X'}} F ^\flat ( \D _{X'} \otimes _{\B _{X'}} \E ' ) .}
\ar[l] _-\sim
 }
\end{equation}
\normalsize
La commutativité du rectangle du bas découle d'un calcul local immédiat.
Le diagramme \ref{4-5prediag2} est donc commutatif.

Enfin, pour tous $\D _{X'}$-modules à gauche $\E '$ et $\FF '$, on dispose du diagramme commutatif :
\begin{equation}
  \label{4-5prediag3}
  \xymatrix @R=0,3cm {
  {F ^* \E ' \otimes _{\B _X}   F ^* F ^\flat \D _{X'} \otimes _{\D _X} F ^* \FF ' }
\ar[d] _-\sim
&&
  {F ^* \E ' \otimes _{\B _X}  \D _X \otimes _{\D _X} F ^* \FF ' }
\ar[ll] _-\sim
\\
{F ^* \E ' \otimes _{\B _X}   F ^*  \D _{X'} \otimes _{\D _{X'}}  \FF ' }
&&
  {F ^* \E ' \otimes _{\B _X}   F^* \FF ' } \ar[ll] _-\sim  \ar[u] _-\sim
\\
  {F ^* (\E ' \otimes _{\B _{X'}}  \D _{X'} ) \otimes _{\D _{X'}}  \FF ' }
\ar[u] _-\sim
&
  {F ^* (\E ' \otimes _{\B _{X'}}  \D _{X'}  \otimes _{\D _{X'}}  \FF ' ) }
\ar[l] _-\sim
&
  {F ^* (\E ' \otimes _{\B _{X'}}  \FF ' ).}
\ar[l] _-\sim \ar[u] _-\sim
 }
\end{equation}
\normalsize
En effet, par application de $F ^* \E ' \otimes _{\B _X} -$ au rectangle du haut
de \cite[1.4.15.1]{caro_comparaison} (utilisé pour $\FF '$), on parvient au rectangle
du haut de \ref{4-5prediag3}.
De plus, soit $e ', f', a$ des sections locales respectives de $\E'$, $\FF'$, $\B_X$.
Le morphisme composé $F ^* (\E ' \otimes _{\B _{X'}}  \FF ' )
\rightarrow F ^* \E ' \otimes _{\B _X}   F ^*  \D _{X'} \otimes _{\D _{X'}}  \FF ' $
du rectangle du bas passant par le haut envoie $a \otimes (e'\otimes f')$ sur
$(a \otimes e')\otimes (1 \otimes f')$ puis $(a \otimes e')\otimes (1 \otimes 1) \otimes f'$.
Via le chemin du bas, on calcule
$a \otimes (e'\otimes f')\mapsto  a \otimes (e'\otimes 1 \otimes f') \mapsto
a \otimes (e'\otimes 1 ) \otimes f' \mapsto ( a \otimes e') \otimes (1 \otimes 1 ) \otimes f' $.
D'où la commutativité du carré du bas.

Dans le diagramme ci-dessous,
on vérifie la commutativité du carré du bas par fonctorialité:
\begin{equation}
  \label{4-5prediag4}
  \xymatrix @R=0,3cm {
  {F ^* \E ' \otimes _{\B _X}   F ^* F ^\flat \D _{X'} \otimes _{\D _X} F ^* \FF ' }
\ar@{=}[r]
&
{F ^* \E ' \otimes _{\B _X}   F ^* F ^\flat \D _{X'} \otimes _{\D _X} F ^* \FF ' }
\ar[d] _-\sim
\\
{F ^* F ^\flat (\E ' \otimes _{\B _{X'}} \D _{X'})
  \otimes _{\D _X} F ^* \FF '}
\ar[r] _-\sim
\ar[u] _-\sim
\ar[d] _-\sim
&
{F ^\flat (F ^* \E ' \otimes _{\B _{X}}F ^*  \D _{X'})
  \otimes _{\D _X} F ^* \FF '}
\ar[d] _-\sim
\\
  {F ^* (\E ' \otimes _{\B _{X'}}  \D _{X'} ) \otimes _{\D _{X'}}  \FF ' }
\ar[r] _-\sim
&
{F ^* \E ' \otimes _{\B _X}   F ^*  \D _{X'} \otimes _{\D _{X'}}  \FF ' .}
 }
\end{equation}
Comme les foncteurs $F ^*$ et $F ^\flat$ commutent canoniquement, 
le carré carré du haut est commutatif par définition.
Or, le composé de droite de \ref{4-5prediag4},
$F ^* \E ' \otimes _{\B _X}   F ^* F ^\flat \D _{X'} \otimes _{\D _X} F ^* \FF '
\riso
F ^* \E ' \otimes _{\B _X}   F ^*  \D _{X'} \otimes _{\D _{X'}}  \FF ' $,
est aussi la flèche de gauche du haut de \ref{4-5prediag3}.
En composant les contours de \ref{4-5prediag4} et \ref{4-5prediag3}, on obtient alors le diagramme commutatif suivant dont
la partie droite est celle de \ref{4-5prediag3} et la partie gauche est celle de \ref{4-5prediag4} :
\begin{equation}
  \label{4-5prediag5}
  \xymatrix @R=0,3cm {
  {F ^* \E ' \otimes _{\B _X}   F ^* F ^\flat \D _{X'} \otimes _{\D _X} F ^* \FF ' }
&
 {F ^* \E ' \otimes _{\B _X}  \D _X \otimes _{\D _X} F ^* \FF ' }
\ar[l] _-\sim
\\
{F ^* F ^\flat (\E ' \otimes _{\B _{X'}} \D _{X'})
  \otimes _{\D _X} F ^* \FF '}
\ar[u] _-\sim
\ar[d] _-\sim
&
{F ^* \E ' \otimes _{\B _X}   F^* \FF ' } \ar[u] _-\sim
\\
  {F ^* (\E ' \otimes _{\B _{X'}}  \D _{X'} ) \otimes _{\D _{X'}}  \FF ' }
&
  {F ^* (\E ' \otimes _{\B _{X'}}  \FF ' ).}
\ar[l] _-\sim \ar[u] _-\sim
 }
\end{equation}

\end{vide}

\begin{vide}
Soient $\E '$ et $\FF '$ deux $\D _{X'}$-modules à gauche et
$ \M '$ un $\D_{X'}$-module à droite. Considérons le diagramme ci-dessous :
\begin{equation}
  \label{4-5diag1}
  \xymatrix @R=0,3cm  {
  {(F ^\flat  \M ' \otimes _{\B _X} F ^* \E ')
  \otimes _{\D _X} F ^* \FF '}
  \ar[r]_-\sim \ar[ddd] _-\sim
  &
  {F ^\flat  \M ' \otimes _{\D _X} (\D _X \otimes _{\B _X} F ^* \E ')
  \otimes _{\D _X} F ^* \FF '}
  \ar[d] _-\sim
  \\
&
  {F ^\flat  \M ' \otimes _{\D _X} (F ^* F ^\flat \D _{X'} \otimes _{\B _X} F ^* \E ')
  \otimes _{\D _X} F ^* \FF '}
  \ar[d] _-\sim
  \\
&
  {F ^\flat  \M ' \otimes _{\D _X} F ^* F ^\flat (\D _{X'} \otimes _{\B _{X'}} \E ')
  \otimes _{\D _X} F ^* \FF '}
  \ar[d] _-\sim
  \\
  {F ^\flat (  \M ' \otimes _{\B _{X'}} \E ')
  \otimes _{\D _X} F ^* \FF '}
  \ar[r]_-\sim \ar[d] _-\sim
  &
  { \M ' \otimes _{\D _{X'}} F ^\flat (\D _{X'} \otimes _{\B _{X'}} \E ')
  \otimes _{\D _X} F ^* \FF '}
  \ar[d] _-\sim
  \\
  {(  \M ' \otimes _{\B _{X'}} \E ')
  \otimes _{\D _{X'}} \FF '}
  \ar[r] _-\sim
  &
  { \M ' \otimes _{\D _{X'}} (\D _{X'} \otimes _{\B _{X'}} \E ')
  \otimes _{\D _{X'}} \FF '.}
  }
\end{equation}
En appliquant le foncteur $-\otimes _{\D _X} F ^* \FF '$ au contour de \ref{4-5prediag2}, on obtient le rectangle du haut de \ref{4-5diag1}.
De plus, modulo
$F ^\flat (  \M ' \otimes _{\D _{X'}} \D _{X'} \otimes _{\B _{X'}} \E ')
  \otimes _{\D _X} F ^* \FF ' \liso
   \M ' \otimes _{\D _{X'}} F ^\flat (\D _{X'} \otimes _{\B _{X'}} \E ')
  \otimes _{\D _X} F ^* \FF '$,
on obtient par fonctorialité la commutativité du carré du bas de
  \ref{4-5diag1}. Le diagramme \ref{4-5diag1} est donc commutatif.

En appliquant
$F ^\flat  \M ' \otimes _{\D _X} (-)
  \otimes _{\D _X} F ^* \FF '$ au diagramme signifiant que l'isomorphisme de transposition
  de $\E '$, noté $\gamma _{\E'}$, (voir \cite[1.3.1]{Be2}) est compatible à Frobenius
(\cite[2.1.8.1]{caro_comparaison}), on obtient le rectangle
  du haut du diagramme :
\begin{equation}
\label{4-5diag2}
  \xymatrix @R=0,3cm  {
  {F ^\flat  \M ' \otimes _{\D _X} (\D _X \otimes _{\B _X} F ^* \E ')
  \otimes _{\D _X} F ^* \FF '}
  \ar[r] _-\sim ^{\gamma _{F ^* \E'}} \ar[d] _-\sim
  &
  {F ^\flat  \M ' \otimes _{\D _X} (F ^* \E ' \otimes _{\B _X} \D _X )
  \otimes _{\D _X} F ^* \FF '}\ar[d] _-\sim
\\
  {F ^\flat  \M ' \otimes _{\D _X} (F ^* F ^\flat \D _{X'} \otimes _{\B _X} F ^* \E ')
  \otimes _{\D _X} F ^* \FF '}\ar[d] _-\sim
  &
  {F ^\flat  \M ' \otimes _{\D _X} (F ^* \E ' \otimes _{\B _X} F ^* F ^\flat \D _{X'})
  \otimes _{\D _X} F ^* \FF '}\ar[d] _-\sim
  \\
  {F ^\flat  \M ' \otimes _{\D _X} F ^* F ^\flat (\D _{X'} \otimes _{\B _{X'}} \E ')
  \otimes _{\D _X} F ^* \FF '}
  \ar[r] _-\sim ^{\gamma _{\E'}} \ar[d] _-\sim
  &
  {F ^\flat  \M ' \otimes _{\D _X} F ^* F ^\flat (\E ' \otimes _{\B _{X'}} \D _{X'})
  \otimes _{\D _X} F ^* \FF '}\ar[d] _-\sim
  \\
  { \M ' \otimes _{\D _{X'}} F ^\flat (\D _{X'} \otimes _{\B _{X'}} \E ')
  \otimes _{\D _X} F ^* \FF '}\ar[d] _-\sim
  &
  {F ^\flat  \M ' \otimes _{\D _X} F ^* (\E ' \otimes _{\B _{X'}} \D _{X'})
  \otimes _{\D _{X'}} \FF '}\ar[d] _-\sim
  \\
  { \M ' \otimes _{\D _{X'}} (\D _{X'} \otimes _{\B _{X'}} \E ')
  \otimes _{\D _{X'}} \FF '}
  \ar[r] _-\sim ^{\gamma _{\E'}}
  &
  { \M ' \otimes _{\D _{X'}} (\E ' \otimes _{\B _{X'}} \D _{X'} )
  \otimes _{\D _{X'}} \FF '.}
  }
\end{equation}
\normalsize
Le rectangle du bas de \ref{4-5diag2} est commutatif par fonctorialité. D'où la commutativité
de \ref{4-5diag2}.

En appliquant $F ^\flat  \M ' \otimes _{\D _X}-$ au diagramme commutatif \ref{4-5prediag5},
on obtient le rectangle du haut de :
\begin{equation}
\label{4-5diag3}
  \xymatrix @R=0,3cm  {
  {F ^\flat  \M ' \otimes _{\D _X} (F ^* \E ' \otimes _{\B _X} \D _X )
  \otimes _{\D _X} F ^* \FF '}\ar[d] _-\sim
  &
  {F ^\flat  \M ' \otimes _{\D _X} (F ^* \E ' \otimes _{\B _X} F ^* \FF ')}
  \ar[l] _-\sim
\\
  {F ^\flat  \M ' \otimes _{\D _X} (F ^* \E ' \otimes _{\B _X} F ^* F ^\flat \D _{X'})
  \otimes _{\D _X} F ^* \FF '}
  \\
  {F ^\flat  \M ' \otimes _{\D _X} F ^* F ^\flat (\E ' \otimes _{\B _{X'}} \D _{X'})
  \otimes _{\D _X} F ^* \FF '}
\ar[d] _-\sim \ar[u] _-\sim
  \\
  {F ^\flat  \M ' \otimes _{\D _X} F ^* (\E ' \otimes _{\B _{X'}} \D _{X'})
  \otimes _{\D _{X'}} \FF '}   \ar[d] _-\sim
  &
  {F ^\flat  \M ' \otimes _{\D _X} F ^* (\E ' \otimes _{\B _{X'}} \FF ' )}   \ar[d] _-\sim \ar[l] _-\sim
\ar[uuu] _-\sim
  \\
  { \M ' \otimes _{\D _{X'}} (\E ' \otimes _{\B _{X'}} \D _{X'} )
  \otimes _{\D _{X'}} \FF '}
  &
  { \M ' \otimes _{\D _{X'}} (\E ' \otimes _{\B _{X'}} \FF ').} \ar[l] _-\sim
  }
\end{equation}
Puisque le carré du bas est commutatif, il en dérive celle de \ref{4-5diag3}.

En composant \ref{4-5diag1}, \ref{4-5diag2} et \ref{4-5diag3}, on obtient le diagramme commutatif :
\begin{equation}
  \label{4-5diag4}
  \xymatrix @R=0,3cm {
{(F ^\flat  \M ' \otimes _{\B _X} F ^* \E ')
  \otimes _{\D _X} F ^* \FF '}
  \ar[r] _-\sim\ar[d] _-\sim
  &
{F ^\flat  \M ' \otimes _{\D _X} (F ^* \E ' \otimes _{\B _X} F ^* \FF ')}
  \ar[d] _-\sim
\\
{(  \M ' \otimes _{\B _{X'}} \E ')
  \otimes _{\D _{X'}} \FF '}
\ar[r] _-\sim
  &
  { \M ' \otimes _{\D _{X'}} (\E ' \otimes _{\B _{X'}} \FF ').}
  }
\end{equation}

\end{vide}

\begin{vide}
On bénéficie du diagramme commutatif suivant :
\begin{equation}
  \label{4-5u*}
  \xymatrix @R=0,3cm {
  {\R u _{0*} ( (\widetilde{\omega} _X \otimes _{\B _X} \D _{X \rightarrow Y})
  \otimes ^\L _{\D _X} \D _{X\rightarrow Y})} \ar[r] _-\sim  \ar[d] _-\sim
  &
  {\R u _{0*} ( \widetilde{\omega} _X \otimes _{\D _X} ^\L
  (\D _{X \rightarrow Y} \otimes _{\B _X} \D _{X\rightarrow Y}))} \ar[d] _-\sim
  \\
    {\R u _{0*} (  ( F ^\flat \widetilde{\omega} _{X'} \otimes _{\B _X} F ^* F ^\flat \D _{X '\rightarrow Y'})
  \otimes _{\D _X} ^\L F ^* F ^\flat \D _{X '\rightarrow Y'})} \ar[r] _-\sim  \ar[d] _-\sim
  &
  {\R u _{0*} (F ^\flat \widetilde{\omega} _{X'} \otimes _{\D _X} ^\L
  (F ^* F ^\flat \D _{X '\rightarrow Y'} \otimes _{\B _X} F ^* F ^\flat \D _{X '\rightarrow Y'}))} \ar[d] _-\sim
  \\
    {\R u _{0*} (  (\widetilde{\omega} _{X'} \otimes _{\B _{X'}} F ^\flat \D _{X '\rightarrow Y'})
  \otimes _{\D _{X'}} ^\L F ^\flat \D _{X '\rightarrow Y'})}
  \ar[r] _-\sim \ar[d] _-\sim
  &
  {\R u _{0*} (\widetilde{\omega} _{X'} \otimes _{\D _{X'}} ^\L
  (F ^\flat \D _{X '\rightarrow Y'} \otimes _{\B _{X'}} F ^\flat \D _{X '\rightarrow Y'}))}
\ar[d] _-\sim
  \\
    { \R u _{0*} F ^\flat  _\mathrm{g} F ^\flat _\mathrm{d}  (  (\widetilde{\omega} _{X'} \otimes _{\B _{X'}} \D _{X '\rightarrow Y'})
  \otimes _{\D _{X'}} ^\L \D _{X '\rightarrow Y'})}
  \ar[r] _-\sim
  &
  { \R u _{0*} F ^\flat  _\mathrm{g} F ^\flat _\mathrm{d}  (\widetilde{\omega} _{X'} \otimes _{\D _{X'}} ^\L
  (\D _{X '\rightarrow Y'} \otimes _{\B _{X'}} \D _{X '\rightarrow Y'}))}
 \\
    {F ^\flat _{\mathrm{g}} F ^\flat  _{\mathrm{d}} \R u _{0*} (  (\widetilde{\omega} _{X'} \otimes _{\B _{X'}} \D _{X '\rightarrow Y'})
  \otimes _{\D _{X'}} ^\L \D _{X '\rightarrow Y'})} \ar[r] _-\sim  \ar[u] _-\sim  ^-{\mathrm{proj}}
  &
  {F ^\flat _{\mathrm{g}} F ^\flat  _{\mathrm{d}} \R u _{0*} (\widetilde{\omega} _{X'} \otimes _{\D _{X'}} ^\L
  (\D _{X '\rightarrow Y'} \otimes _{\B _X} \D _{X '\rightarrow Y'})).}
\ar[u] _-\sim  ^-{\mathrm{proj}}
}
\end{equation}
\normalsize
En effet, on obtient par fonctorialité la commutativité des carrés du haut et du bas. 
Pour vérifier celle des deux autres carrés, il suffit de l'établir sans le foncteur $\R u _{0*}$.
En résolvant $\widetilde{\omega} _{X'}$ par des $\D _{X'}$-modules à droite plats,
$F ^\flat \D _{X '\rightarrow Y'}$
par des $(\D _{X'}, u _0 ^{-1} \D _{Y} )$-bimodules plats,
la commutativité du deuxième carré en partant du haut découle alors de celle de
\ref{4-5diag4}.
En résolvant $\widetilde{\omega} _{X'}$ par des $\D _{X'}$-modules à droite plats,
$\D _{X '\rightarrow Y'}$
par des $(\D _{X'}, u _0 ^{-1} \D _{Y'} )$-bimodules plats,
on se ramène à vérifier la commutativité
du deuxième carré du bas sans le symbole $\L$. 
Cette commutativité se vérifie alors par un calcul (voir le calcul des flèches horizontales donné dans le début de la preuve de \ref{47lemm1}).
\end{vide}

\begin{vide}\label{notagama}
Considérons le diagrame suivant : 
\begin{equation}
  \label{56prediag1}
  \xymatrix @R=0,25cm {
{\D _Y\otimes ^{\mathrm{gg}} _{\B _Y} \D _Y}
\ar[d] _-\sim
&
{\D _Y \otimes ^{\mathrm{dg}} _{\B _Y} \D _Y }
\ar[d] _-\sim  \ar[l] _-\sim  ^-{\gamma _{\D _Y }}
\\
{F ^* F ^\flat \D _{Y'} \otimes ^{\mathrm{gg}} _{\B _Y} \D _Y} \ar[d] _-\sim  &
{\D _Y \otimes ^{\mathrm{dg}} _{\B _Y} F ^* F ^\flat \D _{Y'}  }
\ar[d] _-\sim
\ar[l] _-\sim  ^-{\gamma _{F ^* F ^\flat \D _{Y'} }}
\\
{F ^* F ^\flat \D _{Y'} \otimes ^{\mathrm{gg}} _{\B _Y} F ^* F ^\flat \D _{Y'} } \ar[d] _-\sim
&
{F ^* F ^\flat \D _{Y'} \otimes ^{\mathrm{dg}} _{\B _Y} F ^* F ^\flat \D _{Y'}  } \ar[d] _-\sim
\\
{F ^* F ^\flat _{\mathrm{d}} (F ^\flat \D _{Y'} \otimes ^{\mathrm{gg}} _{\B _{Y'}} \D _{Y'} )} \ar[d] _-\sim
&
{F ^* F ^\flat _{\mathrm{g}} (\D _{Y'} \otimes ^{\mathrm{dg}} _{\B _{Y'}} F ^\flat \D _{Y'} )  }
\ar[d] _-\sim  \ar[l] _-\sim  ^-{\gamma _{F ^\flat \D _{Y'} }}
\\
{F ^* F ^\flat _{\mathrm{d}}  F ^\flat _{\mathrm{g}} (\D _{Y'} \otimes ^{\mathrm{gg}} _{\B _{Y'}} \D _{Y'} )}
&
{F ^* F ^\flat _{\mathrm{g}}  F ^\flat _{\mathrm{d}} (\D _{Y'} \otimes ^{\mathrm{dg}} _{\B _{Y'}} \D _{Y'} )  ,}
\ar[l] _-\sim  ^-{\gamma _{\D _{Y'} }}
}
\end{equation}
où $\D _{Y'} \otimes ^{\mathrm{gg}} _{\B _{Y'}} \D _{Y'}$ (resp. $\D _{Y'} \otimes ^{\mathrm{dg}} _{\B _{Y'}} \D _{Y'}$)
signifie que, pour calculer le produit tensoriel,
on utilise la structure de $\D _{Y'}$-module à gauche (resp. à droite) du terme $\D _{Y'} $ de gauche
et la structure de $\D _{Y'}$-module à gauche du terme $\D _{Y'} $ de droite ; de même pour les termes
analogues. On remarque que
$\D _{Y'} \otimes ^{\mathrm{gg}} _{\B _{Y'}} \D _{Y'}$ (resp. $\D _{Y'} \otimes ^{\mathrm{dg}} _{\B _{Y'}} \D _{Y'}$)
est un $( \D _{Y'} \overset{^{\mathrm{d}}}{}, \overset{^{\mathrm{g}}}{} \D _{Y'} , \D _{Y'} \overset{^{\mathrm{d}}}{})$-trimodule
(resp. un $( \overset{^{\mathrm{g}}}{} \D _{Y'} ,\D _{Y'} \overset{^{\mathrm{d}}}{} , \D _{Y'} \overset{^{\mathrm{d}}}{})$-trimodule),
la structure {\og du milieu \fg} étant la structure {\og produit tensoriel \fg} et les deux autres,
structure gauche et structure droite, sont induites par fonctorialité.

Le diagramme \ref{56prediag1} est commutatif. 
En effet, on établit par fonctorialité la commutativité du carré du haut de \ref{56prediag1}.
De plus, le diagramme \cite[2.1.8.1]{caro_comparaison} signifiant la compatibilité à Frobenius de l'isomorphisme
  de transposition du $\D _{Y'}$-module à gauche
  $F ^\flat \D _{Y'}$, noté $\gamma _{F ^\flat \D _{Y'}}$,
correspond au rectangle (du milieu) du diagramme \ref{56prediag1}.
Enfin, pour vérifier la commutativité du carré du bas, par fonctorialité et comme cela est local, il suffit de le calculer dans le cas relevable et sans les foncteurs
{\og $F ^* F ^\flat _{\mathrm{d}}$\fg} (pour les termes de gauche) ou {\og $F ^* F ^\flat _{\mathrm{g}}$\fg}
(pour les termes de droite): Soient $\phi$, $P'$ des sections locales respectives de
$\mathcal{H} om _{\B _{Y'}} (\B _Y, \B _{Y'})$, $\D _{Y'}$.
Par $\D _{Y'}$-linéarité, il suffit de vérifier que $1 \otimes (P ' \otimes \phi)$ s'envoie via
les deux flèches
$\D _{Y'} \otimes ^{\mathrm{dg}} _{\B _{Y'}} F ^\flat \D _{Y'} \riso
F ^\flat _{\mathrm{g}} (\D _{Y'} \otimes ^{\mathrm{gg}} _{\B _{Y'}} \D _{Y'} )$ sur le même élément.
Or, par le chemin du haut,
$1 \otimes (P ' \otimes \phi) \mapsto
(P ' \otimes \phi) \otimes 1 \mapsto (P ' \otimes 1 )\otimes  \phi$
(le premier calcul est la caractérisation de $\gamma _{F ^\flat \D _{Y'}}$ donnée dans \cite[1.3.1]{Be2}).
Le second chemin donne
$1 \otimes (P ' \otimes \phi) \mapsto
(1 \otimes P ' )\otimes  \phi \mapsto (P ' \otimes 1) \otimes \phi $.
D'où la commutativité de \ref{56prediag1}.

Déduisons-en à présent la commutativité du diagramme ci-dessous : 
\begin{equation}
  \label{56}
\xymatrix @R=0,25cm @C=0,5cm {
{\R u _{0*} ( \widetilde{\omega} _{X} \otimes ^\L _{\D _X}
( \D _{X \rightarrow Y} \otimes  _{\B _X} \D _{X \rightarrow Y}))}
&
{\R u _{0*} ( \widetilde{\omega} _{X} \otimes ^\L _{\D _X}
( \D _{X \rightarrow Y} \otimes  _{u ^{-1} _0 \B _Y} u ^{-1} _0 \D _Y))}
\ar[l] _-\sim
\\
{\R u _{0*} ( \widetilde{\omega} _{X} \otimes ^\L _{\D _X}
u ^* _0( \D _Y \otimes ^{\mathrm{gg}} _{\B _Y} \D _Y))}
\ar[d] _-\sim  \ar[u] _-\sim
&
{\R u _{0*} ( \widetilde{\omega} _{X} \otimes ^\L _{\D _X}
u ^* _0( \D _Y \otimes ^{\mathrm{dg}} _{\B _Y} \D _Y))}
\ar[u] _-\sim  \ar[d] _-\sim  \ar[l] _-\sim  ^-{\gamma _{\D _Y }}
\\
{\R u _{0*} ( \widetilde{\omega} _{X} \otimes ^\L _{\D _X}
u ^* _0 F ^* F ^\flat _{\mathrm{d}} F ^\flat _{\mathrm{g}}  ( \D _{Y'} \otimes ^{\mathrm{gg}} _{\B _{Y'}} \D _{Y'}))}
\ar[d] _-\sim
&
{\R u _{0*} ( \widetilde{\omega} _{X} \otimes ^\L _{\D _X}
u ^* _0 F ^* F ^\flat _{\mathrm{g}} F ^\flat _{\mathrm{d}}  ( \D _{Y'} \otimes ^{\mathrm{dg}} _{\B _{Y'}} \D _{Y'} ))}
\ar[l] _-\sim  ^-{\gamma _{\D _{Y'} }}
\ar[d] _-\sim
\\
{\R u _{0*} F ^\flat _{\mathrm{d}} F ^\flat _{\mathrm{g}}  ( \widetilde{\omega} _{X} \otimes ^\L _{\D _X}
u ^* _0 F ^*  ( \D _{Y'} \otimes ^{\mathrm{gg}} _{\B _{Y'}} \D _{Y'}))}
&
{\R u _{0*} F ^\flat _{\mathrm{g}} F ^\flat _{\mathrm{d}} ( \widetilde{\omega} _{X} \otimes ^\L _{\D _X}
u ^* _0 F ^*   ( \D _{Y'} \otimes ^{\mathrm{dg}} _{\B _{Y'}} \D _{Y'} ))}
\ar[l] _-\sim  ^-{\gamma _{\D _{Y'} }}
\\
{F ^\flat _{\mathrm{d}} F ^\flat _{\mathrm{g}}  \R u _{0*} ( \widetilde{\omega} _{X} \otimes ^\L _{\D _X}
u ^* _0 F ^*  ( \D _{Y'} \otimes ^{\mathrm{gg}} _{\B _{Y'}} \D _{Y'}))}
\ar[u] _-\sim  ^-{\mathrm{proj}} \ar[d] _-\sim
&
{F ^\flat _{\mathrm{g}} F ^\flat _{\mathrm{d}} \R u _{0*} ( \widetilde{\omega} _{X} \otimes ^\L _{\D _X}
u ^* _0 F ^*  ( \D _{Y'} \otimes ^{\mathrm{dg}} _{\B _{Y'}} \D _{Y'} ))}
\ar[u] _-\sim  ^-{\mathrm{proj}} \ar[d] _-\sim
\ar[l] _-\sim  ^-{\gamma _{\D _{Y'} }}
\\
{F ^\flat _{\mathrm{d}} F ^\flat _{\mathrm{g}}  \R u _{0*} ( \widetilde{\omega} _{X'} \otimes ^\L _{\D _{X'}}
u ^{\prime *} _0  ( \D _{Y'} \otimes ^{\mathrm{gg}} _{\B _{Y'}} \D _{Y'}))} \ar[d] _-\sim
&
{F ^\flat _{\mathrm{g}} F ^\flat _{\mathrm{d}}  \R u _{0*} ( \widetilde{\omega} _{X'} \otimes ^\L _{\D _{X'}}
u ^{\prime *} _0   ( \D _{Y'} \otimes ^{\mathrm{dg}} _{\B _{Y'}} \D _{Y'} ))}
\ar[l] _-\sim  ^-{\gamma _{\D _{Y'} }}
\ar[d] _-\sim
\\
{F ^\flat _{\mathrm{d}} F ^\flat _{\mathrm{g}}  \R u _{0*} ( \widetilde{\omega} _{X'} \otimes ^\L _{\D _{X'}}
\! \! \!  (\D _{X '\rightarrow Y'} \otimes  _{\B _{X'}} \! \! \!  \D _{X '\rightarrow Y'}))}
&
{F ^\flat _{\mathrm{g}} F ^\flat _{\mathrm{d}} \R u _{0*} ( \widetilde{\omega} _{X'} \otimes ^\L _{\D _{X'}}
\! \! \!  ( \D _{X '\rightarrow Y'} \otimes  _{ u ^{-1} _0 \B _{Y'}} u ^{-1} _0  \! \! \!  \D _{Y'} )).}
\ar[l] _-\sim
}
\end{equation}
\normalsize
Le deuxième carré en partant du haut est commutatif car il est l'image par le foncteur
$\R u _{0*} ( \widetilde{\omega} _{X} \otimes ^\L _{\D _X}
\L u ^* _0( - ))$ de \ref{56prediag1}
(il n'y a pas d'ambiguïté pour calculer $u ^* _0$ car
$\D _{Y'} \otimes ^{\mathrm{dg}} _{\B _{Y'}} \D _{Y'}$ et
$\D _Y \otimes ^{\mathrm{dg}} _{\B _Y} \D _Y$
sont munis d'une unique structure {\it de $\D _{Y'} $-module à gauche}).
La commutativité des carrés du haut et du bas s'obtient par définition
(des morphismes du haut et du bas)
tandis que celle des deux autres s'établit par fonctorialité.
Pour obtenir celle du troisième carré du haut, il suffit de l'établir sans le foncteur $\R u _{0*}$. 
En résolvant platement $\widetilde{\omega} _X$, on se ramène à le vérifier sans le symbole $\L$, ce qui 
découle d'un calcul local (en se ramenant au cas relevable: voir \ref{ca}). 
Le diagramme \ref{56} est donc commutatif.
\end{vide}

\begin{vide}
  \label{4-5u*d=g56vide}
Vérifions à présent que
le morphisme composé de droite de \ref{4-5u*} est égal à celui de gauche de \ref{56}.
On bénéficie du diagramme commutatif ci-après :
\begin{equation}
  \label{4-5u*d=g56}
\xymatrix @R=0,3cm {
{\D _{X \rightarrow Y} \otimes _{\B _X} \D _{X\rightarrow Y}}
\ar[d] _-\sim
&
{u ^* _0( \D _Y \otimes ^{\mathrm{gg}} _{\B _Y} \D _Y)}
\ar[l] _-\sim
\ar[d] _-\sim
\\
{u ^* _0 F ^* F ^\flat \D _{Y'} \otimes  _{\B _{X}} u ^* _0 F ^* F ^\flat \D _{Y'} }
\ar[d] _-\sim
&
{u ^* _0 (F ^* F ^\flat \D _{Y'} \otimes ^{\mathrm{gg}} _{\B _{Y}} F ^* F ^\flat \D _{Y'} )}
\ar[l] _-\sim
\ar[d] _-\sim
\\
{F ^\flat _{\mathrm{d}} F ^\flat _{\mathrm{g}}(
u ^* _0 F ^*  \D _{Y'} \otimes  _{\B _{X}} u ^* _0 F ^* \D _{Y'} )}
\ar[d] _-\sim
&
{F ^\flat _{\mathrm{d}} F ^\flat _{\mathrm{g}}
u ^* _0 (F ^* \D _{Y'} \otimes ^{\mathrm{gg}} _{\B _{Y}} F ^* \D _{Y'} )}
\ar[l] _-\sim
\ar[d] _-\sim
\\
{F ^\flat _{\mathrm{d}} F ^\flat _{\mathrm{g}}(
F ^* \D _{X'\rightarrow Y'} \otimes  _{\B _{X}} F ^* \D _{X'\rightarrow Y'} )}
\ar[d] _-\sim
&
{F ^\flat _{\mathrm{d}} F ^\flat _{\mathrm{g}}
u ^* _0 F ^* (\D _{Y'} \otimes ^{\mathrm{gg}} _{\B _{Y'}}  \D _{Y'} )}
\ar[d] _-\sim
\\
{F ^\flat _{\mathrm{d}} F ^\flat _{\mathrm{g}}
F ^* (  \D _{X'\rightarrow Y'} \otimes  _{\B _{X'}} \D _{X'\rightarrow Y'})}
&
{F ^\flat _{\mathrm{d}} F ^\flat _{\mathrm{g}}
F ^* u ^{\prime *} _0  (\D _{Y'} \otimes ^{\mathrm{gg}} _{\B _{Y'}}  \D _{Y'} ).}
\ar[l] _-\sim
}
\end{equation}
En effet, la commutativité du rectangle (en bas) et du carré du milieu s'établit par un calcul local (en se ramenant au cas relevable: voir \ref{ca})
tandis que celle du carré du haut est fonctorielle.
En appliquant le foncteur $\R u _{0*} ( \widetilde{\omega} _{X} \otimes ^\L _{\D _X} -)$ au diagramme commutatif \ref{4-5u*d=g56}, on vérifie par fonctorialité
l'égalité entre le composé de droite de \ref{4-5u*} et celui de gauche de \ref{56}.
\end{vide}

\begin{vide}
Soient $\M \in D ^- (u ^{-1} _0 \D _Y \overset{^\mathrm{d}}{})$ et $\FF\in D ^- _{\mathrm{qc}} (\overset{^\mathrm{g}}{} \D _Y)$.
On calcule que le diagramme canonique
\begin{equation}
  \label{projproj}
\xymatrix @R=0,3cm {
{\R u _{0*} (\M \otimes ^{\L} _{u ^{-1} _0 \B _Y} u ^{-1} _0 \FF )}
\ar[d] _-\sim
&
{\R u _{0*} (\M )\otimes ^{\L}  _{\B _Y} \FF }
\ar[l] _-\sim  ^-{\mathrm{proj}}
\ar[d] _-\sim
\\
{\R u _{0*} (\M \otimes ^{\L} _{u ^{-1} _0 \D _Y} u ^{-1} _0 (\D _Y \otimes _{\B _Y} \FF) )}
&
{\R u _{0*} (\M ) \otimes ^{\L} _{\D _Y} (\D _Y \otimes _{\B _Y} \FF) ,}
\ar[l] _-\sim  ^-{\mathrm{proj}}
}
\end{equation}
où les flèches horizontales sont les isomorphismes de projection,
est commutatif. De même en remplaçant $\D _Y $ par $\D _{Y'}$ et
$\B _Y$ par $\B _{Y'}$.
On remarque que la commutativité de \ref{projproj}
implique que le morphisme du haut est $\D _Y$-linéaire.

Considérons le diagramme :
  \begin{equation}
  \label{67-1}
\xymatrix @R=0,3cm{
{\R u _{0*} (\widetilde{\omega} _X \otimes ^\L _{\D _X} \D _{X \rightarrow Y}
\otimes _{u ^{-1} _0 \B _Y} u ^{-1} _0 \D _Y )} \ar[d] _-\sim
&
{\R u _{0*} (\widetilde{\omega} _X \otimes ^\L _{\D _X} \D _{X \rightarrow Y})
\otimes _{\B _Y} \D _Y } \ar[d] _-\sim  \ar[l] _-\sim  ^-{\mathrm{proj}} \\
{\R u _{0*} (\widetilde{\omega} _X \otimes ^\L _{\D _X} \D _{X \rightarrow Y}
\otimes _{u ^{-1} _0 \D _Y} u ^{-1} _0 \D ^2 _Y)}
\ar[d] _-\sim
&
{\R u _{0*} (\widetilde{\omega} _X \otimes ^\L _{\D _X} \D _{X \rightarrow Y})
\otimes _{\D _Y}  \D ^2 _Y}
\ar[d] _-\sim \ar[l] _-\sim  ^-{\mathrm{proj}} \\
{\R u _{0*} (\widetilde{\omega} _X
\underset{\D _X}{\overset{\L}{\otimes}}
F ^* F ^\flat \D _{X '\rightarrow Y'}
\underset{u ^{\text{-}1} _0 \D _Y}{\otimes}
u ^{\text{-}1} _0 F ^*  F ^\flat _{\mathrm{d}} F ^\flat _{\mathrm{g}} \D ^2 _{Y'} )}
&
{\R u _{0*} (\widetilde{\omega} _X
\underset{\D _X}{\overset{\L}{\otimes}}
F ^* F ^\flat \D _{X '\rightarrow Y'})
\underset{\D _Y}{\otimes}
F ^* F ^\flat _{\mathrm{d}} F ^\flat _{\mathrm{g}} \D ^2 _{Y'}}
\ar[l] _-\sim  ^-{\mathrm{proj}}
\ar[d] _-\sim
\\
{F ^\flat _{\mathrm{d}} F ^\flat _{\mathrm{g}}  ( \R u _{0*} (\widetilde{\omega} _X
\underset{\D _X}{\overset{\L}{\otimes}}
F ^* F ^\flat \D _{X '\rightarrow Y'}
\underset{u ^{\text{-}1} _0 \D _Y}{\otimes}
u ^{\text{-}1} _0 F ^*  \D ^2 _{Y'} ))}
\ar[u] _-\sim  ^-{\mathrm{proj}}
&
{F ^\flat _{\mathrm{d}} F ^\flat _{\mathrm{g}} (\R u _{0*} (\widetilde{\omega} _X
\underset{\D _X}{\overset{\L}{\otimes}}
F ^* F ^\flat \D _{X '\rightarrow Y'})
\underset{\D _Y}{\otimes} F ^*  \D ^2 _{Y'}),}
\ar[l] _-\sim  ^-{\mathrm{proj}}
}
\end{equation}
\normalsize
où
$\D ^2 _{Y'} :=\D _{Y'} \otimes ^{\mathrm{dg}} _{\B _{Y'}} \D _{Y'}$
et
$\D ^2 _Y :=\D _Y \otimes ^{\mathrm{dg}} _{\B _Y} \D _Y$.
Par \ref{projproj}, le carré du haut est commutatif, tandis que celui du milieu l'est par fonctorialité
et celui du bas par transitivité des morphismes de projection.

\'Etudions le diagramme ci-dessous :
\begin{equation}
\label{67-2}
\xymatrix @R=0,3cm{
{ \R u _{0*} (\widetilde{\omega} _X
\underset{\D _X}{\overset{\L}{\otimes}}
F ^* F ^\flat \D _{X '\rightarrow Y'}
\underset{u ^{\text{-}1} _0 \D _Y}{\otimes}
u ^{\text{-}1} _0 F ^*  \D ^2 _{Y'} )}
&
{\R u _{0*} (\widetilde{\omega} _X
\underset{\D _X}{\overset{\L}{\otimes}}
F ^* F ^\flat \D _{X '\rightarrow Y'})
\underset{\D _Y}{\otimes} F ^*  \D ^2 _{Y'}}
\ar[l] _-\sim  ^-{\mathrm{proj}} \\
{ \R u _{0*} (\widetilde{\omega} _X
\underset{\D _X}{\overset{\L}{\otimes}}
F ^* \D _{X '\rightarrow Y'}
\underset{u ^{\text{-}1} _0 \D _{Y'}}{\otimes}
u ^{\text{-}1} _0 F ^\flat  \D _{Y'}
\underset{u ^{\text{-}1} _0\D _Y}{\otimes}  u ^{\text{-}1} _0 F ^*  \D ^2 _{Y'} )}
\ar[d] _-\sim  \ar[u] _-\sim
&
{ \R u _{0*} (\widetilde{\omega} _X
\underset{\D _X}{\overset{\L}{\otimes}}
F ^* \D _{X '\rightarrow Y'}
\underset{u ^{\text{-}1} _0 \D _{Y'}}{\otimes}
u ^{\text{-}1} _0 F ^\flat  \D _{Y'} ) \underset{\D _Y}{\otimes}  F ^*  \D ^2 _{Y'} }
\ar[l] _-\sim  ^-{\mathrm{proj}}
\ar[u] _-\sim
\\
{ \R u _{0*} (\widetilde{\omega} _X
\underset{\D _X}{\overset{\L}{\otimes}}
F ^* \D _{X '\rightarrow Y'}
\underset{u ^{\text{-}1} _0 \D _{Y'}}{\otimes}
u ^{\text{-}1} _0 ( F ^\flat \D _{Y'} \underset{\D _Y}{\otimes} F ^*  \D ^2 _{Y'} ))}
\ar[d] _-\sim
&
{\R u _{0*} (\widetilde{\omega} _X
\underset{\D _X}{\overset{\L}{\otimes}}
F ^* \D _{X '\rightarrow Y'})
\underset{\D _{Y'}}{\otimes}
(F ^\flat  \D _{Y'} \underset{\D _Y}{\otimes} F ^*  \D ^2 _{Y'})}
\ar[u] _-\sim  ^-{\mathrm{proj}}\ar[l] _-\sim  ^-{\mathrm{proj}} \ar[d] _-\sim  \\
{ \R u _{0*} (\widetilde{\omega} _{X'} \underset{\D _{X'}}{\overset{\L}{\otimes}}
\D _{X '\rightarrow Y'}
\underset{u ^{\text{-}1} _0 \D _{Y'}}{\otimes}
u ^{\text{-}1} _0  \D ^2 _{Y'} )}
&
{\R u _{0*} (\widetilde{\omega} _{X'} \underset{\D _{X'}}{\overset{\L}{\otimes}}
\D _{X '\rightarrow Y'})
\underset{\D _{Y'}}{\otimes} \D ^2 _{Y'}} \ar[l] _-\sim  ^-{\mathrm{proj}} \\
{ \R u _{0*} (\widetilde{\omega} _{X'} \underset{\D _{X'}}{\overset{\L}{\otimes}}
\D _{X '\rightarrow Y'}
\underset{u ^{\text{-}1} _0 \B _{Y'}}{\otimes} u ^{\text{-}1} _0  \D _{Y'} )}
\ar[u] _-\sim
&
{\R u _{0*} (\widetilde{\omega} _{X'} \underset{\D _{X'}}{\overset{\L}{\otimes}}
\D _{X '\rightarrow Y'})
 \underset{\B _{Y'}}{\otimes} \D _{Y'}.}
\ar[l] _-\sim  ^-{\mathrm{proj}} \ar[u] _-\sim
}
\end{equation}
\normalsize
Le deuxième carré du haut de \ref{67-2} est commutatif par transitivité des morphismes de projection, tandis que
celui du bas l'est par \ref{projproj} et les deux autres le sont par fonctorialité.
Le diagramme \ref{67-2} est donc commutatif.
En lui appliquant $F ^\flat _{\mathrm{d}} F ^\flat _{\mathrm{g}}$
puis en le composant avec \ref{67-1}, on parvient au diagramme commutatif :
\begin{equation}
  \label{67}
\xymatrix @R=0,3cm {
{\R u _{0*} (\widetilde{\omega} _X \otimes ^\L _{\D _X} \D _{X \rightarrow Y}
\otimes _{u ^{-1} _0 \B _Y} u ^{-1} _0 \D _Y )}
\ar[d] _-\sim
&
{u _{0+} (\widetilde{\omega} _X ) \otimes _{\B _Y} \D _Y }
\ar[d] _-\sim  \ar[l] _-\sim  ^-{\mathrm{proj}}
\\
{F ^\flat _{\mathrm{d}} F ^\flat _{\mathrm{g}} (
\R u _{0*} (\widetilde{\omega} _{X'} \otimes ^\L _{\D _{X'}} \D _{X '\rightarrow Y'}
\otimes _{u ^{\text{-}1} _0 \B _Y} u ^{\text{-}1} _0  \D _{Y'}))}
&
{F ^\flat _{\mathrm{d}} F ^\flat _{\mathrm{g}} (u _{0+} '  (\widetilde{\omega} _{X'} ) \otimes _{\B _{Y'}} \D _{Y'} ).}
\ar[l] _-\sim  ^-{\mathrm{proj}}
}
\end{equation}
\end{vide}

\begin{vide}
\label{6-6prime1vide}
\'Etablissons maintenant que le morphisme composé de droite de
\ref{56} est égal au morphisme de gauche de \ref{67}.

Grâce à la commutativité du diagramme du haut de (\cite[1.4.15.1]{caro_comparaison}),
pour tout $\D _{Y'}$-module à gauche $\FF'$, le composé
$F ^* \FF '  \riso \D _Y \otimes _{\D _Y} F ^* \FF'
\riso F ^* F ^\flat \D _{Y'} \otimes _{\D _Y} F ^* \FF'
\riso F ^* \D _{Y'} \otimes _{\D _{Y'}} \FF'$ est égal à l'isomorphisme canonique
$F ^* \FF '  \riso F ^* \D _{Y'} \otimes _{\D _{Y'}} \FF'$.
Il s'en suit que l'isomorphisme composé
$u ^* _0 F ^* (\FF')\riso
(u ^* _0 F ^*  \D _{Y'} )\otimes _{u ^{\text{-}1} _0 \D _{Y'}} u ^{\text{-}1} _0 \FF'$
du contour du diagramme
\begin{equation}
  \label{6-6prime3}
\xymatrix @R=0,3cm @C=0,5cm {
{u ^* _0 F ^* (\FF')}
\ar[rr] _-\sim\ar[ddd] _-\sim
&&
{\D _{X \rightarrow Y}
\otimes _{u ^{\text{-}1} _0 \D _Y} u ^{\text{-}1} _0 F ^* \FF '}
\ar[d] _-\sim
\\
&
{F ^*  u ^{\prime *} _0 F ^\flat \D _{Y'}
\otimes _{u ^{\text{-}1} _0 \D _Y} u ^{\text{-}1} _0 F ^* \FF '}
\ar[d] _-\sim
&
{u ^* _0 F ^* F ^\flat \D _{Y'}
\otimes _{u ^{\text{-}1} _0
\D _Y} u ^{\text{-}1} _0 F ^* \FF '}
\ar[d] _-\sim \ar[l] _-\sim
\\
&
{ F ^* u ^{\prime *} _0 \D _{Y'}
\otimes _{u ^{\text{-}1} _0 \D _{Y'}}
u ^{\text{-}1} _0 F ^\flat  \D _{Y'}
\otimes _{u ^{\text{-}1} _0\D _Y}
u ^{\text{-}1} _0 F ^*  \FF '}
\ar[d] _-\sim
&
{ u ^* _0 F ^* \D _{Y'}
\otimes _{u ^{\text{-}1} _0 \D _{Y'}}
u ^{\text{-}1} _0 F ^\flat  \D _{Y'}
\otimes _{u ^{\text{-}1} _0\D _Y} u ^{\text{-}1} _0 F ^*  \FF '}
\ar[l] _-\sim \ar[d] _-\sim
\\
{F ^* u ^{\prime *} _0 (\FF')}
\ar[r] _-\sim
&
{F ^* u ^{\prime *} _0  \D _{Y'} \otimes _{u ^{\text{-}1} _0 \D _{Y'}} u ^{\text{-}1} _0 \FF'}
&
{u ^* _0 F ^*  \D _{Y'} \otimes _{u ^{\text{-}1} _0 \D _{Y'}} u ^{\text{-}1} _0 \FF'}
\ar[l] _-\sim
}
\end{equation}
\normalsize
passant par le haut puis la droite est l'isomorphisme canonique.
Via un calcul local, il en découle la commutativité du contour de \ref{6-6prime3}.
Comme celle des carrés se vérifie par fonctorialité,
\ref{6-6prime3} est commutatif.

Considérons le diagramme suivant :
\begin{equation}
\label{6-6prime2}
\xymatrix @C=-0,5cm @R=0,3cm {
{\R u _{0*} ( \widetilde{\omega} _{X} \otimes ^\L _{\D _X} u ^* _0 F ^*  \D _{Y'} ^2 )}
\ar[dd] _-\sim \ar[r] _-\sim
&
{\R u _{0*} ( \widetilde{\omega} _X
\otimes  ^\L _{\D _X}
\D _{X \rightarrow Y}
\otimes _{u ^{\text{-}1} _0 \D _Y}
u ^{\text{-}1} _0 F ^*   \D ^2 _{Y'} )}
\ar[d] _-\sim
\\
&
{ \R u _{0*} (\widetilde{\omega}  _X
\underset{\D _X}{\overset{\L}{\otimes}}
F ^* \D _{X '\rightarrow Y'}
\underset{u ^{\text{-}1} _0 \D _{Y'}}{\otimes}
u ^{\text{-}1} _0 F ^\flat  \D _{Y'}
\underset{u ^{\text{-}1} _0\D _Y}{\otimes}  u ^{\text{-}1} _0 F ^*  \D ^2 _{Y'} )}
\ar[d] _-\sim
\\
{\R u _{0*} ( \widetilde{\omega}  _{X} \otimes ^\L _{\D _X} F ^* u ^{\prime *} _0 \D _{Y'} ^2 )}
\ar[d] _-\sim \ar[r] _-\sim
&
{ \R u _{0*} (\widetilde{\omega} _X \otimes ^\L _{\D _X}
F ^* \D _{X '\rightarrow Y'} \otimes _{u ^{\text{-}1} _0 \D _{Y'}} u ^{\text{-}1} _0 \D ^2 _{Y'} )}
\ar[d] _-\sim
\\
{\R u _{0*} ( \widetilde{\omega} _{X'} \otimes ^\L _{\D _{X'}} u ^{\prime *} _0 \D _{Y'} ^2 )}
\ar[d] _-\sim \ar[r] _-\sim
&
{ \R u _{0*} (\widetilde{\omega} _{X'} \otimes ^\L _{\D _{X'}}
\D _{X '\rightarrow Y'} \otimes _{u ^{\text{-}1} _0 \D _{Y'}} u ^{\text{-}1} _0 \D ^2 _{Y'} )}
\ar[d] _-\sim
\\
{\R u _{0*} ( \widetilde{\omega} _{X'} \otimes ^\L _{\D _{X'}}
( \D _{X '\rightarrow Y'} \otimes  _{ u ^{-1} _0 \B _{Y'}} u ^{-1} _0  \D _{Y'} ))}
\ar@{=}[r]
&
{\R u _{0*} ( \widetilde{\omega} _{X'} \otimes ^\L _{\D _{X'}}
( \D _{X '\rightarrow Y'} \otimes  _{ u ^{-1} _0 \B _{Y'}} u ^{-1} _0  \D _{Y'} )).}
}
\end{equation}
\normalsize
Comme $\D _{X '\rightarrow Y'} = u ^{\prime *} _0 \D _{Y'}$,
en appliquant $\R u _{0*} ( \widetilde{\omega} _{X} \otimes ^\L _{\D _X} - )$
au trapèze de \ref{6-6prime3} pour $\FF '= \D ^2 _{Y'}$ (ce dernier est muni d'une unique structure canonique de
$\D  _{Y'}$-module à gauche), on obtient
le rectangle du haut de \ref{6-6prime2}.
La commutativité du carré du milieu
s'établit par fonctorialité.
Pour vérifier celle du carré du bas, 
on peut y omettre $\R u _{0*} ( \widetilde{\omega} _{X'} \otimes ^\L _{\D _{X'}} -)$, ce qui résulte alors d'un calcul local immédiat.

En appliquant $F ^\flat _{\mathrm{d}} F ^\flat _{\mathrm{g}} $ à \ref{6-6prime2}
et en le composant au diagramme commutatif ci-dessous
\begin{equation}
\notag
\xymatrix @R=0,3cm @C=0cm {
{\R u _{0*} ( \widetilde{\omega} _{X} \otimes ^\L _{\D _X}
( \D _{X \rightarrow Y} \otimes  _{u ^{-1} _0 \B _Y} u ^{-1} _0 \D _Y))}
\ar@{=}[r]
\ar[d] _-\sim
&
{\R u _{0*} ( \widetilde{\omega} _{X} \otimes ^\L _{\D _X}
( \D _{X \rightarrow Y} \otimes  _{u ^{-1} _0 \B _Y} u ^{-1} _0 \D _Y))}
\ar[d] _-\sim
\\
{\R u _{0*} ( \widetilde{\omega} _{X} \otimes ^\L _{\D _X}
u ^* _0 \D ^2 _Y )}
 \ar[d] _-\sim 
 \ar[r] _-\sim
&
{\R u _{0*} (\widetilde{\omega} _{X} \otimes ^\L _{\D _X} \D _{X \rightarrow Y}
\otimes _{u ^{-1} _0 \D _Y} u ^{-1} _0 \D ^2 _Y)}
\ar[d] _-\sim
\\
{\R u _{0*} ( \widetilde{\omega} _{X} \otimes ^\L _{\D _X}
u ^* _0 F ^* F ^\flat _{\mathrm{g}} F ^\flat _{\mathrm{d}}  \D _{Y'} ^2 )}
\ar[r] _-\sim
&
{\R u _{0*} (\widetilde{\omega} _{X} \otimes ^\L _{\D _X}
\D _{X \rightarrow Y}
\otimes _{u ^{\text{-}1} _0 \D _Y}
u ^{\text{-}1} _0 F ^*  F ^\flat _{\mathrm{d}} F ^\flat _{\mathrm{g}} \D ^2 _{Y'} )}
\\
{F ^\flat _{\mathrm{g}} F ^\flat _{\mathrm{d}}  (\R u _{0*} ( \widetilde{\omega} _{X} \otimes ^\L _{\D _X}
u ^* _0 F ^*  \D _{Y'} ^2 ))}
\ar[u] _-\sim  ^-{\mathrm{proj}}
\ar[r] _-\sim
&
{F ^\flat _{\mathrm{d}} F ^\flat _{\mathrm{g}} \R u _{0*} (\widetilde{\omega} _{X} \otimes ^\L _{\D _X}
\D _{X \rightarrow Y}
\otimes _{u ^{\text{-}1} _0 \D _Y}
u ^{\text{-}1} _0 F ^*   \D ^2 _{Y'} ),}
\ar[u] _-\sim  ^-{\mathrm{proj}}
}
\end{equation}
on obtient un diagramme commutatif dont
le morphisme composé de gauche est égal à celui de droite de
\ref{56} et dont celui de droite s'identifie au morphisme de gauche de \ref{67}.

\end{vide}

\begin{prop}
  \label{prop4-7}
Pour tout $ \M ' \in D ^- ( \D _{X'} \overset{^{\mathrm{d}}}{})$,
on dispose du diagramme commutatif :
\begin{equation}
  \label{4-7}
\xymatrix @R=0,3cm {
{\R \mathcal{H} om _{\D_Y}
(u _{0+}  (F ^\flat  \M ' ),
u _{0+} (\widetilde{\omega} _X \otimes _{\B _X} \D _{X \rightarrow Y} ) )}
\ar[r] _-\sim  &
{\R \mathcal{H} om _{\D_Y}
(u _{0+}  (F ^\flat  \M ' ), u _{0+} (\widetilde{\omega} _X ) \otimes _{\B _Y} \D _{Y})}
\\
{\R \mathcal{H} om _{\D_Y}
(F ^\flat  u _{0+} '   ( \M ' ),
u _{0+} (\widetilde{\omega} _X \otimes _{\B _X} \D _{X \rightarrow Y} ) )}
\ar[r] _-\sim \ar[d] _-\sim \ar[u] _-\sim &
{\R \mathcal{H} om _{\D_Y}
(F ^\flat  u _{0+} '   ( \M ' ), u _{0+} (\widetilde{\omega} _X ) \otimes _{\B _Y} \D _{Y})}
\ar[d] _-\sim \ar[u] _-\sim
\\
{\R \mathcal{H} om _{\D_Y}
(F ^\flat  u _{0+} '   ( \M ' ),
F ^\flat _{\mathrm{g}} F ^\flat _{\mathrm{d}} (
u _{0+} '  (\widetilde{\omega} _{X'} \otimes _{\B _{X'}} \D _{X '\rightarrow Y'} ) ))} \ar[r] _-\sim&
{\R \mathcal{H} om _{\D_Y}
(F ^\flat  u _{0+} '   ( \M ' ), F ^\flat _{\mathrm{g}} F ^\flat _{\mathrm{d}}(u _{0+} '  (\widetilde{\omega} _{X'} ) \otimes _{\B _{Y'}} \D _{Y'}))}
\\
{\R \mathcal{H} om _{\D _{Y'}}
(u _{0+} '   ( \M ' ), F ^\flat _{\mathrm{d}} u _{0+} '  (\widetilde{\omega} _{X'} \otimes _{\B _{X'}} \D _{X '\rightarrow Y'} ) )}
\ar[r] _-\sim \ar[d] _-\sim  \ar[u] _-\sim ^{F ^\flat}
&
{\R \mathcal{H} om _{\D _{Y'}}
(u _{0+} '   ( \M ' ), F ^\flat _{\mathrm{d}}(u _{0+} '  (\widetilde{\omega} _{X '}) \otimes _{\B _{Y'}} \D _{Y'}))}
 \ar[u] _-\sim ^{F ^\flat}
\ar[d] _-\sim
\\
{F ^\flat \R \mathcal{H} om _{\D _{Y'}}
(u _{0+} '   ( \M ' ),
 u _{0+} '  (\widetilde{\omega} _{X'} \otimes _{\B _{X'}} \D _{X '\rightarrow Y'} ) )}
\ar[r] _-\sim &
 {F ^\flat \R \mathcal{H} om _{\D _{Y'}} (u _{0+} '   ( \M ' ), u _{0+} '  (\widetilde{\omega} _{X'} ) \otimes _{\B _{Y'}} \D _{Y'}),}
}
\end{equation}
\normalsize
où les flèches verticales du haut sont induites par \ref{u+comFrob}.
\end{prop}

\begin{proof}
En composant \ref{4-5u*}, \ref{56} et \ref{67}
(grâce à \ref{4-5u*d=g56vide} et \ref{6-6prime1vide} ces compositions ont un sens)
on obtient
(modulo $F ^\flat _{\mathrm{d}} F ^\flat _{\mathrm{g}} \riso F ^\flat _{\mathrm{g}} F ^\flat _{\mathrm{d}}$)
le diagramme commutatif :
\begin{equation}
  \label{47u*}
\xymatrix @R=0,3cm {
{u _{0+} ( \widetilde{\omega} _X \otimes _{\B _X} \D _{X \rightarrow Y}) }
\ar[r] _-\sim \ar[d] _-\sim
&
{u _{0+} (\widetilde{\omega} _X ) \otimes _{\B _Y} \D _Y }
\ar[d]_-\sim
\\
{F ^\flat _{\mathrm{g}} F ^\flat  _{\mathrm{d}} u _{0+} '   (\widetilde{\omega} _{X'} \otimes _{\B _{X'}} \D _{X '\rightarrow Y'})}
\ar[r] _-\sim
&
{F ^\flat _{\mathrm{d}} F ^\flat _{\mathrm{g}} (u _{0+} '  (\widetilde{\omega} _{X'} ) \otimes _{\B _{Y'}} \D _{Y'} ).}
}
\end{equation}
En lui appliquant
$\R \mathcal{H} om _{\D _{Y'}} (F ^\flat u _{0+} '   ( \M ' ), -)$,
il vient
le deuxième carré du haut de \ref{4-7}. Enfin, les trois autres le sont par fonctorialité.
 \end{proof}

\subsection{\label{ccn}Résultat par composition des trois précédents morphismes construits}

\begin{theo}
\label{theo1-7}
Avec les notations du chapitre \ref{ca},
on dispose,
pour tout $ \M ' \in D ^- ( \D _{X'} \overset{^{\mathrm{d}}}{})$,
du diagramme commutatif suivant
\begin{equation}
  \label{1-7}
\xymatrix @R=0,3cm {
{u _{0+} \DD ( F ^\flat  \M ') }
\ar[r]  \ar[d] _-\sim
&
{\R \mathcal{H} om _{\D_Y}
(u _{0+}  (F ^\flat  \M ' ), u _{0+} (\widetilde{\omega} _X ) \otimes _{\B _Y} \D _{Y})[d _{X}]}
\ar[d] _-\sim
\\
{F ^\flat  u _{0+} ' \DD (  \M ')}
\ar[r]
&
{F ^\flat \R \mathcal{H} om _{\D _{Y'}}
(u _{0+} '( \M ' ), u _{0+} '(\widetilde{\omega} _{X'} ) \otimes _{\B _{Y'}} \D _{Y'})[d _{X}].}
}
\end{equation}
\end{theo}
\begin{proof}
  On vérifie que le morphisme composé de droite de \ref{1-2}
est égal, modulo
$\widetilde{\omega} _{X'} \otimes _{\B _{X'}} \D _{X '\rightarrow Y'}
\riso \widetilde{\omega} _{X'} \otimes _{\B _{X'}} \D _{X'} \otimes _{\D _{X'} } \D _{X '\rightarrow Y'}$,
 au morphisme composé de gauche de \ref{2-4}.
De plus, le composé de droite de \ref{2-4} est celui de gauche de \ref{4-7}.
  En composant \ref{1-2}, \ref{2-4} et \ref{4-7}, on aboutit alors au diagramme ci-dessous :
\begin{equation}
  \label{1-7diag1}
\xymatrix @R=0,3cm @C=0,3cm {
{u _{0+} ( \R \mathcal{H} om _{\D _X} ( F ^\flat  \M ', \widetilde{\omega} _X \otimes _{\B _X}
\! \! \! \D _X))}
\ar[r]  \ar[d] _-\sim
&
{\R \mathcal{H} om _{\D_Y}
(u _{0+}  (F ^\flat  \M ' ), u _{0+} (\widetilde{\omega} _X ) \otimes _{\B _Y} \! \! \!\D _{Y})}
\ar[d] _-\sim
\\
{F ^\flat  u _{0+} '  ( \R \mathcal{H} om _{\D _{X'}} (  \M ', \widetilde{\omega} _{X'} \otimes _{\B _{X'}}\! \! \!\D _{X'}))}
\ar[r]
&
{F ^\flat \R \mathcal{H} om _{\D _{Y'}} (u _{0+} '   ( \M ' ), u _{0+} '  (\widetilde{\omega} _{X'} ) \otimes _{\B _{Y'}}\! \! \! \D _{Y'}).}
}
\end{equation}
On conclut via \ref{1-2bis} et \ref{rema1-2}.
 \end{proof}

\begin{vide}
[{\bf Dernière étape de la construction de l'isomorphisme de dualité relative: le morphisme trace.}]
\label{1.5.2}
  D'après le théorème ci-dessus,
pour tout $ \M ' \in D ^- ( \D _{X'} \overset{^{\mathrm{d}}}{})$,
le morphisme canonique
\begin{equation}
  \label{1-7ligne}
u _{0+} ' \DD (  \M ')
\rightarrow
\R \mathcal{H} om _{\D _{Y'}} (u _{0+} '   ( \M ' ), u _{0+} '  (\widetilde{\omega} _{X'} ) \otimes _{\B _{Y'}} \D _{Y'}) [d _{X'}]
\end{equation}
est compatible à Frobenius.
Dans le cas relevable (voir \ref{ca}), 
Virrion a construit un morphisme trace (\cite{Vir04})
$\mathrm{Tr} _+$ :
$u _{+}  (\omega _{X'} ) [d _{X'}] \rightarrow \omega _{Y'} [d _{Y'}]$.
Lorsque $u _{0} ^{-1}\B _{Y}$ et $\O _{X}$ sont Tor-indépendants, il en résulte le morphisme trace 
$\mathrm{Tr} _+$ :
$u _{+}  (\widetilde{\omega} _{X'} ) [d _{X'}] \rightarrow \widetilde{\omega} _{Y'} [d _{Y'}]$, car d'après Virrion (voir \cite{Vir04}), avec cette hypothèse l'image directe commute à l'extension via $\O _{Y}\to \B _{Y}$.
On en déduit par fonctorialité le premier morphisme : 
\begin{gather}
  \R \mathcal{H} om _{\D _{Y'}} (u _{+}   ( \M ' ), u _{+}  (\widetilde{\omega} _{X'} ) \otimes _{\B _{Y'}} \D _{Y'}) [d _{X'}]
\label{1-7ligne2}
\rightarrow
\R \mathcal{H} om _{\D _{Y'}} (u _{+}   ( \M ' ), \widetilde{\omega} _{Y'} \otimes _{\B _{Y'}} \D _{Y'}) [d _{Y'}]
\riso
\DD u _{+} (  \M ').
\end{gather}
En composant \ref{1-7ligne} avec \ref{1-7ligne2},
on obtient le morphisme
$u _{+} \DD (  \M ') \rightarrow \DD u _{+} (  \M ')$ construit par Virrion dans (\cite{Vir04}).
Lorsque
$ \M ' \in D ^\mathrm{b} _\mathrm{parf}
( \D _{X'} \overset{^{\mathrm{d}}}{})$ et $S$ est régulier,
Virrion a prouvé que celui-ci est un isomorphisme.
On l'appelle {\it isomorphisme de dualité relative}.

Pour établir que l'isomorphisme de dualité relative est compatible à Frobenius, il suffit
de démontrer que le morphisme trace
$\mathrm{Tr} _+$ :
$u _{+}  (\widetilde{\omega} _{X'} ) [d _{X'}] \rightarrow \widetilde{\omega} _{Y'} [d _{Y'}]$
l'est. C'est ce que nous ferons au niveau des schémas formels
dans le deuxième chapitre dans le cas d'une immersion fermée.
\end{vide}

\subsection{Compatibilité au changement de niveaux}

Afin de vérifier que l'isomorphisme de dualité relative d'une immersion fermée est compatible à Frobenius, nous aurons besoin de vérifier sa compatibilité au changement de niveau sur les schémas (plus précisément, afin d'obtenir \ref{isodualrelind}). En effet, comme nous travaillerons avec des complexes quasi-cohérents (au sens de Berthelot) sur les $\V$-schémas formels lisses, cette compatibilité sur les schémas sera requise. Notons que
comme,  pour tout $\V$-schéma formel  $\X$ lisse, les extensions $\widehat{\D} ^{(m)} _{\X,\Q}\to \widehat{\D} ^{(m+1)} _{\X,\Q}$ sont plates, la version "formelle tensorisé par $\Q$" de cette section est facile à vérifier (e.g., l'analogue formel tensorisé par $\Q$ de \ref{24niveaulembis} et donc de \ref{theo17niveau} sont aisément validé sans supposer que $u _{0}$ est une immersion fermée). 

Fixons deux entiers naturels $m' \geq m$.
Nous noterons alors dans toutes ces sections, $\D _{?}:= 
\B _{X} ^{(m)} \otimes _{\O _{X}} \D ^{(m)} _{?}$, 
$\widetilde{\D} _{?}:= 
\B _{X} ^{(m')} \otimes _{\O _{X}} \D ^{(m')} _{?}$. 
De plus, on pose 
$\B _{X}:= \B _{X} ^{(m)}$, $u _{0+} :=u _{0+} ^{(m)}$, $\DD:=\DD ^{(m)}$ ;
$\widetilde{\B} _{X}:= \B _{X} ^{(m')}$,
$\widetilde{u} _{0+} :=u _{0+} ^{(m')} $, 
$\widetilde{\DD}:=\DD ^{(m')}$.
De même, en remplaçant $X$ par $Y$.

Nous supposons dans cette section que 
le faisceau $\widetilde{\D} _{Y}$ est Tor-dimension finie à droite et à gauche sur $\D _{Y}$ (e.g., lorsque $\B _{X} ^{(m)} = \O _{X}$ d'après \cite[1.2.6]{Beintro2}).

\begin{vide}
  Soient $ \M \in D ^- (\D _{X} \overset{^{\mathrm{d}}}{})$,
$\NN \in D ^+ (\D _{X} \overset{^{\mathrm{d}}}{}, \D _{X} \overset{^{\mathrm{d}}}{})$,
$\widetilde{\NN} \in D ^+ (\widetilde{\D} _X \overset{^{\mathrm{d}}}{}, \widetilde{\D} _X \overset{^{\mathrm{d}}}{})$
et $\NN \rightarrow \widetilde{\NN}$ un morphisme $\D _{X}$-bilinéaire. On pose
$\widetilde{\M} :=  \M \otimes ^{\L} _{\D _{X}} \widetilde{\D} _X$.
On dispose du morphisme de changement de niveau
$$\R \mathcal{H} om _{ \D _{X}} ( \M, \NN ) \rightarrow
\R \mathcal{H} om _{ \D _{X}} ( \M, \widetilde{\NN} )
 \overset{\sim}{\leftarrow}
\R \mathcal{H} om _{\widetilde{\D} _X } (\widetilde{\M}, \widetilde{\NN} ),$$
où l'isomorphisme de droite se calcule en résolvant $\M$ par des $\D _{X}$-modules à droites plats et
$\widetilde{\NN}$ par des $\widetilde{\D} _X $-bimodules à droites injectifs (on utilise le lemme de \cite[I.1.4]{virrion}).

Il en découle en particulier {\it le morphisme de changement de niveau du foncteur dual}
$\DD (  \M) \rightarrow \widetilde{\DD} (\widetilde{\M})$.
\end{vide}

\begin{vide}
 Soient $\M \in D ^{-} (\D _{X} \overset{^{\mathrm{d}}}{})$ et
$\widetilde{\M} :=  \M \otimes ^{\L} _{\D _{X}} \widetilde{\D} _X$.
On dispose du 
{\it morphisme canonique de changement de niveaux de l'image directe}
$u _{0+} ( \M) \to \widetilde{u} _{0+} (\widetilde{\M})$. Il se construit via l'un des deux composé du diagramme commutatif suivant :
\begin{equation}
\label{chgtniveu+1}
  \xymatrix @R=0,3cm @C=0,5cm {
{u _{0+} ( \M)} \ar@{=}[r] 
&
{ \R u _{0*} ( \M \otimes ^\L _{\D _{X}} \D _{X \rightarrow Y})}
\ar[r] \ar[d]
&
{\R u _{0*} ( \M \otimes _{\D _{X}} ^\L \widetilde{\D} _{X \rightarrow Y}) }
\ar[d] _-\sim
\\
{}
&
{ \R u _{0*} (\widetilde{\M} \otimes ^\L _{\D _{X}} \D _{X \rightarrow Y})} \ar[r] &
{\R u _{0*} (\widetilde{\M} \otimes ^\L _{\widetilde{\D} _X} \widetilde{\D} _{X \rightarrow Y})}
\ar@{=}[r] &
{\widetilde{u} _{0+} (\widetilde{\M})}.}
\end{equation}

\end{vide}

\begin{vide}
  \label{viderema12niveau}
Soient $ \M \in D ^- (\D _{X} \overset{^{\mathrm{d}}}{})$
et $\widetilde{\M} :=  \M \otimes ^{\L} _{\D _{X}} \widetilde{\D} _X$.
L'isomorphisme de transposition $\delta$ (\cite[1.3.1]{Be2}) de
$\widetilde{\omega} ^{(m)} _{X} \otimes _{\B _{X}}\D _{X}$ est compatible au changement de niveaux, i.e.,
le diagramme
\begin{equation}
  \label{transpchgtniv}
\xymatrix @R=0,3cm {
{\widetilde{\omega} ^{(m)} _{X} \otimes _{\B _{X}}\D _{X}} \ar[r] \ar[d] ^-{\delta} _\sim
&
{\widetilde{\omega} ^{(m')} _X \otimes _{\widetilde{\B} _X} \widetilde{\D} _X} \ar[d] ^-{\delta} _\sim
\\
{\widetilde{\omega} ^{(m)} _{X} \otimes _{\B _{X}}\D _{X}} \ar[r]
&
{\widetilde{\omega} ^{(m')} _X \otimes _{\widetilde{\B} _{X}} \widetilde{\D} _X}
}
\end{equation}
est commutatif. En effet, avec les formules \cite[2.2.3.1]{Be1} et \cite[1.3.3.1]{Be2}, cela dérive d'un calcul local.
On en déduit la commutativité de la partie gauche du diagramme suivant :
\small
\begin{equation}
  \label{rema12niveau}
\xymatrix @R=0,3cm {
{u _{0+} \DD ( \M )}
\ar[r] ^-{\sim} _-{\delta} \ar[dd]
&
{\R u _{0*} (\R \mathcal{H} om _{\D _{X}}
(  \M, \widetilde{\omega} ^{(m)} _{X} \otimes _{\B _{X}}\D _{X}) \otimes _{\D _{X} } ^\L \D _{X \rightarrow Y})}
\ar@{=}[r]\ar[d]
&
{\R u _{0*} (\R \mathcal{H} om _{\D _{X}}
(  \M, \widetilde{\omega} ^{(m)} _{X} \otimes _{\B _{X}}\D _{X}) \otimes _{\D _{X} } ^\L \D _{X \rightarrow Y})}
\ar[d]
\\
&
{\R u _{0*} (\R \mathcal{H} om _{\D _{X}}
(  \M, \widetilde{\omega} ^{(m')} _{X} \otimes _{\widetilde{\B} _{X}} \widetilde{\D} _{X}) \otimes _{\D _{X} } ^\L \D _{X \rightarrow Y})}
\ar[r]
&
{\R u _{0*} (\R \mathcal{H} om _{\D _{X}}
(  \M, \widetilde{\omega} ^{(m')} _{X} \otimes _{\widetilde{\B} _{X}} \widetilde{\D} _X) \otimes _{\widetilde{\D} _X } ^\L \widetilde{\D} _{X\rightarrow Y})}
\\
{u _{0+}\widetilde{\DD} ( \widetilde{\M})}
\ar[r] ^-{\sim} _-{\delta}
\ar[d]
&
{\R u _{0*} (\R \mathcal{H} om _{\widetilde{\D} _{X}}
(  \widetilde{\M}, \widetilde{\omega} ^{(m')} _{X} \otimes _{\widetilde{\B} _{X}} \widetilde{\D} _{X}) \otimes _{\D _{X} } ^\L \D _{X \rightarrow Y})}
\ar[u] ^-\sim
\ar[r] \ar[d]
&
{\R u _{0*} (\R \mathcal{H} om _{\widetilde{\D} _X}
( \widetilde{\M}, \widetilde{\omega} ^{(m')} _{X} \otimes _{\widetilde{\B} _{X}} \widetilde{\D} _X)
\otimes _{\widetilde{\D} _X } ^\L \widetilde{\D} _{X\rightarrow Y})}
\ar[u] ^-\sim \ar@{=}[d]
\\
{\widetilde{u} _{0+}\widetilde{\DD} ( \widetilde{\M})}
\ar[r] ^-{\sim} _-{\delta}
&
{\R u _{0*} (\R \mathcal{H} om _{\widetilde{\D} _X}
( \widetilde{\M}, \widetilde{\omega} ^{(m')} _{X} \otimes _{\widetilde{\B} _{X}} \widetilde{\D} _X)
\otimes _{\widetilde{\D} _X }^\L  \widetilde{\D} _{X\rightarrow Y})}
\ar@{=}[r]
&
{\R u _{0*} (\R \mathcal{H} om _{\widetilde{\D} _X}
( \widetilde{\M}, \widetilde{\omega} ^{(m')} _{X} \otimes _{\widetilde{\B} _{X}} \widetilde{\D} _X)
\otimes _{\widetilde{\D} _X }^\L  \widetilde{\D} _{X\rightarrow Y}).}
}
\end{equation}
\normalsize
Celle de la partie droite se vérifie par fonctorialité ou définition.
\end{vide}

\begin{prop}
\label{nota12niveau}
Soient $ \M \in D ^- (\D _{X} \overset{^{\mathrm{d}}}{})$
et $\widetilde{\M} :=  \M \otimes ^{\L} _{\D _{X}} \widetilde{\D} _X$.
On bénéficie du diagramme commutatif suivant :
\begin{equation}
  \label{12niveau}
\xymatrix @R=0,3cm @C=0,5cm {
{\R u _{0*} (\R \mathcal{H} om _{\D _{X}}
(  \M, \widetilde{\omega} ^{(m)} _{X} \otimes _{\B _{X}}\D _{X}) \otimes _{\D _{X} } ^\L \D _{X \rightarrow Y})}
\ar[r] \ar[d]
&
{\R u _{0*} \R \mathcal{H} om _{\D _{X}}
(  \M, \widetilde{\omega} ^{(m)} _{X} \otimes _{\B _{X}} \D _{X \rightarrow Y})}
\ar[d]
\\
{\R u _{0*} (\R \mathcal{H} om _{\D _{X}}
(  \M, \widetilde{\omega} ^{(m')} _{X} \otimes _{\widetilde{\B} _{X}} \widetilde{\D} _X) \otimes _{\widetilde{\D} _X } ^\L  \widetilde{\D} _{X\rightarrow Y})}
\ar[r]
&
{\R u _{0*} \R \mathcal{H} om _{\D _{X}}
(  \M, \widetilde{\omega} ^{(m')} _{X} \otimes _{\widetilde{\B} _{X}} \widetilde{\D} _{X\rightarrow Y})}
\\
{\R u _{0*} (\R \mathcal{H} om _{\widetilde{\D} _X}
( \widetilde{\M}, \widetilde{\omega} ^{(m')} _{X} \otimes _{\widetilde{\B} _{X}} \widetilde{\D} _X)
\otimes _{\widetilde{\D} _X } ^\L \widetilde{\D} _{X\rightarrow Y})}
\ar[r] \ar[u] ^-\sim
&
{\R u _{0*} \R \mathcal{H} om _{\widetilde{\D} _X}
( \widetilde{\M}, \widetilde{\omega} ^{(m')} _{X} \otimes _{\widetilde{\B} _{X}} \widetilde{\D} _{X\rightarrow Y}).}
\ar[u] ^-\sim
}
\end{equation}
\normalsize
\end{prop}
\begin{proof}
Par fonctorialité, il suffit de l'établir sans les foncteurs $\R u _{0*}$.
Pour vérifier la commutativité du carré du haut, on se ramène par fonctorialité à l'établir pour
$\widetilde{\omega} ^{(m')} _{X} \otimes _{\widetilde{\B} _{X}} \widetilde{\D} _X$ à la place de $\widetilde{\omega} ^{(m)} _{X} \otimes _{\B _{X}}\D _{X}$
et $\widetilde{\omega} ^{(m')} _{X} \otimes _{\widetilde{\B} _{X}} \widetilde{\D} _X \otimes _{\D _{X} } ^\L \D _{X \rightarrow Y}$
 à la place de $\widetilde{\omega} ^{(m)} _{X} \otimes _{\B _{X}} \D _{X \rightarrow Y}$.
Pour cela, il s'agit de résoudre 
$\M$ par des $\widetilde{\D} _{X}$-modules à droite plats, 
$\widetilde{\omega} ^{(m')} _{X} \otimes _{\widetilde{\B} _{X}} \widetilde{\D} _X$
par des $\widetilde{\D} _{X}$-bimodules à droite injectifs,
de prendre des résolutions finies de $\D _{X \rightarrow Y}$
par des $(\D _{X}, u _0 ^{-1} \D _{Y} )$-bimodules plats,
de $\widetilde{\D} _{X \rightarrow Y}$ par
des $(\widetilde{\D} _{X}, u _0 ^{-1} \widetilde{\D} _{Y} )$-bimodules plats
(compatibles avec le morphisme
$\D _{X \rightarrow Y}\rightarrow \widetilde{\D} _{X \rightarrow Y}$).

Enfin, pour obtenir celle du carré du bas,
on résout $\M$ par des $\D _X$-modules à droite plats,
$\widetilde{\omega} ^{(m')} _{X} \otimes _{\widetilde{\B} _{X}} \widetilde{\D} _{X}$ par des
$\widetilde{\D} _{X}$-bimodules à droite injectifs,
$\widetilde{\D} _{X \rightarrow Y}$
par des $(\widetilde{\D} _{X}, u _0 ^{-1} \widetilde{\D} _{Y} )$-bimodules plats.
 \end{proof}

Afin d'établir \ref{24niveauprop}, nous aurons besoin du lemme ci-dessous. 

\begin{lemm}
\label{2=2}
Soient $ \M  \in D  ^-   (\D _{X} \overset{^{\mathrm{d}}}{})$,
$\widetilde{\NN} \in D ^+ (\widetilde{\D} _X \overset{^{\mathrm{d}}}{})$
et $\widetilde{\M} :=  \M \otimes ^{\L} _{\D _{X}} \widetilde{\D} _X$.
On dispose du diagramme commutatif
\begin{equation}
\label{2=2diag}
\xymatrix @R=0,3cm {
{\R \mathcal{H} om _{\D _{X}} (  \M, \widetilde{\NN})  }
\ar[r]
&
{\R \mathcal{H} om _{u _0 ^{\text{-}1} \D _{Y}} (  \M \otimes _{\D _{X}} ^\L \D _{X \rightarrow Y},
\widetilde{\NN}\otimes _{\D _{X}}  ^\L \D _{X \rightarrow Y})  }
\ar[d]
\\
&
{\R \mathcal{H} om _{u _0 ^{\text{-}1} \D _{Y}} (  \M \otimes _{\D _{X}}  ^\L \D _{X \rightarrow Y},
\widetilde{\NN}\otimes _{\widetilde{\D} _{X}}  ^\L  \widetilde{\D} _{X \rightarrow Y})  }
\\
{\R \mathcal{H} om _{\widetilde{\D} _{X}} (  \widetilde{\M}, \widetilde{\NN})  }
\ar[r]
\ar[uu] _-\sim
&
{\R \mathcal{H} om _{u _0 ^{\text{-}1} \widetilde{\D} _{Y}}
(\widetilde{\M} \otimes _{\widetilde{\D} _{X}}  ^\L \widetilde{\D} _{X \rightarrow Y},
\widetilde{\NN}\otimes _{\widetilde{\D} _{X}}  ^\L \widetilde{\D} _{X \rightarrow Y} ) .}
\ar[u] _-\sim
}
\end{equation}

\end{lemm}
\begin{proof}
Donnons d'abord la construction de l'isomorphisme de droite du bas de \ref{2=2diag}.
Comme $u _0 ^{-1} \D _Y$ et $u _0 ^{-1} \widetilde{\D} _Y$
sont $u _0 ^{-1} \B  _Y$-plats à gauche (et à droite), il s'en suit :
$\D _{X \rightarrow Y} \otimes _{u _0 ^{-1}\D _Y } ^\L  u _0 ^{-1} \widetilde{\D} _Y
\riso \widetilde{\D} _{X \rightarrow Y}$.
D'où :
\begin{equation}
\label{chgtniveu+1bis}
\M \otimes _{\D _{X}} ^\L
\D _{X \rightarrow Y} \otimes _{u _0 ^{-1}\D _Y } ^\L  u _0 ^{-1} \widetilde{\D} _Y
\riso
\M \otimes _{\D _{X}} ^\L \widetilde{\D} _{X \rightarrow Y}
\riso
\widetilde{\M} \otimes _{\widetilde{\D} _{X}} ^\L \widetilde{\D} _{X \rightarrow Y}.
\end{equation}
On en déduit l'isomorphisme de droite du bas de \ref{2=2diag}. 

En résolvant 
$\M$ par des $\D _X$-modules à droite plats,
$\widetilde{\NN}$ par des
$\widetilde{\D} _{X}$-modules à droite injectifs,
$\D _{X \rightarrow Y}$
par des $(\D _{X}, u _0 ^{-1} \D _{Y} )$-bimodules plats,
$\widetilde{\D} _{X \rightarrow Y}$
par des $(\widetilde{\D} _{X}, u _0 ^{-1} \widetilde{\D} _{Y} )$-bimodules plats (et ce de façon bornée et compatible), on se ramène à établir la commutativité de \ref{2=2diag} sans les symboles $\R$ et $\L$, ce qui se vérifie à la main.
 \end{proof}

\begin{prop}  \label{24niveauprop}
Soient $ \M \in D ^- (\D _{X} \overset{^{\mathrm{d}}}{})$
et $\widetilde{\M} :=  \M \otimes ^{\L} _{\D _{X}} \widetilde{\D} _X$.
On bénéficie du diagramme commutatif ci-après :
  \begin{equation}
  \label{24niveau}
\xymatrix @R=0,3cm {
{\R u _{0*} \R \mathcal{H} om _{\D _{X}} (  \M, \widetilde{\omega} ^{(m)} _{X} \otimes _{\B _{X}} \D _{X \rightarrow Y})  }
\ar[r] \ar[d]
&
{\R u _{0*} \R \mathcal{H} om _{u _0 ^{\text{-}1} \D _{Y}}
(  \M \otimes _{\D _{X}} ^\L \D _{X \rightarrow Y} ,
(\widetilde{\omega} ^{(m)} _{X} \otimes _{\B _{X}} \D _{X \rightarrow Y})
\otimes _{\D _{X}}  ^\L  \D _{X \rightarrow Y})  }
\ar[d]
\\
{\R u _{0*} \R \mathcal{H} om _{\D _{X}} (  \M, \widetilde{\omega} ^{(m')} _{X} \otimes _{\widetilde{\B} _{X}} \widetilde{\D} _{X \rightarrow Y})  }
\ar[r]
&
{\R u _{0*} \R \mathcal{H} om _{u _0 ^{\text{-}1} \D _{Y}} (  \M \otimes _{\D _{X}}  ^\L  \D _{X \rightarrow Y},
(\widetilde{\omega} ^{(m')} _{X} \otimes _{\widetilde{\B} _{X}} \widetilde{\D} _{X \rightarrow Y})
\otimes _{\D _{X}}  ^\L  \D _{X \rightarrow Y}      )   }
\ar[d]
\\
&
{\R u _{0*} \R \mathcal{H} om _{u _0 ^{\text{-}1} \D _{Y}} (  \M \otimes _{\D _{X}}  ^\L  \D _{X \rightarrow Y},
(\widetilde{\omega} ^{(m')} _{X} \otimes _{\widetilde{\B} _{X}} \widetilde{\D} _{X \rightarrow Y})
\otimes _{\widetilde{\D} _{X}} ^\L   \widetilde{\D} _{X \rightarrow Y}      )   }
\\
{\R u _{0*} \R \mathcal{H} om _{\widetilde{\D} _X}
( \widetilde{\M}, \widetilde{\omega} ^{(m')} _{X} \otimes _{\widetilde{\B} _{X}} \widetilde{\D} _{X \rightarrow Y})  }
\ar[r] \ar[uu] ^-\sim
&
{\R u _{0*} \R \mathcal{H} om _{u _0 ^{\text{-}1} \widetilde{\D} _{Y}}
(  \widetilde{\M} \otimes _{\widetilde{\D} _{X}}  ^\L  \widetilde{\D} _{X \rightarrow Y} ,
(\widetilde{\omega} ^{(m')} _{X} \otimes _{\widetilde{\B} _{X}} \widetilde{\D} _{X \rightarrow Y})
\otimes _{\widetilde{\D} _{X}}  ^\L  \widetilde{\D} _{X \rightarrow Y}      )   ,}
\ar[u] ^-\sim
}
\end{equation}
\normalsize
dont l'isomorphisme de droite du bas est induit par \ref{chgtniveu+1bis}.

\end{prop}
\begin{proof}
La commutativité du rectangle du bas dérive de \ref{2=2} tandis que celle du
carré supérieur se vérifie par fonctorialité.
 \end{proof}

\begin{lemm}
\label{24niveaulembis}
On suppose que $u _{0}$ est une immersion fermée. 
Soient $\NN \in D ^-  (u _0 ^{\text{-}1} \D _{Y})$, 
$\widetilde{\NN} \in  D ^+  (u _0 ^{\text{-}1} \widetilde{\D} _{Y},u _0 ^{\text{-}1} \widetilde{\D} _{Y})$. 
Le diagramme suivant est commutatif : 
  \begin{equation}
  \label{24niveaubisdiaglemm}
\xymatrix @R=0,3cm @C=0,5cm {
{ u _{0*} \R \mathcal{H} om _{u _0 ^{\text{-}1} \D _{Y}} (  \NN , \widetilde{\NN}   )   }
\ar[r] ^-{ u _{0*}}
&
{\R \mathcal{H} om _{\D _{Y}} (    u _{0*}( \NN ),  u _{0*}(\widetilde{\NN}  )   )} 
\\
&
{\R \mathcal{H} om _{\widetilde{\D} _{Y}} 
(  u _{0*} (\NN   )             \otimes _{ \D _{Y}} ^{\L}\widetilde{\D} _{Y} ), 
 u _{0*} ( \widetilde{\NN} ) )}
\ar[u] ^-\sim 
\\
{u _{0*} \R \mathcal{H} om _{u _0 ^{\text{-}1} \widetilde{\D} _{Y}}
 (  \NN                \otimes _{u _0 ^{\text{-}1} \D _{Y}} ^{\L} u _0 ^{\text{-}1} \widetilde{\D} _{Y} , 
 \widetilde{\NN}  )     }
\ar[r] ^-{u _{0*}}
\ar[uu] ^-\sim
&
{\R \mathcal{H} om _{\widetilde{\D} _{Y}} (  u _{0*} (\NN                \otimes _{u _0 ^{\text{-}1} \D _{Y}} ^{\L} u _0 ^{\text{-}1} \widetilde{\D} _{Y} ), 
u _{0*} ( \widetilde{\NN} )  ). }
\ar[u] ^-\sim _{\mathrm{proj}}
}
\end{equation}

\end{lemm}

\begin{proof}
Il suffit de résoudre $\NN$ par des $u _0 ^{\text{-}1} \D _{Y} $-modules plats et $\widetilde{\NN}$ par des $u _0 ^{\text{-}1} \widetilde{\D} _{Y} $-bimodules injectifs.
\end{proof}

\begin{rema}
Il est possible que l'hypothèse de \ref{24niveaulembis} que $u _{0}$ est une immersion fermée soit inutile. 
\end{rema}

\begin{prop}
\label{24niveaupropbis}
On suppose que $ u _{0}$ est une immersion fermée. 
Soient $ \M \in D ^- (\D _{X} \overset{^{\mathrm{d}}}{})$
et $\widetilde{\M} :=  \M \otimes ^{\L} _{\D _{X}} \widetilde{\D} _X$.
Le diagramme suivant
  \begin{equation}
  \label{24niveaubis}
\xymatrix @R=0,3cm @C=0,5cm {
{u _{0*} \R \mathcal{H} om _{u _0 ^{\text{-}1} \D _{Y}}
(  \M \otimes _{\D _{X}}  ^\L  \D _{X \rightarrow Y} ,
(\widetilde{\omega} ^{(m)} _{X} \otimes _{\B _{X}} \D _{X \rightarrow Y})
\otimes _{\D _{X}} ^\L  \D _{X \rightarrow Y})  }
\ar[d]
\ar[r] ^-{u _{0*}}
&
{\R \mathcal{H} om _{\D _{Y}}
( u _{0+}( \M ),
u _{0+} (\widetilde{\omega} ^{(m)} _{X} \otimes _{\B _{X}} \D _{X \rightarrow Y}))}
\ar[d]
\\
{u _{0*} \R \mathcal{H} om _{u _0 ^{\text{-}1} \D _{Y}} (  \M \otimes _{\D _{X}} ^\L  \D _{X \rightarrow Y},
(\widetilde{\omega} ^{(m')} _{X} \otimes _{\widetilde{\B} _{X}} \widetilde{\D} _{X \rightarrow Y})
\otimes _{\widetilde{\D} _{X}}  ^\L  \widetilde{\D} _{X \rightarrow Y}      )   }
\ar[r] ^-{u _{0*}}
&
{\R \mathcal{H} om _{\D _{Y}} (   u _{0+}( \M ),
 \tilde{u} _{0+}(\widetilde{\omega} ^{(m')} _{X} \otimes _{\widetilde{\B} _{X}} \widetilde{\D} _{X \rightarrow Y})  )   }
\\
{u _{0*} \R \mathcal{H} om _{u _0 ^{\text{-}1} \widetilde{\D} _{Y}}
(  \widetilde{\M} \otimes _{\widetilde{\D} _{X}}  ^\L  \widetilde{\D} _{X \rightarrow Y} ,
(\widetilde{\omega} ^{(m')} _{X} \otimes _{\widetilde{\B} _{X}} \widetilde{\D} _{X \rightarrow Y})
\otimes _{\widetilde{\D} _{X}} ^\L  \widetilde{\D} _{X \rightarrow Y}      )   }
\ar[r] ^-{u _{0*}}
\ar[u] ^-\sim
&
{\R \mathcal{H} om _{\widetilde{\D} _{Y}} (   \tilde{u} _{0+}( \widetilde{\M} ),
 \tilde{u} _{0+}(\widetilde{\omega} ^{(m')} _{X} \otimes _{\widetilde{\B} _{X}} \widetilde{\D} _{X \rightarrow Y})  )   }
\ar[u] ^-\sim
}
\end{equation}
\normalsize
est commutatif.
\end{prop}
\begin{proof}
 La commutativité du carré du haut de \ref{24niveaubis} se prouve par fonctorialité.
Traitons à présent celle du bas.
On vérifie avec \ref{chgtniveu+1} que l'on dispose du diagramme commutatif canonique :
\begin{equation}
  \label{chgtniveu+1bis2}
\xymatrix @R=0,3cm {
{u _{0+} ( \M )}
\ar@{=}[r] 
&
{u _{0*} ( \M \otimes ^\L _{\D _{X}} \D _{X \rightarrow Y})}
\ar[r] \ar[d] \ar[rd] 
&
{u _{0*} ( \M \otimes ^\L  _{\D _{X}} \D _{X \rightarrow Y}) \otimes _{\D _Y } ^\L \widetilde{\D} _Y}
\ar[d] ^{\mathrm{proj}} _-\sim
\\
{\tilde{u} _{0+} ( \widetilde{\M} )}
\ar@{=}[r]
&
{u _{0*} (\widetilde{\M} \otimes _{\widetilde{\D} _{X}} ^\L \widetilde{\D} _{X \rightarrow Y})}
&
{u _{0*} ( \M \otimes _{\D _{X}} ^\L
\D _{X \rightarrow Y} \otimes _{u _0 ^{-1}\D _Y } ^\L  u _0 ^{-1} \widetilde{\D} _Y),}
\ar[l] ^-\sim
}
\end{equation}
dont l'isomorphisme du bas est l'image par $u _{0*}$ de \ref{chgtniveu+1bis}, la flèche de gauche est
le morphisme de changement de niveaux de l'image directe de \ref{chgtniveu+1}.
Pour vérifier la commutativité du carré du bas de \ref{24niveaubis}, 
il suffit alors d'appliquer le lemme \ref{24niveaulembis} au cas où 
$\NN=\M \otimes _{\D _{X}} ^\L  \D _{X \rightarrow Y}$ et 
$\widetilde{\NN}:(\widetilde{\omega} ^{(m')} _{X} \otimes _{\widetilde{\B} _{X}} \widetilde{\D} _{X \rightarrow Y}) \otimes _{\widetilde{\D} _{X}}  ^\L  \widetilde{\D} _{X \rightarrow Y} $.

 \end{proof}

\begin{lemm}
  Soient $ \M$ (resp. $\widetilde{\M}$) un $\D _X$-module à droite (resp. $\widetilde{\D} _X$-module à droite),
  $\E$ et $\FF$ (resp. $\widetilde{\E}$ et $\widetilde{\FF}$) deux $\D _X$-modules à gauche (resp. $\widetilde{\D} _X$-modules à gauche).
  Soient $a$ : $\E \rightarrow \widetilde{\E}$, $b$ : $\FF \rightarrow \widetilde{\FF}$, $c$ : $\M \rightarrow \widetilde{\M}$ des
  morphisme $\D _X$-linéaires.
On dispose alors du diagramme commutatif suivant :
\begin{equation}
  \label{47lemm1}
\xymatrix @R=0,3cm {
{(  \M \otimes _{\B _{X}} \E ) \otimes _{\D _{X}} \FF}
\ar[r] _-\sim \ar[d]
&
{  \M \otimes _{\D _{X}}( \E \otimes _{\B _{X}} \FF )}
\ar[d]
\\
{(  \widetilde{\M} \otimes _{\widetilde{\B} _{X}} \widetilde{\E} ) \otimes _{\widetilde{\D} _X} \widetilde{\FF}}
\ar[r] _-\sim
&
{\widetilde{\M} \otimes _{\widetilde{\D} _{X}}( \widetilde{\E}  \otimes _{\widetilde{\B} _{X}} \widetilde{\FF}).}
}
\end{equation}

\end{lemm}

\begin{proof}
L'isomorphisme du haut est le composé
$(  \M \otimes _{\B _{X}} \E ) \otimes _{\D _{X}} \FF \riso
  \M \otimes _{\D _{X}} (\D _X \otimes _{\B _{X}} \E ) \otimes _{\D _{X}} \FF
  \underset{\gamma _{\E}}{\riso}
\M \otimes _{\D _{X}} (\E  \otimes _{\B _{X}} \D _X ) \otimes _{\D _{X}} \FF
\riso
\M \otimes _{\D _{X}} (\E  \otimes _{\B _{X}} \FF)$.
  Soient $m$, $f$ et $e$ des sections locales respectives de
$ \M$, $\FF$ et $\E$.
Comme $\gamma _{\E} ( 1 \otimes e) = e \otimes 1$ (voir \cite[1.3.1]{Be2}),
on calcule que $(m\otimes e) \otimes f$ s'envoie via cet isomorphisme
sur $m\otimes (e \otimes f)$. De même pour l'isomorphisme du bas. Il en résulte que quelque soit le chemin de \ref{47lemm1}
$(  \M \otimes _{\B _{X}} \E ) \otimes _{\D _{X}} \FF \rightarrow
\widetilde{\M} \otimes _{\widetilde{\D} _{X}}( \widetilde{\E}  \otimes _{\widetilde{\B} _{X}} \widetilde{\FF})$
 choisi, $(m\otimes e) \otimes f$ s'envoie sur $c(m)\otimes (a(e) \otimes b(f))$.
 \end{proof}

\begin{prop}
Soient $ \M \in D ^- (\D _{X} \overset{^{\mathrm{d}}}{})$
et $\widetilde{\M} :=  \M \otimes ^{\L} _{\D _{X}} \widetilde{\D} _X$.
Le diagramme commutatif ci-contre
\begin{equation}
  \label{47niveau}
\xymatrix @R=0,3cm @C=0,5cm {
{\R \mathcal{H} om _{\D _{Y}} ( u _{0+} ( \M), u _{0+}( \widetilde{\omega} ^{(m)} _{X} \otimes _{\B _{X}} \D _{X \rightarrow Y}))}
\ar[r] _-\sim \ar[d]
&
{\R \mathcal{H} om _{\D _{Y}} ( u _{0+} ( \M),  u _{0+} ( \widetilde{\omega} ^{(m)} _{X} )\otimes _{\B _Y} \D _{Y})}
\ar[d]
\\
{\R \mathcal{H} om _{\D _{Y}}
( u _{0+} ( \M), \widetilde{u} _{0+} ( \widetilde{\omega} ^{(m')} _{X} \otimes _{\widetilde{\B} _{X}} \widetilde{\D} _{X \rightarrow Y}))}
\ar[r] _-\sim
&
{\R \mathcal{H} om _{\D _{Y}}
( u _{0+} ( \M),  \widetilde{u} _{0+} ( \widetilde{\omega} ^{(m')} _X )\otimes _{\widetilde{\B} _{Y}} \widetilde{\D} _{Y})}
\\
{\R \mathcal{H} om _{\widetilde{\D} _Y}
( \widetilde{u} _{0+} (\widetilde{\M}), \widetilde{u} _{0+} ( \widetilde{\omega} ^{(m')} _{X} \otimes _{\widetilde{\B} _{X}} \widetilde{\D} _{X \rightarrow Y}))}
\ar[r] _-\sim \ar[u] _-\sim
&
{\R \mathcal{H} om _{\widetilde{\D} _Y}
( \widetilde{u} _{0+} (\widetilde{\M}),  \widetilde{u} _{0+} ( \widetilde{\omega} ^{(m')} _X )\otimes _{\widetilde{\B} _{Y}} \widetilde{\D} _{Y}).}
\ar[u] _-\sim
}
\end{equation}
\normalsize
est commutatif.
\end{prop}
\begin{proof}
En résolvant $\D _{X \rightarrow Y}$
par des $(\D _{X}, u _0 ^{-1} \D _{Y} )$-bimodules plats,
$\widetilde{\D} _{X \rightarrow Y}$
par des $(\widetilde{\D} _{X}, u _0 ^{-1} \widetilde{\D} _{Y} )$-bimodules plats (et ce de façon compatible),
il découle de \ref{47lemm1} (qui implique aussi que $-\otimes _{\B _X}-$ préserve la $\D _X$-platitude à gauche) la commutativité du carré :
\begin{equation}
  \label{47-1}
\xymatrix @R=0,3cm {
{ \R u _{0*} ( ( \widetilde{\omega} ^{(m)} _{X} \otimes _{\B _{X}} \D _{X \rightarrow Y}) \otimes _{\D _{X}} ^\L \D _{X \rightarrow Y})}
\ar[r]_-\sim  \ar[d]
&
{ \R u _{0*} ( \widetilde{\omega} ^{(m)} _{X} \otimes ^\L  _{\D _{X}} ( \D _{X \rightarrow Y} \otimes _{\B _X} \D _{X \rightarrow Y}))}
\ar[d]
\\
{ \R u _{0*} ( ( \widetilde{\omega} ^{(m')} _X \otimes _{\widetilde{\B} _X}
\widetilde{\D} _{X \rightarrow Y}) \otimes _{\widetilde{\D} _X} ^\L  \widetilde{\D} _{X \rightarrow Y})}
\ar[r] _-\sim
&
{ \R u _{0*} ( \widetilde{\omega} ^{(m')} _X \otimes ^\L  _{\widetilde{\D} _X}
( \widetilde{\D} _{X \rightarrow Y} \otimes _{\widetilde{\B} _{X}} \widetilde{\D} _{X \rightarrow Y})).}
}
\end{equation}

De plus, comme l'isomorphisme de transposition est compatible au changement de niveaux, on obtient la commutativité du carré gauche du diagramme ci-après
\begin{equation}
  \label{47-2}
\xymatrix @C=0,3cm  @R=0,3cm {
  { \R u _{0*} ( \widetilde{\omega} ^{(m)} _{X} \otimes _{\D _{X}} ^\L ( \D _{X \rightarrow Y} \otimes _{\B _X} \!\! \D _{X \rightarrow Y}))}
\ar[r] _-\sim\ar[d]
&
{ \R u _{0*}  ( \widetilde{\omega} ^{(m)} _{X} \otimes _{\D _{X}} ^\L  (\D _{X \rightarrow Y}
\underset{u ^{\text{-1}}\B _Y}{\otimes} \!\!\!\!
u ^{\text{-1}} \D _{Y}))}
\ar[d]
&
{ \R u _{0*}  ( \widetilde{\omega} ^{(m)} _{X} \otimes _{\D _{X}} ^\L  \!\! \D _{X \rightarrow Y} )\otimes _{\B _Y} \D _{Y}}
\ar[d] \ar[l] _-\sim ^-{\mathrm{proj}}
\\
{ \R u _{0*} ( \widetilde{\omega} ^{(m')} _X \otimes _{\widetilde{\D} _X} ^\L
( \widetilde{\D} _{X \rightarrow Y} \otimes _{\widetilde{\B} _X}\!\!  \widetilde{\D} _{X \rightarrow Y}))}
\ar[r] _-\sim
&
{ \R u _{0*}  ( \widetilde{\omega} ^{(m')} _X \otimes _{\widetilde{\D} _X} ^\L (\widetilde{\D} _{X \rightarrow Y} \!\!\!\!
\underset{u ^{\text{-1}}\widetilde{\B} _Y}{\otimes} \!\!\!\!
 u ^{\text{-1}} \widetilde{\D} _{Y}))}
&
{ \R u _{0*}  ( \widetilde{\omega} ^{(m')} _X \otimes _{\widetilde{\D} _X}^\L  \!\!  \widetilde{\D} _{X\rightarrow Y} )
\otimes _{\widetilde{\B} _{Y}} \widetilde{\D} _{Y}.}
\ar[l] _-\sim ^-{\mathrm{proj}}
}
\end{equation}
\normalsize
Le carré de droite est commutatif par fonctorialité. 
En mettant bout à bout \ref{47-1} et \ref{47-2}, on obtient la commutativité de
\begin{equation}
  \label{47niveaupre}
\xymatrix @R=0,3cm {
{ u _{0+}( \widetilde{\omega} ^{(m)} _{X} \otimes _{\B _{X}} \D _{X \rightarrow Y}) }
\ar[r] _-\sim \ar[d]
&
{ u _{0+} ( \widetilde{\omega} ^{(m)} _{X} )\otimes _{\B _Y} \D _{Y}}
\ar[d]
\\
{ \widetilde{u} _{0+}( \widetilde{\omega} ^{(m')} _{X} \otimes _{\widetilde{\B} _{X}} \widetilde{\D} _{X \rightarrow Y}) }
\ar[r] _-\sim
&
{ \widetilde{u} _{0+} ( \widetilde{\omega} ^{(m')} _X )\otimes _{\widetilde{\B} _{Y}} \widetilde{\D} _{Y}.}
}
\end{equation}

En appliquant $\R \mathcal{H} om _{\D _{Y}} ( u _{0+} ( \M), -)$
à \ref{47niveaupre}, on obtient le carré du haut de \ref{47niveau}.
Comme le carré du bas de \ref{47niveau} est commutatif par fonctorialité, \ref{47niveau} est commutatif.
 \end{proof}

\begin{theo}
  \label{theo17niveau}
  On suppose que $u _{0}$ est une immersion fermée. 
  Soient $ \M \in D ^- (\D _{X} \overset{^{\mathrm{d}}}{})$
et $\widetilde{\M} :=  \M \otimes ^{\L} _{\D _{X}} \widetilde{\D} _X$.
On dispose du diagramme commutatif suivant :
\begin{equation}
  \label{17niveau}
\xymatrix @R=0,3cm {
{u _{0+} \DD ( \M )}
\ar[r] \ar[d]
&
{\R \mathcal{H} om _{\D _{Y}} ( u _{0+} ( \M),  u _{0+} ( \widetilde{\omega} ^{(m)} _{X} )\otimes _{\B _Y} \D _{Y}) [ d_X]}
\ar[d]
\\
{\widetilde{u} _{0+}\widetilde{\DD} ( \widetilde{\M})}
\ar[r]
&
{\R \mathcal{H} om _{\widetilde{\D} _Y}
( \widetilde{u} _{0+} (\widetilde{\M}),  \widetilde{u} _{0+} ( \widetilde{\omega} ^{(m')} _X )\otimes _{\widetilde{\B} _{Y}} \widetilde{\D} _{Y})[ d_X].}
}
\end{equation}

\end{theo}

\begin{proof}
Il suffit de composer \ref{rema12niveau}, \ref{12niveau}, \ref{24niveau}, \ref{24niveaubis}, \ref{47niveau}.
 \end{proof}

\begin{vide}
Avec les notations \ref{theo17niveau},
plaçons-nous dans le {\og cas relevable\fg}. Supposons $u$ propre, $u _{0} ^{-1}\B _{Y}$ et $\O _{X}$ Tor-indépendants.
Virrion a construit le morphisme trace 
$\mathrm{Tr} _+$ :
$u _{+}  (\omega _{X} ) [d _{X}] \rightarrow \omega _{Y} [d _{Y}]$ et vérifié sa compatibilité au changement de niveaux (voir \cite{Vir04}).
Il en résulte qu'il en est de même du morphisme déduit par extension
 $\mathrm{Tr} _+$ : $u _{+} ( \widetilde{\omega} ^{(m)} _{X} ) [d _X]\rightarrow \widetilde{\omega} ^{(m)} _{Y} [d _Y ]$ (voir \ref{1.5.2}). 
On en déduit que le carré supérieur du diagramme ci-après
\begin{equation}
  \label{78niveau}
\xymatrix @R=0,3cm {
{\R \mathcal{H} om _{\D _{Y}} ( u _{+} ( \M),  u _{+} ( \widetilde{\omega} ^{(m)} _{X} )\otimes _{\B _Y} \D _{Y})[d _X]}
\ar[r] ^-{\mathrm{Tr} _+} \ar[d]
&
{\R \mathcal{H} om _{\D _{Y}} ( u _{+} ( \M),   \widetilde{\omega} ^{(m)} _Y \otimes _{\B _Y} \D _{Y})[d _Y ]}
\ar[d]
\\
{\R \mathcal{H} om _{\D _{Y}}
( u _{+} ( \M),  \widetilde{u} _{+}( \widetilde{\omega} ^{(m')} _X )\otimes _{\widetilde{\B} _{Y}} \widetilde{\D} _{Y})[d _X]}
\ar[r]^-{\mathrm{Tr} _+}
&
{\R \mathcal{H} om _{\D _{Y}}
( u _{+} ( \M),  \widetilde{\omega} ^{(m')} _Y \otimes _{\widetilde{\B} _{Y}} \widetilde{\D} _{Y})[d _Y ]}
\\
{\R \mathcal{H} om _{\widetilde{\D} _Y}
( \widetilde{u} _{+}(\widetilde{\M}),  \widetilde{u} _{+}( \widetilde{\omega} ^{(m')} _X )\otimes _{\widetilde{\B} _{Y}} \widetilde{\D} _{Y})[d _X]}
\ar[r] ^-{\mathrm{Tr} _+}
\ar[u] _-\sim
&
{\R \mathcal{H} om _{\widetilde{\D} _Y}
( \widetilde{u} _{+}(\widetilde{\M}),  \widetilde{\omega} ^{(m')} _Y \otimes _{\widetilde{\B} _{Y}} \widetilde{\D} _{Y})[d _Y ]}
\ar[u] _-\sim
}
\end{equation}
est commutatif. 
La commutativité du carré inférieur étant fonctorielle, il en dérive celle de \ref{78niveau}.
En composant \ref{17niveau} et \ref{78niveau}, on obtient le carré de gauche du diagramme
\begin{equation}
  \label{18niveau}
\xymatrix @R=0,3cm @C=0,4cm {
{u _{+} \DD ( \M)}
\ar[r] \ar[d]
&
{\R \mathcal{H} om _{\D _{Y}} ( u _{+} ( \M),   \widetilde{\omega} ^{(m)} _Y \otimes _{\B _Y} \D _{Y})[d _Y]}
\ar[d] \ar[r] _-\sim
&
{\DD u _{+}  ( \M)}
\ar[d]
\\
{\widetilde{u} _{+} \widetilde{\DD} (\widetilde{\M} )}
\ar[r]
&
{\R \mathcal{H} om _{\widetilde{\D} _Y}
( \widetilde{u} _{+}(\widetilde{\M}),  \widetilde{\omega} ^{(m')} _Y \otimes _{\widetilde{\B} _{Y}} \widetilde{\D} _{Y})[d _Y]}
\ar[r] _-\sim
&
{\widetilde{\DD}  \widetilde{u} _{+} (\widetilde{\M} ),}
}
\end{equation}
dont les flèches horizontales de droite sont induites par les isomorphismes de transposition.
Il dérive de \ref{transpchgtniv} la commutativité du carré de droite de \ref{18niveau}.
Lorsque $\M\in D ^\mathrm{b} _\mathrm{coh} (\D _{X} \overset{^{\mathrm{d}}}{})$,
la commutativité du diagramme \ref{18niveau} indique que l'isomorphisme de dualité relative
$\chi $ : $u _{+} \DD ( \M) \riso \DD u _{+}  ( \M)$ est compatible au changement de niveau.
\end{vide}

\begin{vide}
Supposons que $u _{0}$ soit une immersion fermée. 
Soient $m'\geq m$ deux entiers, $ \M ' \in D ^- (\D _{X' } ^{(m)} \overset{^{\mathrm{d}}}{})$
et $\NN ' :=  \M ' \otimes ^{\L} _{\D _{X'} ^{(m)} } \D ^{(m')}  _{X'}$.
On dispose du diagramme commutatif :
\small
\begin{equation}
  \label{cubecgtnivfrob}
\xymatrix @R=0,3cm @C=-1,3cm{
&
{u _{0+} ^{(m'+s)} \DD  ^{(m'+s)} ( F ^\flat \NN ')[-d _{X}]}
\ar[rr]
\ar'[d][dd]
&&
{\R  \underset{\D _{Y} ^{(m' +s)}}{\mathcal{H} om} (
u _{0+} ^{(m'+s)} F ^\flat \NN ' , u _{0+} ^{(m'+s)} (\widetilde{\omega} ^{(m'+s)} _X) \underset{\B _Y}{\otimes} \D _{Y} ^{(m' +s)})}
\ar[dd]
\\
{u _{0+} ^{(m+s)} \DD  ^{(m+s)} ( F ^\flat \M ')[-d _{X}]}
\ar[rr] \ar[dd] \ar[ur]
&&
{\R  \underset{\D _{Y} ^{(m +s)}}{\mathcal{H}om} (
u _{0+} ^{(m+s)} F ^\flat \M ' , u _{0+} ^{(m+s)} (\widetilde{\omega} ^{(m+s)} _X) \underset{\B _Y}{\otimes} \D _Y ^{(m +s)})}
\ar[dd] \ar[ur]
\\
&
{F ^\flat  u _{0+} ^{\prime (m')} \DD  ^{(m')} ( \NN ')[-d _{X}]}
\ar'[r][rr]
&&
{F ^\flat \R  \underset{\D _{Y'} ^{(m' )}}{\mathcal{H}om} (
u _{0+} ^{\prime (m')} \NN ' , u _{0+} ^{\prime (m')} (\widetilde{\omega} ^{(m')} _{X'}) \underset{\B _{Y'}}{\otimes} \D _{Y'} ^{(m' )})}
\\
{F ^\flat  u _{0+} ^{\prime (m)} \DD  ^{(m)} ( \M ')[-d _{X}]}
\ar[rr] \ar[ur]
&&
{F ^\flat \R  \underset{\D _{Y'} ^{(m )}}{\mathcal{H}om} (
u _{0+} ^{\prime (m)} (  \M ') , u _{0+} ^{\prime (m)} (\widetilde{\omega} ^{(m)} _{X'}) \underset{\B _{Y'}}{\otimes} \D _{Y'} ^{(m )}),}
\ar[ur]
}
\end{equation}
\normalsize
En effet, la commutativité des faces horizontales résultent de \ref{theo17niveau},
celle de devant et de derrière découlent de \ref{theo1-7}.
Pour établir la comutativité du carré de gauche de \ref{cubecgtnivfrob}, il s'agit de vérifier que les isomorphismes de commutation à Frobenius de l'image directe et du foncteur
dual sont compatibles aux morphismes de changement de niveau.
Or, en reprenant la construction de l'isomorphisme $\B _{X} ^{(m+s)} \otimes _{\O _{X}}\D _X ^{(m+s)}\riso F ^* F ^\flat (\B _{X'} ^{(m)} \otimes _{\O _{X'}} \D _{X'} ^{(m)})$ de \cite[2.5.2]{Be2},
on vérifie que la commutativité du diagramme canonique 
\begin{equation}
\label{FFchgtniv}
\xymatrix @R=0,3cm {
{\B _{X} ^{(m+s)} \otimes _{\O _{X}}\D _X ^{(m+s)}}
\ar[r] ^-{\sim} \ar[d]
&
{F ^* F ^\flat (\B _{X'} ^{(m)} \otimes _{\O _{X'}} \D _{X'} ^{(m)})}
\ar[d]
\\
{\B _{X} ^{(m'+s)} \otimes _{\O _{X}}\D _X ^{(m'+s)}}
\ar[r] ^-{\sim} 
&
{F ^* F ^\flat (\B _{X'} ^{(m')} \otimes _{\O _{X'}} \D _{X'} ^{(m')}).}
}
\end{equation}
On en déduit que l'isomorphisme de commutation à Frobenius du foncteur dual (voir \ref{defDualcomFrob}) est compatible au changement de niveau. Il résulte de la commutativité de \ref{FFchgtniv} que l'isomorphisme $F ^\flat  \M ' \otimes _{\B _{X} ^{(m+s)} \otimes _{\O _{X}}\D _X ^{(m+s)}} F ^* \E ' \riso  \M ' \otimes _{\B _{X'} ^{(m')} \otimes _{\O _{X'}} \D _{X'} ^{(m')}} \E '$ fonctoriel en $\M '\in D ^- (\B _{X'} ^{(m')} \otimes _{\O _{X'}} \D _{X'} ^{(m')}\overset{^{\mathrm{d}}}{})$
et $\E '\in D ^- (\overset{^{\mathrm{g}}}{} \B _{X'} ^{(m')} \otimes _{\O _{X'}} \D _{X'} ^{(m')} )$
est compatible au changement de niveaux (pour le voir, il s'agit de reprendre sa construction dans \cite[2.5.7]{Be2}). 
De plus, comme l'isomorphisme de projection est compatible au changement de niveaux, il en est alors de même pour l'isomorphisme de commutation de l'image directe à Frobenius (construit en \ref{u+comFrob}). 
On procède de même pour vérifier la commutativité du diagramme de droite. 
\end{vide}

\section{\label{cf}Isomorphisme de dualité relative: cas des schémas formels}

Soit $\V$ un anneau de valuation discrète complet d'inégales caractéristiques $(0,\,p)$,
de corps résiduel $k$.
On fixe $\pi$ une uniformisante de $\V$, $s$ un entier et
$\sigma$ : $\V \riso \V$ un relèvement de la puissance $s$-ème de l'automorphisme de Frobenius de $k$.
Conformément aux notations du chapitre \ref{ca}, pour tout $k$-schéma lisse $X$,
$F$ sera le Frobenius relatif $F ^s _{X /S }$ où $S =\Spec k$.

Lorsque nous considérerons des produits de
$\V$-schémas formels (resp. $k$-schémas), nous omettrons d'indiquer $\Spf \V$ (resp. $\Spec k$).
Enfin, pour tout faisceau abélien $\M$, on posera $\M _\Q := \M \otimes _{\Z} \Q$.

\subsection{Notations}
Nous garderons les notations suivantes : Soient $\X$ et $\Y$ deux $\V$-schémas formels lisses,
$u _0$ : $X \rightarrow Y$ un morphisme propre entre les fibres spéciales,
$X _{i}$ et $Y _{i}$ les schémas déduits respectivement de $\X$ et $\Y$ par réduction modulo $\pi ^{i+1}$ (e.g. $X _{0} =X$),
$T$ un diviseur de $Y$ tel que $T _X := u _0 ^{-1} (T) $ soit un diviseur de $X$.
On désigne par
$\X '$ (resp. $\Y '$) le $\V$-schéma formel déduit par changement de base par $\sigma$ de $\X$ (resp. $\Y$),
$ u' _0$ : $ X' \hookrightarrow Y'$ l'immersion fermée induite,
$T'$ et $T _{X'}$ les images inverses par $\sigma$ respectives de
$T$ et $T _X$ sur $\Y'$ et $\X'$.

Les définitions et notations concernant les complexes quasi-cohérents seront celles de
\cite[3.2 et 4.2]{Beintro2}.

\begin{vide}

$\bullet $ Pour tout entier $m$, notons
$\widetilde{\omega} _{\X} ^{(m)} :=\omega _{\X} \otimes _{\O _{\X}} \widehat{\mathcal{B}} ^{(m)} _{\X} (T_X)$,
$\widetilde{\omega} _{\Y} ^{(m)} :=\omega _{\X} \otimes _{\O _{\Y}} \widehat{\mathcal{B}} ^{(m)} _{\Y} (T)$,
$\omega _{\X} (\hdag T _X) _{\Q}:=\omega _{\X,\Q} \otimes _{\O _{\X,\Q}} \O _{\X} (\hdag T _X) _{\Q}$
et
$\omega _{\Y} (\hdag T) _{\Q}:=\omega _{\Y,\Q} \otimes _{\O _{\Y,\Q}} \O _{\Y} (\hdag T ) _{\Q}$,
où $\widehat{\B} ^{(m)}  _{\X} ( T _X)$  et $\widehat{\B} ^{(m)}  _{\Y} ( T )$ ont été construits en \cite[4.2.4]{Be1}.
On notera $\smash{\widetilde{\D}} _{\X} ^{(m)}:=
\widehat{\B} ^{(m)}  _{\X} (T _X) \widehat{\otimes} _{\O _{\X}}\widehat{\D} ^{(m)} _{\X}$,
$\smash{\widetilde{\D}} _{\Y} ^{(m)}:=
\widehat{\B} ^{(m)}  _{\Y} (T ) \widehat{\otimes} _{\O _{\Y}}\widehat{\D} ^{(m)} _{\Y}$,
$\smash{\widetilde{\D}} _{X _{i}} ^{(m)}:=
\B ^{(m)}  _{X _{i}} (T _X) \otimes _{\O _{\X}} \D ^{(m)} _{X _{i}}$,
$\smash{\widetilde{\D}} _{Y _{i}} ^{(m)}:=
\B ^{(m)}  _{Y _{i}} (T ) \otimes _{\O _{\X}} \D ^{(m)} _{Y _{i}}$. 
Soient
$\M ^{(m)} \in D ^{\mathrm{b}} _{\mathrm{coh}}
(\smash{\widetilde{\D}} _{\X} ^{(m)}   \overset{^\mathrm{d}}{})$
et $ \M ^{(m)} _i :=  \M ^{(m)} \otimes ^\L _{\widetilde{\D} ^{(m)} _{\X}} \widetilde{\D} ^{(m)} _{X _i}$.
On définit,
de manière analogue à \cite[1.1]{caro_courbe-nouveau},
{\it l'image directe de $\M ^{(m)} $ par $u _0$ à singularités surconvergentes le long de $T$ de niveau $m$} en posant :
\begin{equation}\label{utmhat+}
u  _{0T+} ^{(m)} ( \M^{(m)} ):=
\R u _{0*} (\M ^{(m)} \otimes ^{\L}
_{ \smash{\widetilde{\D}} _{\X} ^{(m)} }
(\widehat{\B} ^{(m)}  _{\X} (T _X) \widehat{\otimes} _{\O _{\X}}\widehat{\D} ^{(m)} _{\X \rightarrow \Y} )).
\end{equation}
Conformément aux notations du précédent chapitre, l'image directe de $\M ^{(m)} _{i}$ par $u _0$ à singularités surconvergentes le long de $T$ de niveau $m$ est défini en posant :
$u  _{0T+} ^{(m)} ( \M^{(m)} _{i} ):=
\R u _{0*} (\M ^{(m)} _{i} \otimes ^{\L}
_{ \smash{\widetilde{\D}} _{X _{i} } ^{(m)} }
(\B ^{(m)}  _{X _{i} } (T _X) \otimes _{\O _{X _{i} }} \D ^{(m)} _{X _{i} \rightarrow Y _{i}} )).$
De manière analogue à \cite[3.5.5]{Beintro2}, on dispose aussi de l'isomorphisme canonique :
$$ u ^{(m)} _{0T+} (\M ^{(m)}) \riso \R \underset{\underset{i}{\longleftarrow}}{\lim} \, u ^{(m)} _{0T+}  (\M ^{(m)} _i) .$$

Le {\it foncteur dual sur $\X$ de niveau $m$ et à singularités surconvergentes le long de $T_X$} de $\M ^{(m)} $ se note :
\begin{equation}\label{Dtdag}
 \DD _{\X, T _X} ^{(m)} ( \M ^{(m)}) :=\widetilde{\omega} _{\X} ^{(m)} 
 \otimes _{\widehat{\mathcal{B}} ^{(m)} _{\X} (T_X)} \R \mathcal{H}om _{ \smash{\widetilde{\D}} _{\X} ^{(m)} }
( \M ^{(m)} ,   \smash{\widetilde{\D}} _{\X} ^{(m)}  [d _{\X}]).
\end{equation}

$\bullet $ Soit
$\M \in D ^{\mathrm{b}} _{\mathrm{coh}} ( \D ^{\dag} _{\X } (\hdag T _X ) _{\Q} \overset{^\mathrm{d}}{})$.
On définit {\it l'image directe de $\M$ par $u _0$ à singularités surconvergentes le long de $T$} en posant
(on s'inspire de  \cite[1.1]{caro_courbe-nouveau} et \cite[4.3.7]{Beintro2}) :
\begin{equation}\label{utdag+}
u  _{0T+}( \M ):=
\R u _{0*} (\M \otimes ^{\L} _{ \D ^{\dag} _{\X } (\hdag T _X) _{\Q}}
\D ^{\dag} _{\X \rightarrow \Y}  ( \hdag T ) _{\Q}),
\end{equation}
où $\D ^{\dag} _{\X \rightarrow \Y}  ( \hdag T ) _{\Q}:=
\underset{\longrightarrow}{\lim} _m
\widehat{\B} ^{(m)}  _{\X} (T _X) \widehat{\otimes} _{\O _{\X}}\widehat{\D} ^{(m)} _{\X \rightarrow \Y,\Q} $.

Le {\it foncteur dual sur $\X$ à singularités surconvergentes le long de $T_X$} de $\M $ se note :
\begin{equation}\label{Dtdag}
\DD _{\X, T _X} ( \M) := 
\omega _{\X} (\hdag T _X) _{\Q} 
\otimes _{\O _{\X} (\hdag T _X) _{\Q}} 
\R \mathcal{H}om _{\D ^{\dag} _{\X } (\hdag T _X ) _{\Q}} ( \M ,  \D ^{\dag} _{\X } (\hdag T _X ) _{\Q}  [d _{\X}]).
\end{equation}
Si aucune confusion n'est possible, afin d'alléger on pourra aussi le noter $\DD _{T _X}$.

Enfin, l'image directe extraordinaire par $u _{0}$  à singularités surconvergentes le long de $T$ de $\M $ se note : 
\begin{equation}\label{imag-direc-extra}
u _{0T !} ( \M) =\DD _{\Y, T }  \circ u _{0 T +}  \circ \DD _{\X, T _X} ( \M).
\end{equation}

Pour toutes les opérations cohomologiques définies ci-dessus,
lorsque le diviseur $T$ est vide, il sera omis.

\end{vide}

\subsection{Compatibilité à Frobenius du morphisme trace dans le cas d'une immersion fermée relevable \label{subsecfrobtr}}

On suppose dans cette section que $u _{0}$ se relève en une immersion fermée $u\,:\,\X \hookrightarrow \Y$. 
Nous prouvons ici que le morphisme trace de Virrion est compatible à Frobenius. 
Pour cela, nous construisons un morphisme trace compatible à Frobenius puis nous vérifions que ces deux morphismes traces sont égaux.

On note
$u _i$ : $X _i \rightarrow Y _i$ les réductions modulo $\pi ^{i+1}$ de $u$ (de même en rajoutant des primes).
De plus, les images directes et images inverses extraordinaires
de niveau $m$ (sans diviseur) par $u$ (resp. $u _i$) seront notées
$u _+ ^{(m)}$ et $u ^{!(m)}$  (resp. $u _{i+} ^{(m)}$ et $u _i ^{!(m)}$).

\begin{nota}
Nous noterons $\R \underline{\Gamma} ^{(m)} _{X _i}$ le foncteur cohomologique local défini par Berthelot (voir \cite[4.4.4]{Beintro2}. Remarquons que comme $X _{i}$ est lisse, ce foncteur est canoniquement isomorphe
au foncteur cohomologique local défini dans \cite[2]{caro_surcoherent}.
\end{nota}

\begin{vide}\label{221}
Soient $ \M _i $ un $\D _{X _i} ^{(m)}$-module quasi-cohérent à droite
et
$ \NN _i $ un $\D _{Y _i} ^{(m)}$-module quasi-cohérent à droite.

On a l'isomorphisme canonique :
$\mathcal{H} ^0 u _i ^{!(m)} (  \NN _i) \riso u ^{-1} _0 \mathcal{H} om _{O _{Y _i}} ( u _{0*} \O _{X _i},  \NN _i)$.
On en déduit l'inclusion canonique
$\mathrm{Tr}\,:\,u _{0*} \mathcal{H} ^0 u _i ^{!(m)} (  \NN _i)
\rightarrow
\NN _i$. Lorsque $\NN _i$ est $\O _{Y_i}$-cohérent, ce dernier est {\it le morphisme trace}
(voir plus généralement le cas des morphismes finis voire propres dans \cite{HaRD}).

On note
$\mathrm{T}$ : $ u _{0*} (\M _i) \rightarrow u _{i+} ^{(m)} (\M _i)
=
u _{0*} ( \M _i \otimes _{\D _{{X _i}} ^{(m)}} u ^* _i \D ^{(m)} _{{Y _i}})$
le morphisme canonique.

Via l'isomorphisme canonique
$$u _{i+} ^{(m)} \mathcal{H} ^0 u _i ^{!(m)} (  \NN _i) \riso
u _{0*} ( u ^{-1} _0\mathcal{H} om _{O _{Y _i}} ( u _{0*} \O _{X _i},  \NN _i)
\otimes _{\D _{{X _i}} ^{(m)}} u ^* _i \D ^{(m)} _{{Y _i}}),$$
le morphisme
$\mathcal{H} om _{O _{Y _i}} ( u _{0*} \O _{X _i},  \NN _i)
\rightarrow
\underline{\Gamma} ^{(m)} _{X _i} ( \NN _i)$ se factorise
par
$u _{i+} ^{(m)} \mathcal{H} ^0 u _i ^{!(m)} (  \NN _i) \rightarrow \underline{\Gamma} ^{(m)} _{X _i} ( \NN _i)$
(on pourra consulter \cite{Becohdiff} pour la version {\og à gauche \fg}).
En outre, celui-ci s'inscrit dans le diagramme commutatif de gauche :
\begin{equation} \label{predefcanm}
  \xymatrix @R=0,3cm {
{u _{i+} ^{(m)} \mathcal{H} ^0 u _i ^{!(m)} (  \NN _i) }
\ar[r] ^-{\mathrm{can}}
\ar[dr] ^-{\mathrm{can}}
&
{\underline{\Gamma} ^{(m)} _{X _i} ( \NN _i)}
\ar[d] ^-{\mathrm{can}}
\\
{u _{0*} \mathcal{H} ^0 u _i ^{!(m)} (  \NN _i) }
\ar[r] ^-{\mathrm{Tr}} \ar[u] ^-{\mathrm{T}}
&
{ \NN _i.}
}
\end{equation}
\end{vide}

\begin{vide}
On identifiera $\M _i$ et $\mathcal{H} ^0 u _i ^{!(m)} u _{0*} (\M _i)$
via l'isomorphisme 
$\mathcal{H} ^0 u _i ^{!(m)} u _{0*} (\M _i) \riso \M _i$.
Avec cette identification, 
$\mathcal{H} ^0 u _i ^{!(m)} (\mathrm{Tr})$ :
$\mathcal{H} ^0 u _i ^{!(m)} u _{0*} \mathcal{H} ^0 u _i ^{!(m)} (  \NN _i)
\rightarrow
\mathcal{H} ^0 u _i ^{!(m)} (\NN _i)$ (resp.
$\mathcal{H} ^0 u _i ^{!(m)} (\mathrm{T})$ :
$\mathcal{H} ^0 u _i ^{!(m)} u _{0*} (\M _i)
\rightarrow \mathcal{H} ^0 u _i ^{!(m)} u _{i+} ^{(m)} (\M _i)$)
devient l'identité
de $\mathcal{H} ^0 u _i ^{!(m)} (\NN _i)$
(resp. le morphisme canonique
$\M _i \rightarrow \mathcal{H} ^0 u _i ^{!(m)} u _{i+} ^{(m)} (\M _i)$).
On calcule d'ailleurs que le morphisme canonique $\M _i \rightarrow \mathcal{H} ^0 u _i ^{!(m)} u _{i+} ^{(m)} (\M _i)$ est un isomorphisme compatible à Frobenius.

En appliquant $\mathcal{H} ^0 u _i ^{!(m)}$ à \ref{predefcanm}
(et via l'identification entre les deux termes $\mathcal{H} ^0 u _i ^{!(m)} u _{0*} (\M _i) \riso \M _i$),
on obtient le diagramme commutatif 
\begin{equation}
\label{predefcanm2}
  \xymatrix @R=0,3cm {
{\mathcal{H} ^0 u _i ^{!(m)} u _{i+} ^{(m)} \mathcal{H} ^0 u _i ^{!(m)} (  \NN _i) }
\ar[r] ^-{\mathrm{can}}
\ar[dr] ^-{\mathrm{can}}
&
{\mathcal{H} ^0 u _i ^{!(m)} \underline{\Gamma} ^{(m)} _{X _i} ( \NN _i)}
\ar[d] ^-{\mathrm{can}}
\\
{ \mathcal{H} ^0 u _i ^{!(m)} (  \NN _i) }
\ar@{=}[r]  \ar[u] ^-{\mathrm{can}} _-\sim
&
{\mathcal{H} ^0 u _i ^{!(m)} ( \NN _i).}
}
\end{equation}

\end{vide}

\begin{vide}
\label{vide-pre2defcanm}
Pour tout $\NN _i \in D ^+ _\mathrm{qc} (\D _{Y _i} ^{(m)} \overset{^d}{})$,
il dérive de \ref{predefcanm} et \ref{predefcanm2} les diagrammes commutatifs :
\begin{equation} \label{pre2defcanm}
  \xymatrix @R=0,3cm {
{u _{i+} ^{(m)} u _i ^{!(m)} (  \NN _i) }
\ar[r] ^-{\mathrm{can}} _-\sim
\ar[dr] ^-{\mathrm{can}}
&
{\R  \underline{\Gamma} ^{(m)} _{X _i} ( \NN _i)}
\ar[d] ^-{\mathrm{can}}
\\
{u _{0*}  u _i ^{!(m)} (  \NN _i) }
\ar[r] ^-{\mathrm{Tr}} \ar[u] ^-{\mathrm{T}}
&
{ \NN _i,}
}
\hfill
\hspace{1cm}
  \xymatrix @R=0,3cm {
{ u _i ^{!(m)} u _{i+} ^{(m)}  u _i ^{!(m)} (  \NN _i) }
\ar[r] ^-{\mathrm{can}} _-\sim
\ar[dr] ^-{\mathrm{can}}
&
{ u _i ^{!(m)} \R  \underline{\Gamma} ^{(m)} _{X _i} ( \NN _i)}
\ar[d] ^-{\mathrm{can}}
\\
{  u _i ^{!(m)} (  \NN _i) }
\ar@{=}[r]  \ar[u] ^-{\mathrm{can}} _-\sim
&
{u _i ^{!(m)} ( \NN _i),}
}
\end{equation}
dont les flèches horizontales du haut sont, par passage de gauche à droite de \cite[4.4.5]{Beintro2}),
des isomorphismes.
\end{vide}

\begin{vide}
\label{+!gammmfrob}
Soit $\NN \in D ^\mathrm{b} _\mathrm{coh} ( \D ^\dag _{\Y} (\hdag T) _{\Q} \overset{^{\mathrm{d}}}{})$
tel que $ u _{T} ^! (\NN) \in
D ^\mathrm{b} _\mathrm{coh} ( \D ^\dag _{\X} (\hdag T _{X}) _{\Q} \overset{^{\mathrm{d}}}{})$.
Rappelons que, d'après \cite[1.1.9 et 1.1.10]{caro_courbe-nouveau},  dans la catégorie des complexes quasi-cohérents les foncteurs $u ^! $ et $ u _{T} ^!$ (resp. $u _{+} $ et $u _{T+}$) sont (en omettant d'indiquer les foncteurs oublis) canoniquement isomorphes (sur les catégories de la forme ${\underset{^{\longrightarrow }}{LD }}  ^\mathrm{b} _{\Q, \mathrm{qc}} (\smash{\widetilde{\D} } _{\X} ^{(\bullet)} \overset{^{\mathrm{d}}}{})$ et ${\underset{^{\longrightarrow }}{LD }}  ^\mathrm{b} _{\Q, \mathrm{qc}} (\smash{\widetilde{\D} } _{\Y} ^{(\bullet)} \overset{^{\mathrm{d}}}{})$). Avec l'équivalence de catégories
$\underset{\longrightarrow}{\lim}$ :
$\smash[b]{\underset{^{\longrightarrow }}{LD }}  ^\mathrm{b} _{\Q, \mathrm{coh}}
(\smash{\widetilde{\D}} _{\Y} ^{(\bullet)} \overset{^{\mathrm{d}}}{})
\cong
D ^{\mathrm{b}} _\mathrm{coh} ( \D ^\dag _{\Y} (\hdag T) _{\Q})$ (voir \cite[4.2.4]{Beintro2} ou \cite[1.1.3]{caro_courbe-nouveau}),
on déduit alors du diagramme de gauche de \ref{pre2defcanm} 
l'isomorphisme canonique $\alpha\,:\,u _{T+} u _{T} ^! (\NN) \riso
\R \underline{\Gamma} _X ^\dag (\NN)$.
On vérifie en outre la compatibilité à Frobenius de $\alpha$ de la manière suivante :
d'après l'analogue $p$-adique de Berthelot du théorème de Kashiwara (voir \cite[5.3.3]{Beintro2}),
$\alpha$ est compatible à Frobenius si et seulement si $u _{T} ^! (\alpha)$ l'est.
Or, par complétion et passage à la limite sur le niveau, il découle du diagramme de droite de \ref{pre2defcanm}
le suivant :
\begin{equation}
  \label{pre3defcanm}
  \xymatrix @R=0,3cm {
{ u  _{T} ^!u _{T+} ^{}  u  _{T} ^!(  \NN ) }
\ar[r] ^-{u  _{T} ^!(\alpha)} _-\sim
&
{ u  _{T} ^!\R  \underline{\Gamma} ^{\dag} _{X } ( \NN )}
\ar[d] ^-{\mathrm{can}}
\\
{  u  _{T} ^!(  \NN ) }
\ar@{=}[r]  \ar[u] ^-{\mathrm{can}} _-\sim
&
{u  _{T} ^!( \NN ).}
}
\end{equation}
Par \cite[5.3.3]{Beintro2} (resp. \cite[1.2.11]{caro_surcoherent}), la flèche de gauche (resp. de droite) est compatible à Frobenius.
Comme le diagramme \ref{pre3defcanm} est commutatif, on obtient la compatibilité à Frobenius de $u _{T} ^! (\alpha)$. D'où celle de $\alpha$.
\end{vide}

\begin{vide}
\label{videdefcandag}
On dispose des isomorphismes canoniques
$u _i ^{!(m)} (\omega _{Y _i} [d _{Y _i}]) \riso \omega _{X _i} [d _{X _i}]$
compatibles à Frobenius. Il en découle 
$u  _{T} ^!(\omega _{\Y} (\hdag T) _{\Q} [d _{Y }]) \riso \omega _{\X} (\hdag T _{X}) _{\Q} [d _{X }]$.
Avec \ref{+!gammmfrob}, on obtient alors
l'isomorphisme canonique compatible à Frobenius :
$u _{T+} (\omega _{\X} (\hdag T _{X}) _{\Q} [d _X]) \riso \R \underline{\Gamma} _X ^\dag (\omega _{\Y} (\hdag T) _{\Q} [d _Y])$.
D'après \cite[1.2.11]{caro_surcoherent}, il en est de même de
$\R \underline{\Gamma} _X ^\dag (\omega _{\Y} (\hdag T) _{\Q} [d _Y])
\rightarrow \omega _{\Y} (\hdag T) _{\Q} [d _Y]$.
Le morphisme canonique $\mathrm{Tr} ' _{T+}$, que l'on définit par composition de la façon suivante :
\begin{equation}\label{can<->tr}
\mathrm{Tr} ' _{T+} \ : \
u _{T+} (\omega _{\X} (\hdag T _{X}) _{\Q} [d _X]) \riso \R \underline{\Gamma} _X ^\dag (\omega _{\Y} (\hdag T) _{\Q} [d _Y])
\rightarrow \omega _{\Y} (\hdag T) _{\Q} [d _Y],
\end{equation}
est ainsi compatible à Frobenius.

Il dérive du diagramme de gauche de \ref{pre2defcanm} par complétion, tensorisation par $\Q$ et passage
à la limite sur le niveau,
que $\mathrm{Tr} ' _{T+} $ s'inscrit dans le diagramme commutatif
\begin{equation}
\label{defcandag}
\xymatrix @R=0,3cm @C =2cm {
{u _{T+} (\omega _{\X} (\hdag T _{X}) _{\Q} )}
\ar@<1ex>[r] ^-{\mathrm{Tr} ' _{T+}}
\ar@<-1ex>[r] _-{\mathrm{Tr} _{T+}}
&
{\omega _{\Y} (\hdag T) _{\Q} [d _Y -d _X]}
\\
{u _{0*} (\omega _{\X} (\hdag T _{X}) _{\Q} ) }
\ar[r] _-{\mathrm{Tr}}
\ar[u] ^-{\mathrm{T}}
&
{\omega _{\Y} (\hdag T) _{\Q} [d _Y -d _X],}
\ar@{=}[u]
}
\end{equation}
où $\mathrm{Tr} _{T+}$ désigne le morphisme trace de \cite[1.2.5,6]{caro_courbe-nouveau} 
qui est obtenu par extension des scalaires via $\D ^\dag _{\Y,\Q} \to \D ^\dag _{\Y } (\hdag T ) _{\Q}$ à partir de celui construit par Virrion sans diviseur dans \cite{Vir04}.

\end{vide}

\begin{theo}\label{tr+compfrob}
  Avec les notations de \ref{videdefcandag}, les deux morphismes
$\mathrm{Tr} _{T+}$, $\mathrm{Tr} ' _{T+}$ : $u _{T+} (\omega _{\X} (\hdag T _{X}) _{\Q} )\rightarrow \omega _{\Y} (\hdag T) _{\Q} [d _Y -d _X]$
sont identiques. En particulier, le morphisme trace $\mathrm{Tr} _{T+}$ est compatible à Frobenius.
\end{theo}

\begin{proof}
  En appliquant $\mathcal{H} ^0 u  _{T} ^!$ à \ref{defcandag}, on obtient le diagramme commutatif :
\begin{equation}
\label{u!defcandag}
\xymatrix @R=0,3cm {
{\mathcal{H} ^0 u _{T} ^! u _{T+} (\omega _{\X} (\hdag T _{X}) _{\Q} )}
\ar@<1ex>[r] ^-{\mathcal{H} ^0 u  _{T} ^!(\mathrm{Tr} ' _{T+})} \ar@<-1ex>[r] _-{\mathcal{H} ^0 u  _{T} ^!(\mathrm{Tr} _{T+})}
&
{\mathcal{H} ^0 u  _{T} ^!(\omega _{\Y} (\hdag T) _{\Q} [d _Y-d _X])}
\\
{\mathcal{H} ^0 u _{T} ^! u _{0*} (\omega _{\X} (\hdag T _{X}) _{\Q}) }
\ar[r] _-{\mathcal{H} ^0 \mathrm{Tr}} \ar[u] ^-{\mathcal{H} ^0 \mathrm{T}}
&
{\mathcal{H} ^0 u  _{T} ^!(\omega _{\Y} (\hdag T) _{\Q} [d _Y-d _X]).}
\ar@{=}[u]
}
\end{equation}
Pour tout $\D ^\dag _{\X} (\hdag T _{X}) _{\Q}$-module (ou complexe) à droite $ \M$, on calcule
que le morphisme canonique
$\mathcal{H} ^0 u _{T} ^! u _{0*} ( \M) \rightarrow \mathcal{H} ^0 u _{T} ^! u _{T+} ( \M)$
est canoniquement isomorphe à l'identité de $ \M$.
La commutativité de \ref{u!defcandag} implique alors
que les deux morphismes
$\mathcal{H} ^0 u  _{T} ^!(\mathrm{Tr} ' _{T+})$,
$\mathcal{H} ^0 u  _{T} ^!(\mathrm{Tr} _{T+})$ :
$\mathcal{H} ^0 u _{T} ^! u _{T+} (\omega _{\X} (\hdag T _{X}) _{\Q} )
\rightarrow \mathcal{H} ^0 u  _{T} ^!(\omega _{\Y} (\hdag T) _{\Q}) [d _Y-d _X]$
sont égaux.
Dans cette égalité, le symbole $\mathcal{H} ^0$ étant inutile (pour le terme de gauche, cela découle de \cite[5.3.3]{Beintro2}),
on obtient $u _{T} ^! (\mathrm{Tr} ' _{T+}) =u  _{T} ^!(\mathrm{Tr} _{T+})$.

Il résulte alors du lemme \ref{lem-tr+compfrob} ci-après que $\mathrm{Tr} ' _{T+}=\mathrm{Tr} _{T+}$.
 \end{proof}

\begin{lemm}
\label{lem-tr+compfrob}
  Soient $\E \in D ^\mathrm{b} _\mathrm{coh} (\D ^\dag _{\X} (\hdag T _{X}) _{\Q})$,
$\FF  \in D ^\mathrm{b} _\mathrm{coh} (\D ^\dag _{\Y} (\hdag T) _{\Q})$ tel que $u _{T} ^! (\FF ) \in
D ^\mathrm{b} _\mathrm{coh} (\D ^\dag _{\X} (\hdag T _{X}) _{\Q})$.
L'application $\mathrm{Hom}  _{D ^\mathrm{b} _\mathrm{coh} (\D ^\dag _{\Y} (\hdag T) _{\Q})}
(u _{T+} (\E), \FF ) \rightarrow
\mathrm{Hom}  _{D ^\mathrm{b} _\mathrm{coh} (\D ^\dag _{\X} (\hdag T _{X}) _{\Q})} ( u_{T} ^! u _{T+} (\E), u _{T} ^! (\FF ))$
induite par $u  _{T} ^!$ est injective.
\end{lemm}
\begin{proof}
Soient $f _1, f _2$ : $u _{T+} (\E) \rightarrow \FF $ deux morphismes $\D ^\dag _{\Y} (\hdag T) _{\Q}$-linéaires
tels que $u _{T} ^! (f _1) =u _{T} ^! (f _2)$.
D'après le théorème de Berthelot-Kashiwara, le morphisme canonique 
$u _{T+} u _{T} ^! u _{T+} (\E) \to u _{T+} (\E) $ est un isomorphisme.
Pour tout $i=1,2$, on obtient par fonctorialité
$f _i$ : $u _{T+} (\E) \liso u _{T+} u _{T} ^! ( u _{T+} (\E))
\overset{u _{T+} u _{T} ^!  (f_i)}{\longrightarrow}
u _{T+} u _{T} ^!(\FF )\to \FF$.
Comme $u _{T} ^! (f _1)= u _{T} ^! (f _2)$, il en découle $f _1 =f _2$.
 \end{proof}

\subsection{Construction et compatibilité à Frobenius du morphisme trace dans le cas d'une immersion fermée non-relevable\label{subsecfrobtr-sing}}

Dans cette section, on suppose que $u _{0}$ est une immersion fermée.

\begin{vide}
[Isomorphisme de recollement pour l'image directe]
\label{tau-u+}
Soit
$\M \in D ^{\mathrm{b}} _{\mathrm{coh}} ( \D ^{\dag} _{\X } (\hdag T _X ) _{\Q} \overset{^\mathrm{d}}{})$.
Il existe une base d'ouverts $\widetilde{\Y}$ de $\Y$
de fibres spéciales $\widetilde{Y}$
telle que, en notant $\widetilde{\X}$ l'ouvert de $\X$ d'espace topologique sous-jacent $u _0 ^{-1} (\widetilde{Y})$
et $\widetilde{X}$ sa fibre spéciale,
le morphisme
$\widetilde{X} \rightarrow \widetilde{Y}$ induit par $u _0$ se relève en un morphisme
$\widetilde{u}$ : $\widetilde{\X} \rightarrow \widetilde{\Y}$
(en effet, il suffit de prendre une base d'ouverts affines de $\Y$).
En notant $\widetilde{T} :=\widetilde{Y} \cap T$,
l'isomorphisme canonique
$\D ^{\dag} _{\X \rightarrow \Y}  ( \hdag T ) _{\Q} | {\widetilde{\X}}
\riso \D ^{\dag} _{\widetilde{\X} \rightarrow \widetilde{\Y}}  ( \hdag \widetilde{T} ) _{\Q}$ induit le suivant
\begin{equation}
  \label{u0+=u+dag}
u  _{0T+}( \M ) | {\Y'} \riso \widetilde{u} _{\widetilde{T}+} (\M | {\X'}),
\end{equation}
où le terme de droite (i.e., l'image directe dans le cas relevable) est défini dans \cite[1.1]{caro_courbe-nouveau}.
Si $\widetilde{v}$ : $\widetilde{\X} \rightarrow \widetilde{\Y}$ est un second relèvement de
$\widetilde{X} \rightarrow \widetilde{Y}$,
il en dérive un isomorphisme
\begin{equation}
  \label{u0+=u+dagrecol}
\tau _{\widetilde{v},\widetilde{u}} \ :\
\widetilde{u} _{\widetilde{T}+} (\M | {\X'}) \riso \widetilde{v} _{\widetilde{T}+}(\M | {\X'})
\end{equation}
vérifiant la condition de transitivité
$\tau _{\widetilde{w},\widetilde{v}} \circ \tau _{\widetilde{v},\widetilde{u}} = \tau _{\widetilde{w},\widetilde{u}}$,
où $\widetilde{w}$ est un troisième relèvement.
\end{vide}

\begin{vide}
\label{u-u'-tau}
On suppose que $u _{0}$ se relève en deux immersions fermées $u, u'\, : \, \X \rightarrow \Y$. On note $\iota \, :\, X \hookrightarrow \X $, $\rho \, : \, X \hookrightarrow \Y \times \Y$ les morphismes canoniques. Soit $ \NN _i $ un $\D _{Y _i} ^{(m)}$-module à droite quasi-cohérent injectif.
Par construction de l'isomorphisme canonique $\tau\, : \ u ^{!} _{i} (\NN _{i}) \riso u ^{\prime !} _{i} (\NN _{i})$ (de manière analogue à \cite[2.1.5]{Be2}), on obtient la commutativité du carré droite du diagramme canonique : 
\begin{equation}
\label{diag-comm-tau-rho}
\xymatrix @R=0,3cm { 
{ u ^{!} _{i} (\NN _{i})} 
\ar[d] ^-{\sim} _-{\tau} 
&
{ \iota ^{!} (u ^{!} _{i} (\NN _{i})) }
\ar@{_{(}->}[l]
\ar[d] ^-{\sim} _-{\tau} 
&
{ \rho ^{!}  (\NN _{i}) }
\ar[l] ^-{\sim}
\ar@{=}[d]
\\
{ u ^{\prime !} _{i} (\NN _{i})}
&
{ \iota ^{!} (u ^{\prime !} _{i} (\NN _{i})) } 
\ar@{_{(}->}[l]
&
{ \rho ^{!}( \NN _{i}) .}
\ar[l] ^-{\sim}
}
\end{equation}
Il en résulte le diagramme commutatif suivant : 
\begin{equation}
\label{Gamma-u+u!-tau}
\xymatrix @R=0,3cm { 
{\underline{\Gamma} ^{(m)} _{X} ( \NN _i)} 
\ar@{=}[d]
& 
{u _{0*} \rho ^{!} (\NN _{i}) } 
\ar@{^{(}->}[d]
\ar@{_{(}->}[l]
&
{u _{0*} \rho ^{!} (\NN _{i}) } 
\ar@{=}[l]
\ar@{^{(}->}[d]
\ar@{^{(}->}[r]
&
{\underline{\Gamma} ^{(m)} _{X} ( \NN _i)} 
\ar@{=}[d]
\\ 
{\underline{\Gamma} ^{(m)} _{X} ( \NN _i)} 
\ar@{=}[d]
& 
{u _{0*} u _{i} ^{!} (\NN _{i}) } 
\ar[r] ^-{\sim} _-{\tau}
\ar@{_{(}->}[l]
\ar@{^{(}->}[d]
&
{u _{0*}u _{i} ^{\prime !} (\NN _{i}) } 
\ar@{^{(}->}[r]
\ar@{^{(}->}[d]
&
{\underline{\Gamma} ^{(m)} _{X} ( \NN _i)} 
\ar@{=}[d]
\\
{\underline{\Gamma} ^{(m)} _{X} ( \NN _i)} 
& 
{u _{i +}  u _{i} ^{!} (\NN _{i}) } 
\ar[r] ^-{\sim} _-{\tau}
\ar[l] ^-{\sim}
&
{u ' _{i +}  u _{i} ^{\prime !} (\NN _{i}) } 
\ar[r] ^-{\sim}
& 
{\underline{\Gamma} ^{(m)} _{X} ( \NN _i) } 
.}
\end{equation}
\end{vide}

\begin{vide}
\label{Gammaiso-u0}
Soit $ \NN _i $ un $\D _{X _i} ^{(m)}$-module à droite quasi-cohérent injectif.
Comme le diagramme \ref{Gamma-u+u!-tau} est commutatif, l'isomorphisme de $\D _{X _i} ^{(m)}$-modules $u _{i +}  u _{i} ^{!} (\NN _{i}) \riso \underline{\Gamma} ^{(m)} _{X} ( \NN _i)$ valable lorsque $u _{0}$ se relève en une immersion fermée $u \,:\, \X \hookrightarrow \Y$ s'étend par recollement au cas où $u _{0}$ est quelconque. On dispose ainsi de l'isomorphisme canonique de $\D _{X _i} ^{(m)}$-modules 
$u _{0+}  u _{0} ^{!} (\NN _{i}) \riso \underline{\Gamma} ^{(m)} _{X} ( \NN _i)$. 

En prenant des résolutions injectives, on obtient alors, pour tout $\NN _i \in D ^+ _\mathrm{qc} (\D _{Y _i} ^{(m)} \overset{^d}{})$, l'isomorphisme canonique $u _{0+}  u _{0} ^{!} (\NN _{i}) \riso \R\underline{\Gamma} ^{(m)} _{X} ( \NN _i)$ s'inscrivant dans le diagramme commutatif : 
\begin{equation}
\label{u0+u0!=Gamma}
\xymatrix @R=0,3cm{
{u _{0+}  u _{0} ^{!} (\NN _{i}) } 
\ar[r] ^-{\sim}
\ar[d] ^-{\mathrm{can}}
& 
{\R\underline{\Gamma} ^{(m)} _{X} ( \NN _i) } 
\ar[d] ^-{\mathrm{can}}
\\ 
{\NN _{i}} 
\ar@{=}[r] 
& 
{\NN _{i} } 
}
\end{equation}

$\bullet$ Les isomorphismes de
\cite[1.1.9 et 1.1.10]{caro_courbe-nouveau} s'étendent au cas non relevable.
En particulier (en omettant les foncteurs oublis), les foncteurs
$u _0 ^! $ et $ u _{0T} ^!$ (resp. $u _{0+} $ et $u _{0T+}$)
sont canoniquement isomorphes (sur les catégories de la forme ${\underset{^{\longrightarrow }}{LD }}  ^\mathrm{b} _{\Q, \mathrm{qc}} (\smash{\widetilde{\D} } _{\X} ^{(\bullet)} \overset{^{\mathrm{d}}}{})$ et ${\underset{^{\longrightarrow }}{LD }}  ^\mathrm{b} _{\Q, \mathrm{qc}} (\smash{\widetilde{\D} } _{\Y} ^{(\bullet)} \overset{^{\mathrm{d}}}{})$).

$\bullet $ Soit $\NN \in D ^\mathrm{b} _\mathrm{coh} ( \D ^\dag _{\Y } (\hdag T) _{\Q} \overset{^{\mathrm{d}}}{})$
tel que $ u  _{0T }^! (\NN) \in
D ^\mathrm{b} _\mathrm{coh} ( \D ^\dag _{\X} (\hdag T _{X}) _{\Q} \overset{^{\mathrm{d}}}{})$.
On déduit de \ref{u0+u0!=Gamma} l'isomorphisme canonique : $u _{0T+} u _{0 T} ^! (\NN) \riso
\R \underline{\Gamma} _{X} ^\dag (\NN)$.
De manière analogue à \ref{+!gammmfrob}, on vérifie que celui-ci est compatible à Frobenius.
En effet, le théorème de Berthelot-Kashiwara reste valable dans le cas non relevable, e.g., si $\NN$ est à support dans $X$, le morphisme canonique $u _{0T+} u _{0 T} ^! (\NN)\to \NN$ est un isomorphisme (en effet, comme cela est local, on se ramène au cas où $u _{0}$ se relève et où le diviseur $T$ est vide). 

$\bullet $ Comme $u _{0T} ^! (\omega _{\Y} (\hdag T ) _{\Q}[d _Y]) \riso \omega _{\X} (\hdag T _X) _{\Q}[d _X]$, 
on construit alors un morphisme trace compatible à Frobenius en posant 
\begin{equation}
\label{can<->tr-gen}
\mathrm{Tr} _{T+} \ : \
u _{0T+} (\omega _{\X} (\hdag T _X) _{\Q}[d _X]) \riso \R \underline{\Gamma} _X ^\dag (\omega _{\Y} (\hdag T ) _{\Q} [d _Y]) \rightarrow \omega _{\Y} (\hdag T ) _{\Q}[d _Y].
\end{equation}
D'après \ref{tr+compfrob}, lorsque $u$ se relève on retrouve le morphisme trace de  \cite[1.2.5,6]{caro_courbe-nouveau} (grâce à \ref{tr+compfrob}).
\end{vide}

\subsection{Construction et compatibilité à Frobenius de l'isomorphisme de dualité relative dans le cas d'une immersion fermée non relevable}

On suppose que $u _{0}$ est une immersion fermée. 
De plus, nous supposerons le diviseur $T$ vide jusqu'au théorème \ref{isodualrelfrob} (nous nous ramènerons dans la preuve de \ref{isodualrelfrob} au cas où le diviseur est vide).

\begin{vide}
\label{formelvide1}
  Soient $ \M \in D ^\mathrm{b} _\mathrm{coh} ( \widehat{\D} ^{(m)} _{\X} \overset{^{\mathrm{d}}}{})$
et
$\NN \in D ^\mathrm{b} _{\mathrm{qc},\mathrm{qc}}
(\widehat{\D} ^{(m)} _{\X } \overset{^{\mathrm{d}}}{} ,\widehat{\D} ^{(m)} _{\X } \overset{^{\mathrm{d}}}{})$.
Par localisation et grâce au lemme sur les foncteurs way-out à droite (e.g., voir \cite{HaRD}), on établit que le morphisme canonique
$$\R \mathcal{H} om _{\widehat{\D} ^{(m)} _{\X}} ( \M,\NN) \rightarrow
\R \mathcal{H} om _{\widehat{\D} ^{(m)} _{\X}} ( \M,\NN)
\smash{\widehat{\otimes}} ^\L _{\widehat{\D} ^{(m)} _{\X}}
\widehat{\D} ^{(m)} _{\X}
:=
\R \underset{\underset{i}{\longleftarrow}}{\lim}
\R \mathcal{H} om _{\widehat{\D} ^{(m)} _{\X}} ( \M,\NN)
\otimes ^\L _{\widehat{\D} ^{(m)} _{\X}} \D ^{(m)} _{X _i}
$$
\normalsize
est un isomorphisme.
De même, en notant
$ \M _i :=  \M \otimes ^\L _{\widehat{\D} ^{(m)} _{\X}} \D ^{(m)} _{X _i}$
et
$\NN _i := \NN \otimes ^\L _{\widehat{\D} ^{(m)} _{\X}} \D ^{(m)} _{X _i}$,
le morphisme canonique
$\R \mathcal{H} om _{\widehat{\D} ^{(m)} _{\X}} ( \M,\NN)
\otimes ^\L _{\widehat{\D} ^{(m)} _{\X}} \D ^{(m)} _{X _i}
\rightarrow
\R \mathcal{H} om _{\D ^{(m)} _{X _i}} ( \M _i ,\NN _i)$
est un isomorphisme. D'où :
$$\R \mathcal{H} om _{\widehat{\D} ^{(m)} _{\X}} ( \M,\NN)
\riso
\R \underset{\underset{i}{\longleftarrow}}{\lim} \,
\R \mathcal{H} om _{\D ^{(m)} _{X _i}} ( \M _i ,\NN _i).$$

\end{vide}

\begin{vide}
\label{formelvidepre2}
Soit $\lambda \,:\, \N \to \N$ une application croissante telle que $\lambda (m)\geq m$ pour tout $m$.
    Soit
$ \M ^{ (\bullet)} \in \smash[b]{\underset{^{\longrightarrow }}{LD }}  ^\mathrm{b} _{\Q}
(\lambda ^{*}\smash{\widehat{\D}} _{\X } ^{(\bullet)} \overset{^{\mathrm{d}}}{})$
tel que $ \M ^{ (m)} \in D ^\mathrm{b} _{\mathrm{coh}}
(\smash{\widehat{\D}} _{\X } ^{(\lambda (m) ) } \overset{^{\mathrm{d}}}{})$ et tel que, pour tous entiers $m '\geq m$,
$ \M ^{ (m)} \otimes ^\L _{\widehat{\D} ^{(\lambda (m) ) } _{\X }} \widehat{\D} ^{(\lambda (m') )} _{\X} \rightarrow  \M ^{ (m')}$
soit un isomorphisme.
De plus, donnons-nous
$ \NN ^{ (\bullet)} \in \smash[b]{\underset{^{\longrightarrow }}{LD }}  ^\mathrm{b} _{\Q}
(\lambda ^{*}\smash{\widehat{\D}} _{\X } ^{(\bullet)} \overset{^{\mathrm{d}}}{},
\lambda ^{*}\smash{\widehat{\D}} _{\X } ^{(\bullet)} \overset{^{\mathrm{d}}}{})$ tel que
$\NN ^{ (m)} \in D  ^\mathrm{b} _{\mathrm{qc},\mathrm{qc}}
(\smash{\widehat{\D}} _{\X } ^{(\lambda (m) ) } \overset{^{\mathrm{d}}}{},
\smash{\widehat{\D}} _{\X } ^{(\lambda (m) ) } \overset{^{\mathrm{d}}}{})$.
Comme $D ^\mathrm{b} _{\mathrm{coh}}
(\smash{\widehat{\D}} _{\X } ^{(\lambda (m))} \overset{^{\mathrm{d}}}{})=D ^\mathrm{b} _{\mathrm{parf}}
(\smash{\widehat{\D}} _{\X } ^{(\lambda (m))} \overset{^{\mathrm{d}}}{})$ (remarquons que l'on a besoin a priori d'éviter d'ajouter des singularités surconvergentes, i.e., $T$ doit être vide), 
par localisation et dévissage, on vérifie alors que le complexe
\begin{equation}
  \R \mathcal{H} om _{\smash[b]{\underset{^{\longrightarrow }}{LD }}  ^\mathrm{b} _{\Q, \mathrm{qc}}
(\lambda ^{*}\smash{\widehat{\D}} _{\X} ^{(\bullet)} \overset{^{\mathrm{d}}}{})}
( \M ^{ (\bullet)} , \NN ^{ (\bullet)} ) :=
( \R \mathcal{H} om _{\widehat{\D} ^{(\lambda (m))} _{\X}} ( \M ^{(m)} ,\NN ^{(m)})) _{m \in \N}
\end{equation}
est un objet de
$\smash[b]{\underset{^{\longrightarrow }}{LD }}  ^\mathrm{b} _{\Q, \mathrm{qc}}
(\lambda ^{*}\smash{\widehat{\D}} _{\X } ^{(\bullet)} \overset{^{\mathrm{d}}}{})$.

En effet,
comme
$ \M ^{ (m)} \otimes ^\L _{\widehat{\D} ^{(\lambda (m) ) } _{\X }} \widehat{\D} ^{(\lambda (m') )} _{\X} \rightarrow  \M ^{ (m')}$ est un isomorphisme,
les morphismes de transition
sont construits pour tous $m' \geq m$ via les morphismes :
\begin{equation}
   \R \mathcal{H} om _{\widehat{\D} ^{(\lambda (m))} _{\X}} ( \M ^{(m)} ,\NN ^{(m)})
\rightarrow
   \R \mathcal{H} om _{\widehat{\D} ^{(\lambda (m))} _{\X}} ( \M ^{(m)} ,\NN ^{(m')})
\liso
\R \mathcal{H} om _{\widehat{\D} ^{(\lambda (m'))} _{\X}} ( \M ^{(m ')} ,\NN ^{(m')}).
\end{equation}
La quasi-cohérence résulte alors de \ref{formelvide1}.

\end{vide}

\begin{vide}
\label{formelvide2}
Fixons $\lambda \,:\, \N \to \N$ une application croissante telle que $\lambda (m)\geq m$ pour tout $m$.
Notons $\mu \,:\, \N \to \N$ l'application définie par $\mu (m) = \lambda (m)+s$ pour tout $m$. 
 Soit
$ \M ^{ (\bullet)} \in \smash[b]{\underset{^{\longrightarrow }}{LD }}  ^\mathrm{b} _{\Q}
(\lambda ^{*}\smash{\widehat{\D}} _{\X } ^{(\bullet)} \overset{^{\mathrm{d}}}{})$
tel que $ \M ^{ (m)} \in D ^\mathrm{b} _{\mathrm{coh}}
(\smash{\widehat{\D}} _{\X } ^{(\lambda (m) ) } \overset{^{\mathrm{d}}}{})$ et tel que, pour tous entiers $m '\geq m$,
$ \M ^{ (m)} \otimes ^\L _{\widehat{\D} ^{(\lambda (m) ) } _{\X }} \widehat{\D} ^{\lambda (m') } _{\X} \rightarrow  \M ^{ (m')}$
soit un isomorphisme.
D'après le cas algébrique (\ref{theo1-7}), pour tous entiers $i$ et $m$,
on dispose du diagramme commutatif
\begin{equation}
\label{isodualrelindpre}
\xymatrix @R=0,3cm @C=0,5cm {
{  u _{0+}^{(\mu (m))} \DD _{X _i} ^{(\mu (m))}( F ^\flat \M _i ^{\prime (m)} )}
\ar[r] \ar[d] _-\sim
&
{\R \mathcal{H} om _{\D ^{(\mu (m))} _{Y _i}}
(u _{0+} ^{(\mu (m))} ( F ^\flat \M _i ^{\prime (m)}),
u _{0+}^{(\mu (m))} (\omega _{X _i} )\otimes ^\L _{\O _{Y _i}} \D _{Y _i} ^{(\mu (m))}) [d _X]}
\ar[d] _-\sim
\\
{ F ^\flat  u _{0+}^{\prime (\lambda (m))} \DD _{X '_i} ^{(\lambda (m))}( \M _i ^{\prime (m)} )}
\ar[r]
&
{F ^\flat  \R \mathcal{H} om _{\D ^{(\lambda (m))} _{Y '_i}}
(u _{0+} ^{\prime (\lambda (m))} ( \M _i ^{\prime (\lambda (m))}),
u _{0+}^{\prime (\lambda (m))} (\omega _{X '_i} )\otimes ^\L _{\O _{Y '_i}} \D _{Y '_i} ^{(\lambda (m))}) [d _X],}
}
\end{equation}
\normalsize
compatible au changement de niveaux (voir le cube commutatif \ref{cubecgtnivfrob}).

Notons ici $u _{0+} ^{\prime (\bullet)}= ( u _{0+} ^{\prime (\lambda (m))}) _{m\in \N}\, :\,
\smash[b]{\underset{^{\longrightarrow }}{LD }}  ^\mathrm{b} _{\Q}
(\lambda ^{*}\smash{\widehat{\D}} _{\X' } ^{(\bullet)} \overset{^{\mathrm{d}}}{})
\to \smash[b]{\underset{^{\longrightarrow }}{LD }}  ^\mathrm{b} _{\Q}
(\lambda ^{*}\smash{\widehat{\D}} _{\Y'} ^{(\bullet)} \overset{^{\mathrm{d}}}{})$,
$\DD _{\X '} ^{(\bullet )}= (\DD _{\X '} ^{( \lambda (m) )}) _{m\in \N } \, :\,
\smash[b]{\underset{^{\longrightarrow }}{LD }}  ^\mathrm{b} _{\Q}
(\lambda ^{*}\smash{\widehat{\D}} _{\X '} ^{(\bullet)} \overset{^{\mathrm{d}}}{})
\to \smash[b]{\underset{^{\longrightarrow }}{LD }}  ^\mathrm{b} _{\Q}
(\lambda ^{*}\smash{\widehat{\D}} _{\X' } ^{(\bullet)} \overset{^{\mathrm{d}}}{})$,
$u _{0+} ^{(\bullet)}= ( u _{0+} ^{(\mu (m))}) _{m\in \N}\, :\,
\smash[b]{\underset{^{\longrightarrow }}{LD }}  ^\mathrm{b} _{\Q}
(\mu ^{*}\smash{\widehat{\D}} _{\X } ^{(\bullet)} \overset{^{\mathrm{d}}}{})
\to \smash[b]{\underset{^{\longrightarrow }}{LD }}  ^\mathrm{b} _{\Q}
(\mu ^{*}\smash{\widehat{\D}} _{\Y} ^{(\bullet)} \overset{^{\mathrm{d}}}{})$,
$\DD _{\X } ^{(\bullet ) }= (\DD _{\X } ^{( \mu (m) )}) _{m\in \N } \, :\,
\smash[b]{\underset{^{\longrightarrow }}{LD }}  ^\mathrm{b} _{\Q}
(\mu ^{*}\smash{\widehat{\D}} _{\X } ^{(\bullet)} \overset{^{\mathrm{d}}}{})
\to \smash[b]{\underset{^{\longrightarrow }}{LD }}  ^\mathrm{b} _{\Q}
(\mu ^{*}\smash{\widehat{\D}} _{\X } ^{(\bullet)} \overset{^{\mathrm{d}}}{})$ etc. 
De plus, notons $\omega^{(\bullet)} _{\X'}$ le système inductif constant de $\omega _{\X'}$ dans 
$\smash[b]{\underset{^{\longrightarrow }}{LD }}  ^\mathrm{b} _{\Q}
(\lambda ^{*}\smash{\widehat{\D}} _{\X } ^{(\bullet)} \overset{^{\mathrm{d}}}{})$,
$\omega^{(\bullet)} _{\X}$ le système inductif constant de $\omega _{\X}$ dans 
$\smash[b]{\underset{^{\longrightarrow }}{LD }}  ^\mathrm{b} _{\Q}
(\mu ^{*}\smash{\widehat{\D}} _{\X } ^{(\bullet)} \overset{^{\mathrm{d}}}{})$ etc.

En appliquant le foncteur $\R \underset{\underset{i}{\longleftarrow}}{\lim}$ au diagramme \ref{isodualrelindpre},
via \ref{formelvide1} et \ref{formelvidepre2},
on obtient la commutativité du suivant :
\begin{equation}
\label{isodualrelind}
\xymatrix @R=0,3cm @C=0,5cm {
{ u _{0+} ^{(\bullet )} \DD _{\X} ^{(\bullet )} (F ^\flat \M ^{\prime (\bullet)}) }
\ar[r] \ar[d] _-\sim
&
{\R \mathcal{H} om _{\smash[b]{\underset{^{\longrightarrow }}{LD }}  ^\mathrm{b} _{\Q, \mathrm{qc}}
(\mu ^{*}\smash{\widehat{\D}} _{\Y} ^{(\bullet)} \overset{^{\mathrm{d}}}{})}
(  u _{0+} ^{(\bullet)} ( F ^\flat \M ^{\prime (\bullet)}) ,
  u _{0+} ^{(\bullet)} (\omega^{(\bullet)} _{\X})
\smash{\widehat{\otimes}} ^\L  _{\O _{\Y}} \mu ^{*}\widehat{\D} ^{(\bullet)} _{\Y}[d _X])}
\ar[d] _-\sim
\\
{ F ^\flat  u _{0+} ^{\prime (\bullet)} \DD _{\X '} ^{(\bullet )} (\M ^{\prime (\bullet)}) }
\ar[r]
&
{F ^\flat \R \mathcal{H} om _{\smash[b]{\underset{^{\longrightarrow }}{LD }}  ^\mathrm{b} _{\Q, \mathrm{qc}}
(\lambda ^{*}\smash{\widehat{\D}} _{\Y'} ^{(\bullet)} \overset{^{\mathrm{d}}}{})}
(  u _{0+} ^{\prime (\bullet)} ( \M ^{\prime (\bullet)}) ,
  u _{0+} ^{\prime (\bullet)} (\omega ^{(\bullet)} _{\X'})
\smash{\widehat{\otimes}} ^\L  _{\O _{\Y'}} \lambda ^{*} \widehat{\D} ^{(\bullet)} _{\Y'}[d _X]).}
}
\end{equation}

\normalsize
D'après \ref{can<->tr-gen}
et grâce à l'équivalence de catégories
$\underset{\longrightarrow}{\lim}$ :
$\smash[b]{\underset{^{\longrightarrow }}{LD }}  ^\mathrm{b} _{\Q, \mathrm{coh}}
(\smash{\widehat{\D}} _{\X} ^{(\bullet)} \overset{^{\mathrm{d}}}{})
\cong
D ^{\mathrm{b}} _\mathrm{coh} ( \D ^\dag _{\X,\Q})$ (voir \cite[4.2.4]{Beintro2}),
le morphisme trace
$\mathrm{Tr} _+$ :
$u _{0+} ^{(\bullet)} (\omega _{\X}^{(\bullet)}) [d _X] \rightarrow \omega _{\Y}^{(\bullet)} [d _Y]$,
est compatible à Frobenius, i.e.,
le diagramme ci-dessous
\begin{equation}
  \label{7-finformelpre1}
\xymatrix @R=0,3cm {
{u _{0+} ^{(\bullet)} ( \omega _{\X }^{(\bullet)} [d _{X}]) }
\ar[r] ^-{\mathrm{Tr} _+}
\ar[d] _-\sim
&
{\omega _{\Y} ^{(\bullet)} [d _{Y}]}
\ar[d] _-\sim
\\
{F ^\flat u _{0+} '(\omega _{\X '} ^{(\bullet)} [d _{X'}]) }
\ar[r] ^-{\mathrm{Tr} _+}
&
{F ^\flat\omega _{\Y '} ^{(\bullet)}[d _{Y'}]}
}
\end{equation}
est commutatif.
En appliquant $-\smash{\widehat{\otimes}} ^\L  _{\O _{\Y}} \mu ^{*}\widehat{\D} ^{(\bullet)} _{\Y}$
à \ref{7-finformelpre1},
on obtient le carré du haut de :
\begin{equation}
  \label{7-finformelpre2}
\xymatrix @R=0,3cm {
{u _{0+} ^{(\bullet)} ( \omega _{\X} ^{(\bullet)})
\smash{\widehat{\otimes}} ^\L  _{\O _{\Y}} \mu ^{*} \widehat{\D} ^{(\bullet)} _{\Y} [d _{X}]}
\ar[r] ^-{\mathrm{Tr} _+}
\ar[d] _-\sim
&
{\omega _{\Y} ^{(\bullet)} \smash{\widehat{\otimes}} ^\L  _{\O _{\Y}} \mu ^{*} \widehat{\D} ^{(\bullet)} _{\Y}[d _{Y}]}
\ar[d] _-\sim
\\
{F ^\flat u _{0+} ^{\prime (\bullet)} ( \omega _{\X '} ^{(\bullet)})  \smash{\widehat{\otimes}} ^\L  _{\O _{\Y}} \mu ^{*} \widehat{\D} ^{(\bullet)} _{\Y}[d _{X}]}
\ar[r] ^-{\mathrm{Tr} _+}
\ar[d] _-\sim
&
{F ^\flat \omega _{\Y'} ^{(\bullet)} \smash{\widehat{\otimes}} ^\L  _{\O _{\Y}} \mu ^{*} \widehat{\D} ^{(\bullet)} _{\Y}[d _{Y'}]}
\ar[d] _-\sim
\\
{F ^\flat u _{0+} ^{\prime (\bullet)} ( \omega _{\X '} ^{(\bullet)})
\smash{\widehat{\otimes}} ^\L  _{\O _{\Y}} F ^* F^\flat  \lambda ^{*} \widehat{\D} ^{(\bullet)} _{\Y'}[d _{X}]}
\ar[r] ^-{\mathrm{Tr} _+}
\ar[d] _-\sim
&
{F ^\flat \omega _{\Y'} ^{(\bullet)} \smash{\widehat{\otimes}} ^\L  _{\O _{\Y}} F ^* F^\flat  \lambda ^{*} \widehat{\D} ^{(\bullet)} _{\Y'}[d _{Y'}]}
\ar[d] _-\sim
\\
{F ^\flat _\mathrm{g} F ^\flat _\mathrm{d} (
u _{0+} ^{\prime (\bullet)}
( \omega _{\X '} ^{(\bullet)}) \smash{\widehat{\otimes}} ^\L  _{\O _{\Y'}} \lambda ^{*} \widehat{\D} ^{(\bullet)} _{\Y'})[d _{X}]}
\ar[r] ^-{\mathrm{Tr} _+}
&
{F ^\flat _\mathrm{g} F ^\flat _\mathrm{d}
(\omega _{\Y'} ^{(\bullet)}\smash{\widehat{\otimes}} ^\L  _{\O _{\Y'}} \lambda ^{*} \widehat{\D} ^{(\bullet)} _{\Y'})[d _{Y'}].}
}
\end{equation}
On remarque que l'on peut remplacer $\smash{\widehat{\otimes}} ^\L  $ par $\otimes $ dans la partie droite de \ref{7-finformelpre2}. 
On vérifie par fonctorialité la commutativité des deux autres carrés de \ref{7-finformelpre2}.
D'où celle de \ref{7-finformelpre2}.
En lui appliquant
$\R \mathcal{H} om _{\smash[b]{\underset{^{\longrightarrow }}{LD }}  ^\mathrm{b} _{\Q, \mathrm{qc}}
(\mu ^{*}\smash{\widehat{\D}} _{\Y} ^{(\bullet)} \overset{^{\mathrm{d}}}{})}
(u _{0+} ^{(\bullet)} ( F ^\flat \M ^{\prime (\bullet)}) ,-)$, noté ici plus simplement
$\R \mathcal{H} om
(u _{0+} ^{(\bullet)} ( F ^\flat \M ^{\prime (\bullet)}) ,-)$,
on arrive au carré supérieur ci-dessous :
\begin{equation}
    \label{7-finformelpre3}
\xymatrix @R=0,3cm {
{\R \mathcal{H} om
(  u _{0+} ^{(\bullet)}  F ^\flat \M ^{\prime (\bullet)} ,
u _{0+} ^{ (\bullet)} ( \omega _{\X } ^{(\bullet)})
\smash{\widehat{\otimes}} ^\L  _{\O _{\Y}} \mu ^{*} \widehat{\D} ^{(\bullet)} _{\Y})
[d _{X}]}
\ar[r] ^-{\mathrm{Tr} _+}
\ar[d] _-\sim
&
{\R \mathcal{H} om
(u _{0+} ^{(\bullet)}  F ^\flat \M ^{\prime (\bullet)},
\omega _{\Y } ^{(\bullet)} \otimes  _{\O _{\Y}} \mu ^{*} \widehat{\D} ^{(\bullet)} _{\Y}[d _{Y}])}
\ar[d] _-\sim
\\
{\R \mathcal{H} om
(u _{0+} ^{(\bullet)}  F ^\flat \M ^{\prime (\bullet)},
F ^\flat _\mathrm{g} F ^\flat _\mathrm{d} ( u _{0+} ^{\prime (\bullet)} ( \omega _{\X '} ^{(\bullet)})
\smash{\widehat{\otimes}} ^\L  _{\O _{\Y'}} \lambda ^{*} \widehat{\D} ^{(\bullet)} _{\Y'})) [d _{X}]}
\ar[r] ^-{\mathrm{Tr} _+}
&
{\R \mathcal{H} om ( u _{0+} ^{(\bullet)}  F ^\flat \M ^{\prime (\bullet)},
F ^\flat _\mathrm{g} F ^\flat _\mathrm{d}
(\omega _{\Y'} ^{(\bullet)}\otimes  _{\O _{\Y'}} \lambda ^{*} \widehat{\D} ^{(\bullet)} _{\Y'})
)[d _{Y}]}
\\
{\R \mathcal{H} om
(F ^\flat  u _{0+} ^{\prime (\bullet)} \M ^{\prime (\bullet)},
F ^\flat _\mathrm{g} F ^\flat _\mathrm{d} ( u _{0+} ^{\prime (\bullet)} ( \omega _{\X '} ^{(\bullet)})
\smash{\widehat{\otimes}} ^\L  _{\O _{\Y'}} \lambda ^{*} \widehat{\D} ^{(\bullet)} _{\Y'}) )[d _{X}]}
\ar[r] ^-{\mathrm{Tr} _+} \ar[u] _-\sim
&
{\R \mathcal{H} om
(F ^\flat  u _{0+} ^{\prime (\bullet)} \M ^{\prime (\bullet)},
F ^\flat _\mathrm{g} F ^\flat _\mathrm{d}
(\omega _{\Y'} ^{(\bullet)}\otimes  _{\O _{\Y'}} \lambda ^{*} \widehat{\D} ^{(\bullet)} _{\Y'}))[d _{Y}]}
\ar[u] _-\sim
\\
{F ^\flat \R \mathcal{H} om
(u _{0+} ^{\prime (\bullet)} ( \M ^{\prime (\bullet)}),
u _{0+} ^{\prime (\bullet)} ( \omega _{\X '} ^{(\bullet)})
\smash{\widehat{\otimes}} ^\L  _{\O _{\Y'}} \lambda ^{*} \widehat{\D} ^{(\bullet)} _{\Y'} )[d _{X}]}
\ar[r] ^-{\mathrm{Tr} _+} \ar[u] _-\sim ^{F ^\flat}
&
{F ^\flat \R \mathcal{H} om
(u _{0+} ^{\prime (\bullet)} ( \M ^{\prime (\bullet)}),
\omega _{\Y'} ^{(\bullet)}\otimes  _{\O _{\Y'}} \lambda ^{*} \widehat{\D} ^{(\bullet)} _{\Y'})[d _{Y}].}
\ar[u] _-\sim ^{F ^\flat}
}
\end{equation}
\normalsize
Puisque les deux carrés du bas sont commutatifs par fonctorialité,
\ref{7-finformelpre3} est commutatif.
De plus, comme pour \ref{1-2bis}, on dispose du diagramme commutatif :
\begin{equation}
    \label{7-finformelpre4}
\xymatrix @R=0,3cm @C=0,5cm {
{\R \mathcal{H} om _{\smash[b]{\underset{^{\longrightarrow }}{LD }}  ^\mathrm{b} _{\Q, \mathrm{qc}}
(\mu ^{*}\smash{\widehat{\D}} _{\Y} ^{(\bullet)} \overset{^{\mathrm{d}}}{})}
(u _{0+} ^{(\bullet)}  F ^\flat \M ^{\prime (\bullet)},
\omega _{\Y } ^{(\bullet)} \otimes  _{\O _{\Y}} \mu ^{*} \widehat{\D} ^{(\bullet)} _{\Y}[d _{Y}])}
\ar[d] _-\sim
\ar[r] _-\sim
&
{\DD _{\Y } ^{(\bullet)} \circ u _{0+} ^{(\bullet)} (F ^\flat  \M  ^{\prime (\bullet )})}
\ar[d] _-\sim
\\
{F ^\flat \R \mathcal{H} om _{\smash[b]{\underset{^{\longrightarrow }}{LD }}  ^\mathrm{b} _{\Q, \mathrm{qc}}
(\lambda ^{*}\smash{\widehat{\D}} _{\Y'} ^{(\bullet)} \overset{^{\mathrm{d}}}{})}
(u _{0+} ^{\prime (\bullet)} ( \M ^{\prime (\bullet)}),
\omega _{\Y'} ^{(\bullet)}\otimes  _{\O _{\Y'}} \lambda ^{*} \widehat{\D} ^{(\bullet)} _{\Y'} )[d _{Y}]}
\ar[r]_-\sim
&
{F ^\flat \DD  _{\Y '} ^{(\bullet )}  \circ u _{0+} ^{\prime (\bullet )} (\M  ^{\prime (\bullet )} ).}
}
\end{equation}
En composant \ref{isodualrelind}, \ref{7-finformelpre3} et \ref{7-finformelpre4},
on obtient le diagramme canonique commutatif :
\begin{equation}
    \label{7-finformelpre5}
\xymatrix @R=0,3cm {
{u _{0+} ^{(\bullet)} \circ \DD _{\X } ^{(\bullet)} ( F ^\flat  \M  ^{\prime (\bullet )}) }
\ar[d] _-\sim
\ar[r] _-\sim
&
{\DD _{\Y } ^{(\bullet)} \circ u _{0+} ^{(\bullet)} (F ^\flat  \M  ^{\prime (\bullet )})}
\ar[d] _-\sim
\\
{F ^\flat  u _{0+}  ^{\prime (\bullet )}\circ \DD _{\X '} ^{(\bullet )} (  \M  ^{\prime (\bullet )})}
\ar[r]_-\sim
&
{F ^\flat \DD  _{\Y '} ^{(\bullet )}  \circ u _{0+} ^{\prime (\bullet )} (\M  ^{\prime (\bullet )} ).}
}
\end{equation}
\end{vide}

\begin{theo}\label{isodualrelfrob}
Pour tout $ \M\in D ^{\mathrm{b}} _\mathrm{coh}(\D ^\dag _{\X } (\hdag T _X) _{\Q} \overset{^{\mathrm{d}}}{})$,
on dispose de l'isomorphisme canonique fonctoriel en $\M$
\begin{equation}
\label{isodualrelfrob-eq}
\chi\ : \ u _{0T +} \circ \DD _{\X, T _X} ( \M) \riso \DD _{\Y, T }  \circ u _{0 T +}  ( \M),
\end{equation}
dit {\og isomorphisme de dualité relative \fg}.
Cet isomorphisme est compatible à Frobenius, i.e.,
pour tout
$ \M' \in D ^{\mathrm{b}} _\mathrm{coh}(\D ^\dag _{\X '} (\hdag T _{X'}) _{\Q} \overset{^{\mathrm{d}}}{})$,
le diagramme
\begin{equation}
  \label{diagdualreldiv}
\xymatrix @R=0,3cm {
{ u _{0T +} \circ \DD _{\X, T _X} (F ^\flat  \M ')}
\ar[r]_-\sim ^-{\chi} \ar[d] _-\sim
&
{\DD _{\Y, T }  \circ u _{0T +}  (F ^\flat  \M ')}
\ar[d] _-\sim
\\
{F ^\flat u _{0 T '+} '\circ \DD _{\X',T _{X'}} ( \M ')}
\ar[r] _-\sim ^-{\chi}
&
{F ^\flat \DD _{\Y',T '}  \circ u _{0T' +} ' ( \M '),}
}
\end{equation}
où les isomorphismes
verticaux dérivent des isomorphismes de commutation à Frobenius de l'image directe
et du foncteur dual, est commutatif.
De plus, l'isomorphisme \ref{isodualrelfrob-eq} vérifie la condition usuelle de transitivité pour le composé de deux immersions fermées. 
\end{theo}

\begin{proof}
Via les équivalences de catégories
$\underset{\longrightarrow}{\lim}$ :
$\smash[b]{\underset{^{\longrightarrow }}{LD }}  ^\mathrm{b} _{\Q, \mathrm{coh}}
(\smash{\widetilde{\D} } _{\X} ^{(\bullet)} \overset{^{\mathrm{d}}}{})
\cong
D ^{\mathrm{b}} _\mathrm{coh} ( \D ^\dag _{\X} (\hdag T _{X} ) _{\Q})$,
$\underset{\longrightarrow}{\lim}$ :
$\smash[b]{\underset{^{\longrightarrow }}{LD }}  ^\mathrm{b} _{\Q, \mathrm{coh}}
(\smash{\widetilde{\D} } _{\Y} ^{(\bullet)} \overset{^{\mathrm{d}}}{})
\cong
D ^{\mathrm{b}} _\mathrm{coh} ( \D ^\dag _{\Y} (\hdag T ) _{\Q})$,
il découle de \ref{can<->tr-gen} que l'on dispose du morphisme trace : 
$\mathrm{Tr} _{T+}$ :
$u _{0T+} ^{(\bullet)} (\widetilde{\omega} _{\X}^{(\bullet)}) [d _X] \rightarrow \widetilde{\omega} _{\Y}^{(\bullet)} [d _Y]$ (avec les notations de \ref{formelvide2}, ici $\lambda (m) =m$).
De manière analogue à \cite[1.2.7]{caro_courbe-nouveau}, on en déduit la construction de l'isomorphisme de dualité relative \ref{isodualrelfrob-eq}. 
Le fait que $\chi$ soit un isomorphisme est local. 
On se ramène au cas où $u _{0}$ se relève et où $T$ est vide (par \cite[4.3.12]{Be2}). 
On retrouve alors la construction de Virrion. Le morphisme $\chi$ est donc un isomorphisme.  

Traitons à présent sa compatibilité à Frobenius. Pour vérifier la commutativité du diagramme \ref{diagdualreldiv}, en appliquant le foncteur $\DD _{\Y, T }$ à \ref{diagdualreldiv}, modulo la commutation à Frobenius du foncteur dual et de l'isomorphisme de bidualité (voir \cite{virrion}), il s'agit d'établir la commutativité du diagramme 
\begin{equation}
  \label{diagdualreldiv-preuve}
\xymatrix @R=0,3cm {
{ u _{0T !}  (F ^\flat  \M ')}
\ar[r]_-\sim ^-{\chi} \ar[d] _-\sim
&
{ u _{0T +}  (F ^\flat  \M ')}
\ar[d] _-\sim
\\
{F ^\flat u _{0 T '!} ' ( \M ')}
\ar[r] _-\sim ^-{\chi}
&
{F ^\flat u _{0T' +} ' ( \M '),}
}
\end{equation}
où le foncteur image directe extraordinaire est défini en \ref{imag-direc-extra}. Comme les foncteurs $F ^\flat $, $u _{0T +}$ et $u _{0T !}$ sont exacts, l'assertion devient locale. On peut supposer que les diviseurs sont vides (par \cite[4.3.12]{Be2}). 
Via l'équivalence de catégories induite par le foncteur $\underset{\longrightarrow}{\lim}$ ci-dessus (avec $T $ vide), il existe $\lambda \,:\, \N \to \N$ une application croissante telle que $\lambda (m)\geq m$ pour tout $m$, 
il existe $ \M ^{ (\bullet)} \in \smash[b]{\underset{^{\longrightarrow }}{LD }}  ^\mathrm{b} _{\Q}
(\lambda ^{*}\smash{\widehat{\D}} _{\X } ^{(\bullet)} \overset{^{\mathrm{d}}}{})$
tel que $ \M ^{ (m)} \in D ^\mathrm{b} _{\mathrm{coh}}
(\smash{\widehat{\D}} _{\X } ^{(\lambda (m) ) } \overset{^{\mathrm{d}}}{})$ et tel que, pour tous entiers $m '\geq m$,
$ \M ^{ (m)} \otimes ^\L _{\widehat{\D} ^{(\lambda (m) ) } _{\X }} \widehat{\D} ^{\lambda (m') } _{\X} \rightarrow  \M ^{ (m')}$
soit un isomorphisme et vérifiant 
$\underset{\longrightarrow}{\lim}  (\M ^{ (\bullet)}) \riso \M '$.
On conclut alors en appliquant le foncteur $\underset{\longrightarrow}{\lim} $ à \ref{7-finformelpre5}.

Concernant la transitivité de l'isomorphisme de dualité relative, on procède de la même manière pour se ramener au cas relevable sans diviseur, cas déjà traité dans \cite{Beintro2}.

 \end{proof}

\begin{vide}
[Compléments]

Soit $\widetilde{\Y}$ un ouvert de $\Y$ de fibre spéciale $\widetilde{Y}$
tel que, en notant $\widetilde{\X}$ l'ouvert de $\X$ d'espace topologique sous-jacent $u _0 ^{-1} (\widetilde{Y})$
et $\widetilde{X}$ sa fibre spéciale,
le morphisme
$\widetilde{X} \rightarrow \widetilde{Y}$ induit par $u _0$ se relève en un morphisme
$\widetilde{u}$ : $\widetilde{\X} \rightarrow \widetilde{\Y}$.
On note $\widetilde{T}: = T \cap \widetilde{Y}$,
$T _{\widetilde{X}} := T _X \cap \widetilde{X}$.

Pour tout $\M \in D ^{\mathrm{b}} _\mathrm{coh}(\D ^\dag _{\X } (\hdag T _{X}) _{\Q} \overset{^{\mathrm{d}}}{})$,
le diagramme
\begin{equation}
  \label{chirest=chiu1}
\xymatrix @R=0,3cm {
{ u _{0T +} \circ \DD _{\X, T _X} (\M ) | {\widetilde{\Y}}}
\ar[r] ^-{\chi} _-\sim
\ar[d] _-\sim
&
{\DD _{\Y, T }  \circ u _{0T +}  (\M )| {\widetilde{\Y}}}
\\
{\widetilde{u} _{\widetilde{T}+}\circ \DD _{\widetilde{\X}, T _{\widetilde{X}}} (\M | {\widetilde{\X}} ) }
\ar[r] ^-{\chi} _-\sim
&
{\DD _{\widetilde{\X}, T _{\widetilde{X}}} \circ \widetilde{u} _{\widetilde{T}+} (\M | {\widetilde{\X}} ) ,}
\ar[u] _-\sim
}
\end{equation}
où les isomorphismes horizontaux sont induits par \ref{u0+=u+dag},
celui du haut (resp. du bas) est l'isomorphisme de dualité relative de \ref{isodualrelfrob}
(resp. de \cite[1.2.7]{caro_courbe-nouveau}), est commutatif.

De plus, si $\widetilde{v}$ : $\widetilde{\X} \rightarrow \widetilde{\Y}$ est un second relèvement de
$\widetilde{X} \rightarrow \widetilde{Y}$, on obtient le diagramme commutatif ci-après
\begin{equation}
  \label{chirest=chiu2}
\xymatrix @R=0,3cm {
{\widetilde{v} _{\widetilde{T}+}\circ \DD _{\widetilde{\X}, T _{\widetilde{X}}} (\M | {\widetilde{\X}} ) }
\ar[r] ^-{\chi} _-\sim
&
{\DD _{\widetilde{\X}, T _{\widetilde{X}}}\circ \widetilde{v} _{\widetilde{T}+} (\M | {\widetilde{\X}} ) }
\ar[d] _-\sim ^-{\tau _{\widetilde{v},\widetilde{u}}}
\\
{\widetilde{u} _{\widetilde{T}+}\circ \DD _{\widetilde{\X}, T _{\widetilde{X}}} (\M | {\widetilde{\X}} ) }
\ar[r] ^-{\chi} _-\sim
\ar[u] _-\sim ^-{\tau _{\widetilde{v},\widetilde{u}}}
&
{\DD _{\widetilde{\X}, T _{\widetilde{X}}} \circ \widetilde{u} _{\widetilde{T}+} (\M | {\widetilde{\X}} ) ,}
}
\end{equation}
où les isomorphismes horizontaux sont les isomorphismes de recollement
de \ref{u0+=u+dagrecol}.
\end{vide}

\subsection{Compatibilité à Frobenius du morphisme d'adjonction dans le cas d'une immersion fermée}
On suppose que $u _{0}$ est une immersion fermée. 
\begin{theo}\label{frobadj}
  Soient $\M '\in D ^{\mathrm{b}} _\mathrm{coh}(\D ^\dag _{\X '} (\hdag T _{X'}) _{\Q} \overset{^{\mathrm{d}}}{})$
et
$\NN ' \in D ^{\mathrm{b}} _\mathrm{coh}(\D ^\dag _{\Y '} (\hdag T ') _{\Q} \overset{^{\mathrm{d}}}{})$.
On dispose du diagramme commutatif :
\begin{equation}\label{frobadjdiag}
  \xymatrix @R=0,3cm {
{ \R \mathcal{H} om _{\D ^\dag _{\Y '} (\hdag T ') _{\Q}} ( u _{0T'+} ' (\M '), \NN ')}
\ar[r] _-\sim \ar[d] _-\sim ^-{F ^{\flat}}
&
{u _{0*}  \R \mathcal{H} om _{\D ^\dag _{\X '} (\hdag T _{X'}) _{\Q}} ( \M ', u _{0T'} ^{\prime !} \NN ')}
\ar[d] _-\sim ^-{F ^{\flat}}
\\
{ \R \mathcal{H} om _{\D ^\dag _{\Y } (\hdag T ) _{\Q}} ( F ^\flat u _{0T'+} '(\M '), F ^\flat \NN ')}
\ar[d] _-\sim
&
{u _{0*}  \R \mathcal{H} om _{\D ^\dag _{\X } (\hdag T _X) _{\Q}} (F ^\flat \M ', F ^\flat u _{0T'} ^{\prime !} \NN ')}
\ar[d] _-\sim
\\
{ \R \mathcal{H} om _{\D ^\dag _{\Y } (\hdag T ) _{\Q}} ( u _{0T+} F ^\flat (\M '), F ^\flat \NN ')}
\ar[r] _-\sim
&
{u _{0*}  \R \mathcal{H} om _{\D ^\dag _{\X } (\hdag T _X) _{\Q}} (F ^\flat \M ', u _{0T} ^! F ^\flat \NN ').}
}
\end{equation}
\end{theo}

\begin{proof}
Pour construire les isomorphismes horizontaux de \ref{frobadjdiag},
on procède, via \ref{isodualrelfrob}, comme pour \cite[1.2.9]{caro_courbe-nouveau}
ou \cite[IV.4.1]{Vir04}.
Il reste à vérifier la commutativité du diagramme \ref{frobadjdiag}
dont les isomorphismes verticaux sont induits fonctoriellement par $F ^\flat$
ou par les isomorphismes de commutation à Frobenius de l'image inverse extraordinaire
ou de l'image directe.
Afin de simplifier les notations (autrement, il s'agit de rajouter dans la preuve des $T$ ou $T _{X}$),
supposons que le diviseur $T$ est vide.
Vérifions d'abord la commutativité du diagramme canonique suivant :
  \begin{equation}\label{frobadjdiag1}
  \xymatrix @R=0,3cm @C=0,5cm  {
{ \R \mathcal{H} om _{\D ^\dag _{\Y ',\Q}} ( u _{0+}' (\M '), \NN ')}
\ar[r] _-\sim \ar[d] _-\sim ^-{F ^{\flat}}
&
{ \R \mathcal{H} om _{\D ^\dag _{\Y ',\Q}} ( u _{0+} '(\M '),
\NN ' \otimes _{\D ^\dag _{\Y ',\Q}} \D ^\dag _{\Y ',\Q} )}
\ar[d] _-\sim ^-{F ^{\flat}}
\\
{ \R \mathcal{H} om _{\D ^\dag _{\Y ,\Q}} ( F ^\flat u _{0+} '(\M '), F ^\flat \NN ')}
\ar@{=}[dd] \ar[r] _-\sim
&
{ \R \mathcal{H} om _{\D ^\dag _{\Y ,\Q}} (F ^\flat u _{0+} '(\M '),
F ^\flat(\NN ' \otimes _{\D ^\dag _{\Y ',\Q}} \D ^\dag _{\Y ',\Q} ))}
\\
&
{ \R \mathcal{H} om _{\D ^\dag _{\Y ,\Q}} (F ^\flat u _{0+} '(\M '),
F ^\flat (\NN ' )\otimes _{\D ^\dag _{\Y ,\Q}} F ^* F ^\flat \D ^\dag _{\Y ',\Q} )}
\ar[u] _-\sim
\\
{ \R \mathcal{H} om _{\D ^\dag _{\Y ,\Q}} ( F ^\flat u _{0+} '(\M '), F ^\flat \NN ')}
\ar[r] _-\sim\ar[d] _-\sim
&
{ \R \mathcal{H} om _{\D ^\dag _{\Y ,\Q}} (F ^\flat u _{0+} '(\M '),
F ^\flat (\NN ' )\otimes _{\D ^\dag _{\Y ,\Q}} \D ^\dag _{\Y ,\Q} )}
\ar[u] _-\sim \ar[d] _-\sim
\\
{ \R \mathcal{H} om _{\D ^\dag _{\Y ,\Q}} ( u _{0+} F ^\flat (\M '), F ^\flat \NN ')}
\ar[r] _-\sim
&
{ \R \mathcal{H} om _{\D ^\dag _{\Y ,\Q}} (u _{0+} F ^\flat  (\M '),
F ^\flat (\NN ' )\otimes _{\D ^\dag _{\Y ,\Q}} \D ^\dag _{\Y ,\Q} ).}
}
\end{equation}
\normalsize
Il découle de la commutativité du rectangle du haut de \cite[1.4.15.1]{caro_comparaison} (et par passage de gauche
à droite),
celle du rectangle du milieu de \ref{frobadjdiag1}. Enfin, celle des carrés s'établit par fonctorialité. 

En utilisant \cite[2.1.17.(iii)]{caro_comparaison}, on parvient au rectangle supérieur du diagramme ci-après :
  \begin{equation}\label{frobadjdiag2}
  \xymatrix @R=0,3cm {
{ \R \mathcal{H} om _{\D ^\dag _{\Y ',\Q}} ( u _{0+} '(\M '),
\NN ' \otimes _{\D ^\dag _{\Y ',\Q}} \D ^\dag _{\Y ',\Q} )}
\ar[d] _-\sim ^-{F ^{\flat}} \ar[r] _-\sim
&
{\NN ' \otimes _{\D ^\dag _{\Y ',\Q}}  ^\L
\R \mathcal{H} om _{\D ^\dag _{\Y' ,\Q}} ( u _{0+} '(\M '), \D ^\dag _{\Y ',\Q} )}
\\
{ \R \mathcal{H} om _{\D ^\dag _{\Y ,\Q}} (F ^\flat u _{0+} ' (\M '),
F ^\flat(\NN ' \otimes _{\D ^\dag _{\Y ',\Q}} \D ^\dag _{\Y ',\Q} ))}
&
{F ^\flat \NN ' \otimes _{\D ^\dag _{\Y ,\Q}}  ^\L F ^* \R \mathcal{H} om _{\D ^\dag _{\Y' ,\Q}}
( u _{0+} ' (\M '), \D ^\dag _{\Y ',\Q} )}
\ar[d] _-\sim ^-{F ^{\flat}}
\ar[u] _-\sim
\\
{ \R \mathcal{H} om _{\D ^\dag _{\Y ,\Q}} (F ^\flat u _{0+} ' (\M '),
F ^\flat (\NN ' )\otimes _{\D ^\dag _{\Y ,\Q}} F ^* F ^\flat \D ^\dag _{\Y' ,\Q} )}
\ar[u] _-\sim \ar[r] _-\sim
&
{F ^\flat \NN ' \otimes _{\D ^\dag _{\Y ,\Q}}  ^\L \R \mathcal{H} om _{\D ^\dag _{\Y ,\Q}}
( F ^{\flat} u _{0+} ' (\M '),  F ^{\flat} F ^*  \D ^\dag _{\Y ',\Q} )}
\\
{ \R \mathcal{H} om _{\D ^\dag _{\Y ,\Q}} (F ^\flat u _{0+} ' (\M '),
F ^\flat (\NN ' )\otimes _{\D ^\dag _{\Y ,\Q}} \D ^\dag _{\Y ,\Q} )}
\ar[u] _-\sim \ar[d] _-\sim \ar[r] _-\sim
&
{F ^\flat \NN ' \otimes _{\D ^\dag _{\Y ,\Q}}  ^\L \R \mathcal{H} om _{\D ^\dag _{\Y ,\Q}}
(  F ^{\flat} u _{0+} ' (\M '), \D ^\dag _{\Y ,\Q} )}
\ar[u] _-\sim\ar[d] _-\sim
\\
{ \R \mathcal{H} om _{\D ^\dag _{\Y ,\Q}} (u _{0+} F ^\flat  (\M '),
F ^\flat (\NN ' )\otimes _{\D ^\dag _{\Y ,\Q}} \D ^\dag _{\Y ,\Q} )}
\ar[r] _-\sim
&
{F ^\flat \NN ' \otimes _{\D ^\dag _{\Y ,\Q}}  ^\L \R \mathcal{H} om _{\D ^\dag _{\Y ,\Q}}
(u _{0+} F ^{\flat}  (\M '), \D ^\dag _{\Y ,\Q} ).}
}
\end{equation}
\normalsize
Les deux carrés l'étant par fonctorialité, on obtient la commutativité
de \ref{frobadjdiag2}.

En appliquant $F ^\flat \NN ' \otimes _{\D ^\dag _{\Y ,\Q}}  ^\L -$ à \ref{diagdualreldiv},
on obtient (modulo les passages de gauche à droite, i.e., modulo
$\omega \otimes _{\O}-$)
le rectangle commutatif (en bas) du diagramme suivant
\begin{equation}\label{frobadjdiag3}
  \xymatrix @R=0,3cm {
{\NN ' \otimes _{\D ^\dag _{\Y ',\Q}}  ^\L
\R \mathcal{H} om _{\D ^\dag _{\Y ',\Q}} ( u _{0+}' (\M '), \D ^\dag _{\Y ',\Q} )}
&
{\NN ' \otimes _{\D ^\dag _{\Y ',\Q}}  ^\L
u _{0+} ' \R \mathcal{H} om _{\D ^\dag _{\X ',\Q}} ( \M ', \D ^\dag _{\X ',\Q} )}
\ar[l] _-\sim ^-{\chi}
\\
{F ^\flat \NN ' \otimes _{\D ^\dag _{\Y ,\Q}}  ^\L F ^* \R \mathcal{H} om _{\D ^\dag _{\Y ',\Q}}
( u _{0+} '(\M '), \D ^\dag _{\Y' ,\Q} )}
\ar[d] _-\sim ^-{F ^{\flat}}
\ar[u] _-\sim
&
{F ^\flat \NN ' \otimes _{\D ^\dag _{\Y ,\Q}}  ^\L
F^* u _{0+} ' \R \mathcal{H} om _{\D ^\dag _{\X ',\Q}} ( \M ', \D ^\dag _{\X ',\Q} )}
\ar[l] _-\sim^-{\chi} \ar[u] _-\sim \ar[d] _-\sim
\\
{F ^\flat \NN ' \otimes _{\D ^\dag _{\Y ,\Q}}  ^\L \R \mathcal{H} om _{\D ^\dag _{\Y ,\Q}}
( F ^{\flat} u _{0+} '(\M '),  F ^{\flat} F ^*  \D ^\dag _{\Y ',\Q} )}
&
{F ^\flat \NN ' \otimes _{\D ^\dag _{\Y ,\Q}}  ^\L
 u _{0+}  F ^* \R \mathcal{H} om _{\D ^\dag _{\X ',\Q}} ( \M ', \D ^\dag _{\X ',\Q} )}
\ar[d] _-\sim ^-{F ^{\flat}}
\\
{F ^\flat \NN ' \otimes _{\D ^\dag _{\Y ,\Q}}  ^\L \R \mathcal{H} om _{\D ^\dag _{\Y ,\Q}}
(  F ^{\flat} u _{0+} '(\M '), \D ^\dag _{\Y ,\Q} )}
\ar[u] _-\sim\ar[d] _-\sim
&
{F ^\flat \NN ' \otimes _{\D ^\dag _{\Y ,\Q}}  ^\L
 u _{0+} \R \mathcal{H} om _{\D ^\dag _{\X ,\Q}} (F ^\flat \M ',F ^\flat F ^* \D ^\dag _{\X ',\Q} )}
\\
{F ^\flat \NN ' \otimes _{\D ^\dag _{\Y ,\Q}}  ^\L \R \mathcal{H} om _{\D ^\dag _{\Y ,\Q}}
(u _{0+} F ^{\flat}  (\M '), \D ^\dag _{\Y ,\Q} )}
&
{F ^\flat \NN ' \otimes _{\D ^\dag _{\Y ,\Q}}  ^\L
 u _{0+} \R \mathcal{H} om _{\D ^\dag _{\X ,\Q}} (F ^\flat \M ', \D ^\dag _{\X ,\Q} ).}
\ar[u] _-\sim
\ar[l] _-\sim ^-{\chi}
}
\end{equation}
\normalsize
La commutativité du carré (supérieur) se vérifie par fonctorialité. D'où celle de \ref{frobadjdiag3}.

En posant $\E ':= \R \mathcal{H} om _{\D ^\dag _{\X ',\Q}} ( \M ', \D ^\dag _{\X ',\Q} )$,
$\E := F ^* \E '$, $\NN := F ^\flat \NN '$,
on dispose du diagramme commutatif
suivant
\begin{equation}\label{frobadjdiag4}
  \xymatrix @R=0,3cm @C=0,5cm {
{\NN ' \otimes _{\D ^\dag _{\Y ',\Q}}  ^\L
u _{0+} '( \R \mathcal{H} om _{\D ^\dag _{\X ',\Q}} ( \M ', \D ^\dag _{\X ',\Q} ) ) }
\ar[r] _-\sim ^-{\mathrm{proj}}
\ar@{=}[d]
&
{ u _{0*} (u _{0} ^{\prime !} \NN ' \otimes ^\L _{\D ^\dag _{\X ',\Q}}
\R \mathcal{H} om _{\D ^\dag _{\X ',\Q}} ( \M ', \D ^\dag _{\X ',\Q} ) )}
\ar@{=}[d]
\\
{\NN ' \otimes _{\D ^\dag _{\Y ',\Q}}  ^\L
u _{0*} ( \D ^\dag _{\Y '\leftarrow \X ',\Q} \otimes _{\D ^\dag _{\X ',\Q}} \E ' ) }
\ar[r] _-\sim ^-{\mathrm{proj}}
&
{ u _{0*} (u _{0} ^{\text{-}1} \NN ' \otimes _{u _{0} ^{\text{-}1} \D ^\dag _{\Y ',\Q}}  ^\L
\D ^\dag _{\Y '\leftarrow \X ',\Q} \otimes ^\L _{\D ^\dag _{\X ',\Q}} \E ' ) }
\\
{\NN  \otimes _{\D ^\dag _{\Y ,\Q}}  ^\L
F^* u _{0*} ( \D ^\dag _{\Y ' \leftarrow \X ',\Q} \otimes _{\D ^\dag _{\X ',\Q}} \E ' )}
\ar[u] _-\sim \ar[d] _-\sim ^-{\mathrm{proj}}
\ar[r] _-\sim ^-{\mathrm{proj}}
&
{ u _{0*} (u _{0} ^{\text{-}1} (\NN
\underset{\D ^\dag _{\Y \Q}}{\overset{\L}{\otimes}}
F^* \D ^\dag _{\Y ',\Q})
\underset{u _{0} ^{\text{-}1} \D ^\dag _{\Y '\Q}}{\overset{\L}{\otimes}}
\D ^\dag _{\Y '\leftarrow \X ',\Q} \underset{\D ^\dag _{\X ',\Q}}{\overset{\L}{\otimes}} \E ' )}
\ar[u] _-\sim \ar[d] _-\sim
\\
{\NN
\! \! \! \!
\underset{\D ^\dag _{\Y ,\Q}}{\overset{\L}{\otimes}}
\! \! \! \!
u _{0*} (u _{0} ^{\text{-}1} F^* \D ^\dag _{\Y ',\Q} \!
\underset{u _{0} ^{\text{-}1} \D ^\dag _{\Y '\Q}}{\overset{\L}{\otimes}}
\!\! \! \D ^\dag _{\Y '\leftarrow \X ',\Q}
\underset{\D ^\dag _{\X ',\Q}}{\overset{\L}{\otimes}}
\! \! \! \! \E ' ) }
\ar[r] _-\sim ^-{\mathrm{proj}}
&
{u _{0*} (u _{0} ^{\text{-}1} \NN \! \! \! \!
\underset{u _{0} ^{\text{-}1} \D ^\dag _{\Y \Q}}{\overset{\L}{\otimes}}
\! \! \! \!
u _{0} ^{\text{-}1}  F ^* \D ^\dag _{\Y ',\Q} \!
\underset{u _{0} ^{\text{-}1} \D ^\dag _{\Y '\Q}}{\overset{\L}{\otimes}}
\! \! \! \!
\D ^\dag _{\Y '\leftarrow \X ',\Q} \underset{\D ^\dag _{\X ',\Q}}{\overset{\L}{\otimes}} \! \! \! \E ' ) }
\\
{ \NN \! \! \! \!
\underset{\D ^\dag _{\Y ,\Q}}{\overset{\L}{\otimes}}
\! \! \! \!
u _{0*} (u _{0} ^{\text{-}1} F^* \D ^\dag _{\Y ',\Q}\!
\underset{u _{0} ^{\text{-}1} \D ^\dag _{\Y '\Q}}{\overset{\L}{\otimes}}
                \!\! \! \! \!
F ^\flat \D ^\dag _{\Y '\leftarrow \X ',\Q}
                        \!\! \!
\underset{\D ^\dag _{\X ,\Q}}{\overset{\L}{\otimes}}
\!\! \E  ) }
\ar[r] _-\sim ^-{\mathrm{proj}} \ar[u] _-\sim \ar[d] _-\sim
&
{u _{0*} (u _{0} ^{\text{-}1} \NN
\! \! \! \!
\underset{u _{0} ^{\text{-}1} \D ^\dag _{\Y \Q}}{\overset{\L}{\otimes}}
\! \! \! \! \!
u _{0} ^{\text{-}1} F ^* \D ^\dag _{\Y ',\Q}
\! \!
\underset{u _{0} ^{\text{-}1} \D ^\dag _{\Y '\Q}}{\overset{\L}{\otimes}}
\!\! \! \! \!
F ^\flat \D ^\dag _{\Y '\leftarrow \X' ,\Q} \underset{\D ^\dag _{\X ,\Q}}{\overset{\L}{\otimes}} \!\! \E  ) }
\ar[u] _-\sim \ar[d] _-\sim
\\
{ \NN  \otimes _{\D ^\dag _{\Y ,\Q}}  ^\L  u _{0*} (
\D ^\dag _{\Y \leftarrow \X ,\Q} \otimes _{\D ^\dag _{\X ,\Q}} \E  ) }
\ar[r] _-\sim ^-{\mathrm{proj}} \ar@{=}[d]
&
{u _{0*} (u _{0} ^{\text{-}1}\NN \underset{u _{0} ^{\text{-}1} \D ^\dag _{\Y \Q}}{\overset{\L}{\otimes}}
\D ^\dag _{\Y \leftarrow \X ,\Q} \otimes _{\D ^\dag _{\X ,\Q}}  ^\L \E  ) }
\ar@{=}[d]
\\
{\NN  \otimes _{\D ^\dag _{\Y ,\Q}}  ^\L  u _{0+}
 F ^* \R \mathcal{H} om _{\D ^\dag _{\X ',\Q}} ( \M ', \D ^\dag _{\X ',\Q} )  }
\ar[r] _-\sim ^-{\mathrm{proj}}
&
{u _{0*} (u _{0} ^{!} ( \NN  ) \otimes _{\D ^\dag _{\X ,\Q}} ^\L
F ^* \R \mathcal{H} om _{\D ^\dag _{\X ',\Q}} ( \M ', \D ^\dag _{\X ',\Q} ) ) }
\\
{\NN  \otimes _{\D ^\dag _{\Y ,\Q}}  ^\L  u _{0+}
\R \mathcal{H} om _{\D ^\dag _{\X ,\Q}} (F ^\flat  \M ', \D ^\dag _{\X ,\Q} )}
\ar[r] _-\sim ^-{\mathrm{proj}}
\ar[u] _-\sim
&
{u _{0*} (u _{0} ^{!} (\NN  ) \otimes _{\D ^\dag _{\X ,\Q}} ^\L
\R \mathcal{H} om _{\D ^\dag _{\X ,\Q}} (F ^\flat \M ', \D ^\dag _{\X ,\Q} ) ) .}
\ar[u] _-\sim
}
\end{equation}
\normalsize
En effet, le troisième carré du haut est commutatif par transitivité des morphismes de
projection tandis que les autres le sont par définition de $\E'$, $\E$, $\NN$ et/ou par fonctorialité.

En appliquant le foncteur $F ^\flat \NN ' \otimes _{\D ^\dag _{\Y ,\Q}}  ^\L -$
à l'isomorphisme de commutation à Frobenius de l'image directe :
$F^* u _{0+} ' \R \mathcal{H} om _{\D ^\dag _{\X ',\Q}} ( \M ', \D ^\dag _{\X ',\Q} )
\riso
 u _{0+}  F ^* \R \mathcal{H} om _{\D ^\dag _{\X ',\Q}} ( \M ', \D ^\dag _{\X ',\Q} )$, on obtient 
le composé de gauche du troisième terme du haut vers le troisième terme du bas de \ref{frobadjdiag4}
ainsi que la deuxième flèche de droite en partant du haut de \ref{frobadjdiag3}.
On en déduit que le composé de droite de \ref{frobadjdiag3} est égal à celui de gauche de \ref{frobadjdiag4}.

De manière analogue au diagramme composé de \ref{frobadjdiag1} et \ref{frobadjdiag2} (sauf leur dernière ligne), on construit le carré du haut sans $u _{0*}$ du diagramme ci-dessous:
\begin{equation}\label{frobadjdiag5}
  \xymatrix @R=0,3cm @C=0,5cm {
{ u _{0*} (u _{0} ^{\prime !} \NN ' \otimes _{\D ^\dag _{\X ',\Q}} ^\L
\R \mathcal{H} om _{\D ^\dag _{\X ',\Q}} ( \M ', \D ^\dag _{\X ',\Q} ) ) }
\ar[r] _-\sim \ar[d] _-\sim
&
{ u _{0*} (\R \mathcal{H} om _{\D ^\dag _{\X ',\Q}} ( \M ', u _{0} ^{\prime !} \NN ' ) ) }
\ar[d] _-\sim ^{F ^\flat}
\\
{ u _{0*} (F ^\flat u _{0} ^{\prime !} \NN ' \otimes _{\D ^\dag _{\X ,\Q}} ^\L
\R \mathcal{H} om _{\D ^\dag _{\X ,\Q}} ( F ^\flat \M ', \D ^\dag _{\X ,\Q} ) ) }
\ar[r] _-\sim \ar[d] _-\sim
&
{ u _{0*} (\R \mathcal{H} om _{\D ^\dag _{\X ,\Q}} ( F ^\flat  \M ', F ^\flat  u _{0} ^{\prime !} \NN ' ) ) }
\ar[d] _-\sim
\\
{ u _{0*} ( u _{0} ^{!} F ^\flat \NN ' \otimes _{\D ^\dag _{\X ,\Q}} ^\L
\R \mathcal{H} om _{\D ^\dag _{\X ,\Q}} ( F ^\flat \M ', \D ^\dag _{\X ,\Q} ) ) }
\ar[r] _-\sim
&
{ u _{0*} (\R \mathcal{H} om _{\D ^\dag _{\X ,\Q}} ( F ^\flat  \M ',  u _{0} ^{!} F ^\flat  \NN ' ) ) .}
}
\end{equation}
Par fonctorialité, on en déduit la commutativité de \ref{frobadjdiag5}.
On vérifie par définition et fonctorialité
que le composé de gauche de \ref{frobadjdiag5} correspond à celui de droite de
\ref{frobadjdiag4}.
En composant \ref{frobadjdiag1}, \ref{frobadjdiag2}, \ref{frobadjdiag3}, \ref{frobadjdiag4}
et \ref{frobadjdiag5}, on obtient le diagramme \ref{frobadj} et donc sa commutativité.
 \end{proof}

\begin{coro}\label{frobadjdiagcor}
    Soient $\M ' \in D ^{\mathrm{b}} _\mathrm{coh}(\D ^\dag _{\X '} (\hdag T _{X'}) _{\Q} \overset{^{\mathrm{d}}}{})$
et
$\NN ' \in D ^{\mathrm{b}} _\mathrm{coh}(\D ^\dag _{\Y '} (\hdag T' ) _{\Q} \overset{^{\mathrm{d}}}{})$.
On bénéficie de {\it l'isomorphisme d'adjonction}
$\mathrm{adj} \  :\ \mathrm{Hom} _{\D ^\dag _{\Y '} (\hdag T ') _{\Q}} ( u _{0T'+} ' (\M '), \NN ')
\riso
{\mathrm{Hom} _{\D ^\dag _{\X '} (\hdag T _{X'}) _{\Q}} ( \M ', u _{0T'} ^{\prime !} \NN ')}$
qui s'inscrit dans le diagramme commutatif suivant :
\begin{equation}\label{frobadjdiagcordiag}
  \xymatrix @R=0,3cm {
{ \mathrm{Hom} _{\D ^\dag _{\Y '} (\hdag T ') _{\Q}} ( u _{0T'+} ' (\M '), \NN ')}
\ar[r] _-\sim ^{\mathrm{adj}} \ar[d] _-\sim ^-{F ^{\flat}}
&
{\mathrm{Hom} _{\D ^\dag _{\X '} (\hdag T _{X'}) _{\Q}} ( \M ', u _{0T'} ^{\prime !} \NN ')}
\ar[d] _-\sim ^-{F ^{\flat}}
\\
{ \mathrm{Hom} _{\D ^\dag _{\Y } (\hdag T ) _{\Q}} ( F ^\flat u _{0T'+} '(\M '), F ^\flat \NN ')}
\ar[d] _-\sim
&
{\mathrm{Hom} _{\D ^\dag _{\X } (\hdag T _X) _{\Q}} (F ^\flat \M ', F ^\flat u _{0T'} ^{\prime !} \NN ')}
\ar[d] _-\sim
\\
{ \mathrm{Hom} _{\D ^\dag _{\Y } (\hdag T ) _{\Q}} ( u _{0T+} F ^\flat (\M '), F ^\flat \NN ')}
\ar[r] _-\sim ^{\mathrm{adj}}
&
{\mathrm{Hom} _{\D ^\dag _{\X } (\hdag T _X) _{\Q}} (F ^\flat  \M ', u _{0T} ^! F ^\flat \NN ').}
}
\end{equation}
\end{coro}
\begin{proof}
  On applique le foncteur $\mathcal{H} ^{0}\circ \R\Gamma (\Y,-)$ à \ref{frobadjdiag}.
 \end{proof}

\begin{rema}\label{remacompisofrobu+}
L'isomorphisme canonique de commutation
à Frobenius de l'image directe a été construit de deux façons différentes.

$(a)$ La première, due à Berthelot (voir \cite[3.4.4]{Be2} pour la version \og à gauche\fg),
se construit directement. C'est cet isomorphisme de commutation à Frobenius de l'image directe
que nous avons utilisé dans ce manuscript.

$(b)$ Une deuxième construction existe lorsque l'on se restreint aux morphismes {\it propres}
(voir \cite[1.2.12]{caro_courbe-nouveau}).
Elle consiste à utiliser
(l'isomorphisme canonique de commutation à Frobenius de l'image inverse extraordinaire et)
l'isomorphisme d'adjonction
entre l'image directe et l'image inverse extraordinaire
(cet isomorphisme d'adjonction n'est valable que pour les morphismes {\it propres})
de telle façon que le diagramme
\ref{frobadjdiagcordiag} soit commutatif.
Cet isomorphisme de commutation à Frobenius de l'image directe donne par construction
la compatibilité à Frobenius des morphismes d'adjonction.

Il résulte alors du corollaire \ref{frobadjdiagcor} que ces deux constructions coïncident dans le cas d'une {\it immersion fermée}.
Nous obtenons ainsi une première {\og unification\fg}. L'extension raisonnable de \ref{frobadjdiagcor} au cas
des morphismes {\it propres} (i.e., l'isomorphisme d'adjonction de \ref{isodualrelpropreadj} est compatible
à Frobenius) garantirait une unification complète des deux constructions.
Remarquons enfin que la difficulté apparaissant dans la preuve de l'extension de \ref{frobadjdiagcor} au cas
des morphismes {\it propres} est la construction et la compatibilité à Frobenius du morphisme trace.

\end{rema}

\subsection{Construction de l'isomorphisme de dualité relative dans le cas d'un morphisme projectif non relevable}

\begin{vide}\label{rematrduarel}
Soient $f$ : $\X \rightarrow \Y $ un morphisme propre de $\V$-schémas formels lisses, $T$ un diviseur
de $Y$ tel que $T _X:= f ^{-1} (T)$ soit un diviseur de $X$.
D'après \cite[1.2.5,6]{caro_courbe-nouveau}, on dispose du morphisme trace 
$\mathrm{Tr}  _{T+} \ : \
u _{T+} (\omega _{\X} (\hdag T _{X}) _{\Q} [d _X]) \rightarrow \omega _{\Y} (\hdag T) _{\Q} [d _Y]$.
Celui-ci est obtenu par extension des scalaires via $\D ^\dag _{\Y,\Q} \to \D ^\dag _{\Y } (\hdag T ) _{\Q}$ à partir du morphisme trace construit par Virrion qui vérifie la condition de transitivité usuelle. Par construction, la transitivité du morphisme de Virrion implique alors celle de $\mathrm{Tr}  _{T+} $. De manière analogue à \cite[3.5.7]{Beintro2}, 
il en dérive que l'isomorphisme de dualité relative
$\chi$ : $f _{T +} \circ \DD _{\X, T _X} ( \M) \riso \DD _{\Y, T }  \circ f _{ T +}  ( \M)$ de \cite[1.2.7]{caro_courbe-nouveau} fonctoriel
en $ \M\in D ^{\mathrm{b}} _\mathrm{coh}(\D ^\dag _{\X } (\hdag T _X) _{\Q} \overset{^{\mathrm{d}}}{})$ vérifie la condition usuelle de transitivité.

\end{vide}

\begin{defi}
\label{quasi-relevables}
Soit $f _0$ : $X  \rightarrow Y $ un morphisme de $k$-schémas lisses. 
On dit que $f _{0}$ est un morphisme {\it propre quasi-relevable} s'il existe un $\V$-schéma formel lisse $\X$ relevant $X$ et si $f _{0}=p _{0} \circ u _{0}$, où $u _{0}\,:\, X\hookrightarrow Z$ est une immersion fermée et $p _0\,:\, Z \to Y$ est un morphisme qui se relève en un morphisme propre $p\,: \, \ZZ \to \Y$ de $\V$-schémas formels lisses. 

Par exemple, $f _0$ est un morphisme propre quasi-relevable s'il existe $\X$, $\Y$ deux $\V$-schémas formels lisses relevants respectivement $X$, $Y$ et si 
\begin{itemize}
\item soit $f _{0}$ est un morphisme projectif, 
\item ou soit il existe une immersion $X\hookrightarrow \PP$ telle que $\PP$ soit un $\V$-schéma formel propre et lisse (e.g. $X $ est quasi-projectif).
\end{itemize}

\end{defi}

\begin{theo}\label{isodualrelpropreadj}
  Soient $\X$ et $\Y$ deux $\V$-schémas formels lisses, $f _0$ : $X  \rightarrow Y $
un morphisme propre quasi-relevable entre leur fibre spéciale (e.g. $f _{0}$ est projectif), $T$ un diviseur de $Y $ tel que
$T _{X}:= f _0 ^{-1}  (T)$ soit un diviseur de $X $.

\begin{enumerate}

\item Il existe alors un isomorphisme canonique fonctoriel en
$ \M\in D ^{\mathrm{b}} _\mathrm{coh}(\D ^\dag _{\X } (\hdag T _X) _{\Q} \overset{^{\mathrm{d}}}{})$,
\begin{equation}
\label{isodualrelpropreadj1} 
\chi\ : \ f _{0T +} \circ \DD _{\X, T _X} ( \M) \riso \DD _{\Y, T }  \circ f _{0 T +}  ( \M),
\end{equation}

dit {\og isomorphisme de dualité relative \fg}. 

\item On bénéficie de plus, pour tout
$\M  \in D ^{\mathrm{b}} _\mathrm{coh}(\D ^\dag _{\X } (\hdag T _{X}) _{\Q} \overset{^{\mathrm{d}}}{})$
et
$\NN  \in D ^{\mathrm{b}} _\mathrm{coh}(\D ^\dag _{\Y } (\hdag T ) _{\Q} \overset{^{\mathrm{d}}}{})$,
 de l'isomorphisme canonique d'adjonction
\begin{equation}
\label{isodualrelpropreadj2} 
\mathrm{adj} \  :\ \mathrm{Hom} _{\D ^\dag _{\Y } (\hdag T ) _{\Q}} ( f _{0T+}  (\M ), \NN )
\riso
{\mathrm{Hom} _{\D ^\dag _{\X } (\hdag T _{X}) _{\Q}} ( \M , f _{0T} ^{!} \NN )}.
\end{equation}

\end{enumerate}

\end{theo}

\begin{proof}

1) a). En choisissant un morphisme propre et lisse $p$ : $\ZZ \rightarrow \Y$
de $\V$-schémas formels lisses, une immersion fermée $u _0 $ : $X  \hookrightarrow Z $ telle que $f _0= p _0 \circ u _0$,
l'isomorphisme de dualité relative à $f _{0}$ est par définition le morphisme rendant commutatif le diagramme ci-dessous :
\begin{equation}
  \label{chipropdef}
\xymatrix @R=0,3cm {
{f _{0T +} \circ \DD _{\X, T _X} ( \M)}
\ar[d] _-\sim
\ar@{.>}[rr] ^-{\chi}
&&
{\DD _{\X, T  }\circ f _{0T +} ( \M)}
\ar[d] _-\sim
\\
{p _{ T+} u _{0T _{Z} +} \DD _{\X, T _X} ( \M)}
\ar[r] _-\sim ^-{\chi}
&
{p _{ T+}  \DD _{\ZZ, T _Z }   u _{0T _{Z} +}  ( \M)}
\ar[r] _-\sim ^-{\chi}
&
{\DD _{\X, T  } p _{0 T+}    u _{0T _{Z} +}  ( \M),}
 }
\end{equation}
où $T _{Z } := p _0 ^{-1} (T )$ et où les isomorphismes de dualité relative du bas sont ceux provenant soit
du cas relevable (\cite[1.2.7]{caro_courbe-nouveau})
soit de celui d'une immersion fermée (\ref{isodualrelfrob}).

1) b). Il reste à présent à vérifier que cet isomorphisme est canonique, i.e., qu'il ne dépend pas du choix de la factorisation
de $f _0$.
Soient
un morphisme propre et lisse $p$ : $\ZZ \rightarrow \Y$ (resp. $p'$ : $\ZZ' \rightarrow \Y$) de $\V$-schémas formels lisses, une immersion fermée $u _0 $ : $X  \hookrightarrow Z $ telle que $f _0= p _0 \circ u _0$
(resp. $u ' _0 $ : $X \hookrightarrow Z '$ telle que $f _0= p '_0 \circ u '_0$).
Le morphisme $u _0$ (resp. $u ' _0$)
se factorise en l'immersion fermée $u '' _0 :=u _0 \times u '_0$ : $X \hookrightarrow Z \times _Y Z'$
(\cite[5.4.5]{EGAI})
suivie de la projection $Z \times _Y Z'\rightarrow  Z$ (resp. $Z \times _Y Z'\rightarrow  Z'$)
qui se relève en la projection
$q $ : $\ZZ \times _{\Y}  \ZZ' \rightarrow \ZZ$
(resp. $q'$ : $\ZZ \times _{\Y}  \ZZ' \rightarrow \ZZ '$).
Comme $ p \circ q = p ' \circ q'$, grâce à la transitivité du cas relevable (voir \ref{rematrduarel}),
on se ramène au cas où
$f _0$ est une immersion fermée (plus précisément aux deux cas où $f _0$ est
égal à $u _0= q _0 \circ u '' _0$ ou à $u '_0=q ' _0 \circ u '' _0$).
Supposons donc que $f _{0}$ est une immersion fermée et 
vérifions que le diagramme \ref{chipropdef} où l'isomorphisme du haut est l'isomorphisme de dualité relative de \ref{isodualrelfrob} est commutatif. Comme on dispose de l'isomorphisme de bidualité (voir \cite{virrion}), il suffit d'établir que le dual de \ref{chipropdef} est commutatif. 
Comme $\M$ est isomorphe à un complexe dont les termes sont des $\D ^\dag _{\X } (\hdag T _{X}) _{\Q}$-modules à droite cohérents, comme les $\D ^\dag _{\X } (\hdag T _{X}) _{\Q}$-modules à droite cohérents (resp. les $\D ^\dag _{\ZZ } (\hdag T _{Z}) _{\Q}$-modules à droite cohérents dans l'image essentielle de $ u_{0T+}$, $ u _{0T!}$)
sont acycliques pour les foncteurs $ f _{0T+}$, $ f _{0T!}$, $ u_{0T+}$, $ u _{0T!}$ (resp. $p _{T+}$, $p _{T!}$), on se ramène au cas où $\M$ est un $\D ^\dag _{\X } (\hdag T _{X}) _{\Q}$-module à droite cohérent.
La commutativité du dual de \ref{chipropdef} est alors local en $\Y$.
On peut donc supposer $\Y$ affine. Comme $\X$ est alors affine, le morphisme $u _0$ se relève.
On conclut grâce à la transitivité de l'isomorphisme de dualité relative dans ce cas (voir \ref{rematrduarel}).

2). Enfin, la construction de l'isomorphisme d'adjonction de \ref{isodualrelpropreadj2} à partir de l'isomorphisme de dualité relative de \ref{isodualrelpropreadj1} est identique à celle de \cite[1.2.10]{caro_courbe-nouveau} (voir aussi la construction dans la preuve de \ref{frobadj}).

 \end{proof}

\bibliographystyle{smfalpha}
\bibliography{Bibliotheque}

\bigskip
\noindent Daniel Caro\\
Laboratoire de Mathématiques Nicolas Oresme\\
Université de Caen
Campus 2\\
14032 Caen Cedex\\
France.\\
email: daniel.caro@math.unicaen.fr

\end{document}